\documentclass[a4paper, reqno]{amsbook}

\usepackage[ngerman,american]{babel}  
\usepackage{amscd}    
\usepackage{amssymb}
\usepackage{color}
\usepackage{dsfont}   
\usepackage{graphicx}
\usepackage{ifpdf}
\usepackage{listings}
\usepackage{float}
\ifpdf
  \usepackage[pdftex,plainpages=false,pdfpagelabels]{hyperref}
  \usepackage{thumbpdf}
\else
  \usepackage{hyperref}
\fi

\input{edge.sty}

\font\cyr=wncyr10 
\newcommand{\lobachevski}{\operatorname{\mbox{\cyr L}}}

\numberwithin{equation}{chapter}

\floatstyle{ruled}
\newfloat{listing}{tbp}{lol}
\floatname{listing}{Listing}

\theoremstyle{plain}

\newtheorem{theorem}{Theorem}
\numberwithin{theorem}{chapter}

\newtheorem{lemma}{Lemma}
\numberwithin{lemma}{chapter}

\newtheorem*{corollary}{Corollary}

\newtheorem*{proposition}{Proposition}

\newtheorem*{flowThm}{Feasible Flow Theorem}

\theoremstyle{definition}

\newtheorem*{definition}{Definition}

\newtheorem*{quasidefinition}{`Definition'}

\theoremstyle{remark}

\newtheorem*{remark}{Remark}

\newcommand{\artanh}{\operatorname{artanh}}

\newcommand{\C}{\mathds{C}}
\newcommand{\Chat}{\widehat{\mathds{C}}}
\newcommand{\R}{\mathds{R}}
\newcommand{\Z}{\mathds{Z}}
\newcommand{\PGL}{\operatorname{\mathit{PGL}}}
\newcommand{\LHS}{\text{LHS}}
\newcommand{\RHS}{\text{RHS}}
\newcommand{\im}{\operatorname{Im}}
\newcommand{\kernel}{\operatorname{Ker}}
\newcommand{\image}{\operatorname{Im}}
\newcommand{\Li}[1][2]{\operatorname{Li}_{#1}}
\newcommand{\Cl}[1][2]{\operatorname{Cl}_{#1}}

\newcommand{\vecE}{\vec{E}}
\newcommand{\vece}{\vec{e}}
\newcommand{\Eint}{E_{\text{\it{}int}}}
\newcommand{\vecEint}{\vecE_{\text{\it{}int}}}

\newcommand{\Seuc}{S_{\text{euc}}}
\newcommand{\SeucL}{S_{\text{euc}}^{\text{\it L}}}
\newcommand{\SeucH}{S_{\text{euc}}^{\text{\it H}}}
\newcommand{\Shyp}{S_{\text{hyp}}}
\newcommand{\ShypL}{S_{\text{hyp}}^{\text{\it L}}}
\newcommand{\ShypH}{S_{\text{hyp}}^{\text{\it H}}}
\newcommand{\Ssph}{S_{\text{sph}}}
\newcommand{\SsphTilde}{\widetilde{S}_{\text{sph}}}
\newcommand{\SsphL}{S_{\text{sph}}^{\text{\it L}}}
\newcommand{\SsphH}{S_{\text{sph}}^{\text{\it H}}}
\newcommand{\Shat}{\widehat{S}}
\newcommand{\Sproto}{S_{\text{proto}}}
\newcommand{\ShatProto}{\Shat_{\text{proto}}}

\DeclareMathOperator{\into}{\mathit{in}}
\DeclareMathOperator{\outof}{\mathit{ex}}
\DeclareMathOperator{\id}{id}

\newcommand{\xor}{\;\mbox{\raisebox{-0.7ex}{$\widehat~$}}\;}

\newcommand{\Sym}{\operatorname{Sym}}
\newcommand{\Alt}{\operatorname{Alt}}

\input{listings_config.sty}
\begin{document}
\ifx\href\undefined\else\hypersetup{linktocpage=true}\fi

\thispagestyle{empty}
\renewcommand{\thepage}{}
  
\begin{center}
  \vspace*{\fill}
  
  {\Large \textbf{Variational Principles for Circle Patterns}}\\
    
  \vspace{\fill}

  vorgelegt von\\
  Dipl.-Math.\ Boris A.~Springborn\\[\baselineskip]  
  von der Fakult\"at II -- Mathematik und Naturwissenschaften\\
  der Technischen Universit\"at Berlin\\
  zur Erlangung des akademischen Grades\\[\baselineskip]
  Doktor der Naturwissenschaften\\
  -- Dr.\ rer.\  nat.\ --\\[\baselineskip]
  genehmigte Dissertation

    

\vspace{\fill}

\begin{tabular}{@{}rl@{}}
\multicolumn{2}{@{}c@{}}{Promotionsausschuss}\\
  \\
  Vorsitzender: & Prof.\ Dr.\ Michael E. Pohst\\
  Gutachter/Berichter:    & Prof.\ Dr.\ Alexander I. Bobenko\\
                          & Prof.\ Dr.\ G\"unter M.\ Ziegler
\end{tabular}
    
\vspace{\baselineskip}
Tag der wissenschaftlichen Aussprache: 27.\ November 2003\\[10ex]
    
\vfill
    
    Berlin 2003\\
    D 83             
  \end{center}


\frontmatter

{\small

\centerline{\Large\bf Abstract}

\vspace{\baselineskip{}}

A Delaunay cell decomposition of a surface with constant curvature gives rise
to a circle pattern, consisting of the circles which are circumscribed to the
facets. We treat the problem whether there exists a Delaunay cell
decomposition for a given (topological) cell decomposition and given
intersection angles of the circles, whether it is unique and how it may be
constructed. Somewhat more generally, we allow cone-like singularities in the
centers and intersection points of the circles. We prove existence and
uniqueness theorems for the solution of the circle pattern problem using a
variational principle. The functionals (one for the euclidean, one for the
hyperbolic case) are convex functions of the radii of the circles. The
critical points correspond to solutions of the circle pattern problem. The
analogous functional for the spherical case is not convex, hence this case is
treated by stereographic projection to the plane. From the existence and
uniqueness of circle patterns in the sphere, we derive a strengthened version
of Steinitz' theorem on the geometric realizability of abstract polyhedra.

We derive the variational principles of Colin~de~Verdi\`ere, Br\"agger, and
Rivin for circle packings and circle patterns from our variational
principles. In the case of Br\"agger's and Rivin's functionals, this requires
a Legendre transformation of our euclidean functional. The respective
Legendre transformations of the hyperbolic and spherical functionals lead to
new variational principles. The variables of the transformed functionals are
certain angles instead of radii. The transformed functionals may be
interpreted geometrically as volumes of certain three-dimensional polyhedra
in hyperbolic space. Leibon's functional for hyperbolic circle patterns
cannot be derived from our functionals. But we construct yet another
functional from which both Leibon's and our functionals can be derived. By
applying the inverse Legendre transformation to Leibon's functional, we
obtain a new variational principle for hyperbolic circle patterns.

We present Java software to compute and visualize circle patterns.

\selectlanguage{ngerman} 

\vspace{\baselineskip}
\centerline{\Large\bf Zusammenfassung}

\vspace{\baselineskip}

Eine Delaunay-Zellzerlegung einer Fl\"ache konstanter Kr\"ummung liefert ein
Kreismuster, welches aus den Kreisen besteht, die den Facetten umschrieben
sind. Wir betrachten das Problem, ob es f\"ur eine vorgegebene (topologische)
Zellzerlegung und vorgegebene Schnittwinkel zwischen den Kreisen eine
entsprechende Delaunay-Zellzerlegung gibt, ob sie eindeutig ist, und wie sie
zu konstruieren ist.  Etwas allgemeiner lassen wir auch kegelartige
Singularit\"aten in den Mittel- und Schnittpunkten der Kreise zu. Wir
beweisen Existenz- und Eindeutigkeitss\"atze f\"ur die L\"osung des
Kreismusterproblems mit Hilfe von Variationsprinzipien. Die Funktionale (eins
f\"ur den euklidischen, eins f\"ur den hyperbolischen Fall) sind konvexe
Funktionen der Radien der Kreise. Kritische Punkte entsprechen L\"osungen
des Kreismusterproblems. Das analoge Funktional f\"ur den sph\"arischen Fall
ist nicht konvex, deshalb wird dieser Fall durch stereographische Projektion
in die Ebene erledigt. Aus der Existenz und Eindeutigkeit von Kreismustern in
der Sph\"are folgern wir eine versch\"arfte Version des Satzes von Steinitz
\"uber die geometrische Realisierbarkeit von abstrakten Polyedern.

Wir leiten die Variationsprinzipien von Colin~de~Verdi\`ere, Br\"agger und
Rivin f\"ur Kreispackungen bzw.\ Kreismuster aus unseren Variationsprinzipien
ab. Im Fall der Funktionale von Br\"agger und Rivin erfordert dies eine
Legendretransformation unseres euklidischen Funktionals. Entsprechende
Legendretransformationen des hyperbolischen und des sph\"arischen Funktionals
liefern neue Variationsprinzipien. Die Variablen der transformierten
Funktionale sind nicht Radien, sondern bestimmte Winkel. Die transformierten
Funktionale besitzen eine geometrische Interpretation als Volumen von
bestimmten dreidimensionalen Polyedern im hyperbolischen Raum. Leibons
Funktional f\"ur hyperbolische Kreismuster l\"asst sich nicht aus unseren
Funktionalen herleiten. Wir konstruieren jedoch ein weiteres Funktional, aus
dem sowohl Leibon's als auch unser Funktional hergeleitet werden kann. Durch
die umgekehrte Legendretransformation von Leibons Funktional erhalten wir ein
neues Variationsprinzip f\"ur hyperbolische Kreismuster.

Wir pr\"asentieren Java Software zur Berechnung und Visualisierung von
Kreismustern. 

}

\selectlanguage{american}

\tableofcontents

\mainmatter{}

\chapter{Introduction}
\label{cha:intro}

\section{Existence and uniqueness theorems}
\label{sec:ex_and_uni}
 
A {\em circle packing}\/ is a configuration of circular discs in a surface
such that the discs may touch but not overlap. We consider only circle
packings in surfaces of constant curvature.
\begin{figure}%
\includegraphics[width=0.475\textwidth, viewport=36 100 180 208, clip=true]{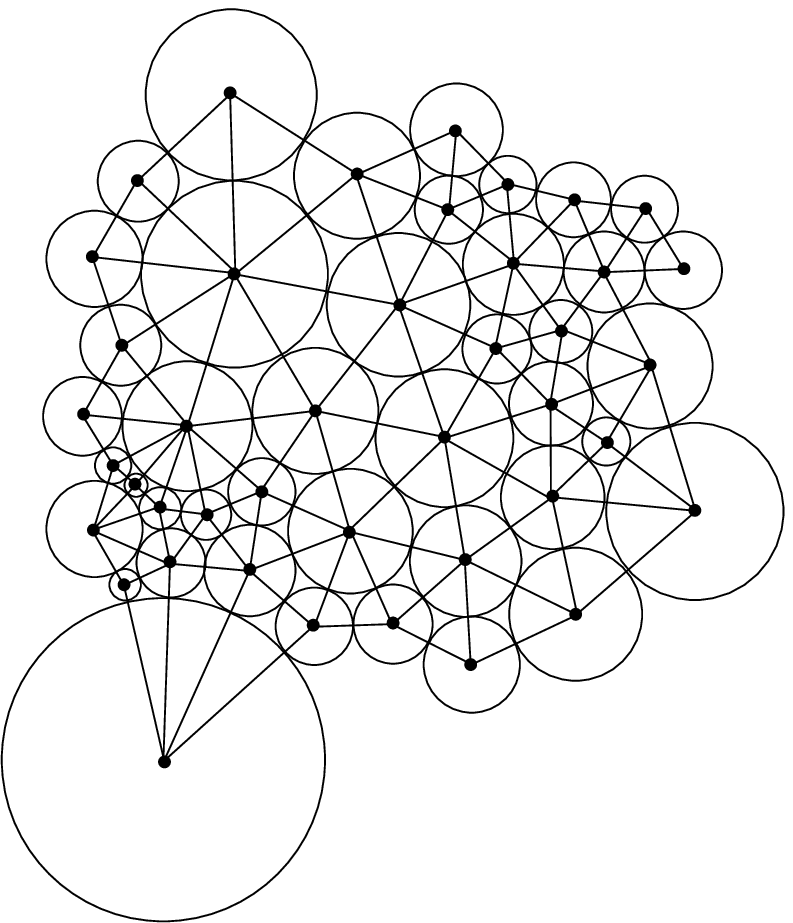}%
\hfill%
\includegraphics[width=0.475\textwidth, viewport=35 60 211 192,clip=true]{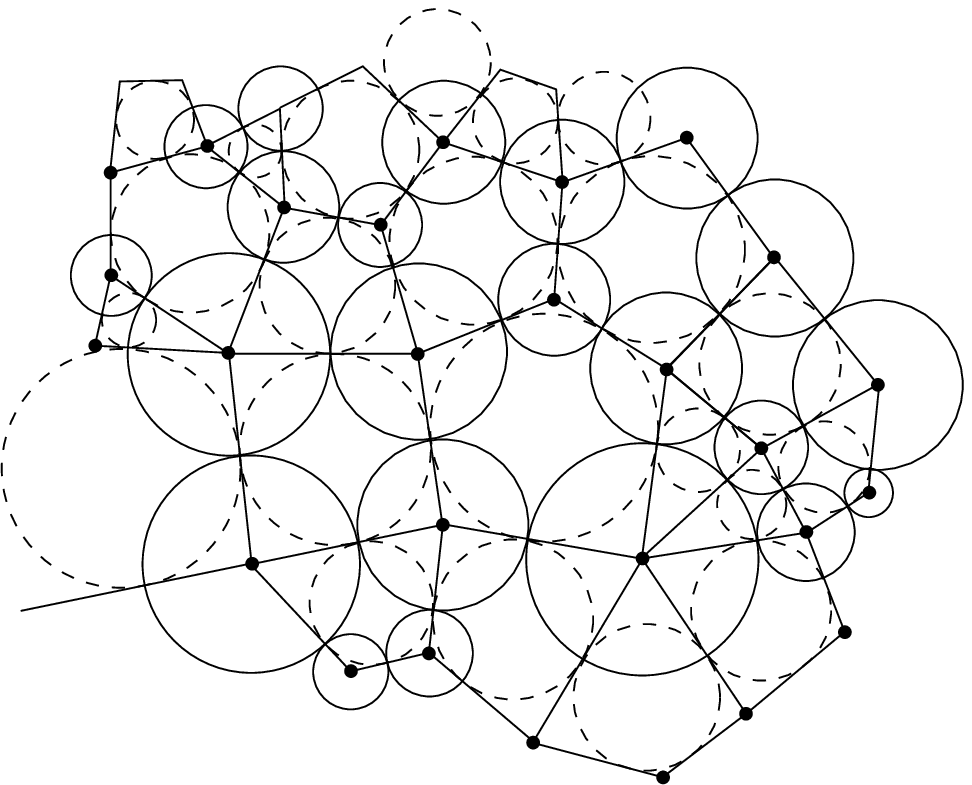}%
\caption{{\em Left:}\/ A circle packing corresponding to a
  triangulation. {\em Right:}\/ A pair or orthogonally intersecting circle
  packings corresponding to a cell
  decomposition.}\label{fig:packing}%
\end{figure}
Connecting the centers of touching discs by geodesics as in
figure~\ref{fig:packing}, one obtains an embedded graph, the {\em adjacency
  graph}\/ of the packing. Consider the case when the adjacency graph
triangulates the surface as in figure~\ref{fig:packing} ({\em left}\/).  The
following theorem of Koebe\ \cite{koebe36} answers the question:
Given an abstract triangulation of the sphere, does there exist a circle
packing whose adjacency graph is a geometric realization of the abstract
triangulation?

\begin{theorem}[Koebe]\label{thm:koebe}
  For every abstract triangulation of the sphere there is a circle packing
  whose adjacency graph is a geometric realization of the triangulation. The
  circle packing corresponding to a triangulation is unique up to M\"obius
  transformations of the sphere.
\end{theorem}

Now consider the case when the adjacency graph of a circle packing gives rise
to a cell decomposition whose faces are not necessarily triangles, as in
figure~\ref{fig:packing} ({\em{}right}\/). To characterize the cell
decompositions of the sphere which correspond to circle packings, we need the
following definition.

\begin{definition}
  A cell complex is called {\em regular}\/ if the characteristic maps, that
  map closed discs onto the closed cells, are homeomorphisms.
  A cell complex is called {\em strongly
    regular}\/ if it is regular and the intersection of two closed cells is
  empty or equal to a single closed cell.
\end{definition}
\begin{remark}
  \label{rem:regular}
  Suppose a cell complex $\Sigma$ is in fact a cell decomposition of a
  compact surface without boundary. One obtains the following conditions for
  $\Sigma$ being regular and strongly regular.
  
  The cell decomposition $\Sigma$ is regular if and only if the following
  conditions hold.
  
  \makebox[2em][l]{\em{}(i)} Each edge is incident with two vertices. (There
  are no loops.)

  \makebox[2em][l]{\em{}(ii)} Each edge is incident with two faces. (There
  are no stalks.)
  
  \makebox[2em][l]{\em{}(iii)} If a vertex $v$ and a face $f$ are incident,
  there are exactly two edges incident with both $v$ and $f$.

  The cell decomposition $\Sigma$ is strongly regular if and only if, in
  addition, the following conditions hold.
  
  \makebox[2em][l]{\em{}(iv)} No two edges are incident with the same two
  vertices.
  
  \makebox[2em][l]{\em{}(v)} No two edges are incident with the same two
  faces.
  
  \makebox[2em][l]{\em{}(vi)} If each of two faces is incident with each of
  two vertices, then there is an edge which is incident with both faces and
  both vertices.

  The above characterization implies: The cell decomposition $\Sigma$ is
  (strongly) regular if and only if its Poincar\'e-dual decomposition is
  (strongly) regular.
\end{remark}

A cell decomposition of the sphere which arises from a circle packing is
strongly regular. Conversely, Koebe's theorem implies that every strongly
regular cell decomposition of the sphere comes from a circle packing. (Simply
triangulate each face by adding an extra vertex inside and connecting it to
the original vertices.) However, the corresponding packing is generally not
unique up to M\"obius transformations.  The following theorem is a
generalization of Koebe's theorem which retains uniqueness by requiring the
existence of a second packing of orthogonally intersecting circles as in
figure~\ref{fig:packing} {\em{}(right)}.

\begin{theorem}\label{thm:ortho}
  For every strongly regular cell decomposition of the sphere, there exists a
  pair of circle packings with the following properties: The adjacency graph
  of the first packing is a geometric realization of the given cell
  decomposition. The adjacency graph of the second packing is a geometric
  realization of the Poincar\'e dual of the given cell decomposition.
  Therefore, to each edge there correspond four circles which touch in pairs.
  It is required that these pairs touch in the same point and intersect each
  other orthogonally.
  
  The pair of circle packings is unique up to M\"obius transformations of the
  sphere.
\end{theorem}

In the case of a circle packing corresponding to a triangulation as in
Koebe's theorem, the second orthogonal packing always exists.  Thus,
theorem~\ref{thm:ortho} is indeed a generalization of Koebe's theorem.
The first published proof is probably due to Brightwell\ and Scheinerman\ 
\cite{brightwell_scheinerman93}. They do not claim to give the first proof.
The works of Thurston \cite{thurston} and Schramm \cite{schramm92} (see
theorem~\ref{thm:schramm} below) indicate that the theorem was well
established at the time.

Associated with a polyhedron in $\R^3$ is a cell decomposition of the sphere
representing its combinatorial type. We say that the polyhedron is a
(geometric) realization of the cell decomposition. Steinitz' representation
theorem for convex polyhedra in $\R^3$ says that a cell decomposition of the
sphere represents the combinatorial type of a convex polyhedron if and only
if it is strongly regular \cite{steinitz22}, \cite{steinitz_rademacher34}.
Theorem~\ref{thm:ortho} implies the following stronger representation theorem
for convex polyhedra in $\R^3$ (see Ziegler~\cite{ziegler95}, theorem~4.13 on
p.~118).

\begin{theorem}
  \label{thm:strong_steinitz}
  For every strongly regular cell decomposition of the sphere there is a
  convex polyhedron in $\R^3$ which realizes it, such that the edges of the
  polyhedron are tangent to the unit sphere. Such a geometric realization is
  unique up to projective transformations which fix the sphere.
  
  Simultaneously, there is a polyhedron with edges tangent to the sphere
  which realizes the Poincar\'e-dual cell decomposition such that corresponding
  edges of the two polyhedra intersect each other orthogonally and touch the
  sphere in their point of intersection.
  
  Among the projectively equivalent polyhedral realizations, there is one and
  up to isometries only one with the property that the barycenter of the
  points where the edges touch the sphere is the center of the sphere. Every
  topological symmetry of the cell decomposition corresponds to an isometry
  of this polyhedral realization.
\end{theorem}

The following much more general theorem is due to Schramm\ \cite{schramm92}.
The proof is based on topological methods which are beyond the scope of this
thesis.

\begin{theorem}[Schramm]
  \label{thm:schramm}
  Let $P$ be a (3-dimensional) convex polyhedron, and let $K\subset\R^3$ be a
  smooth strictly convex body. Then there exists a convex polyhedron
  $Q\subset\R^3$ combinatorially equivalent to P whose edges are tangent to
  $K$.
\end{theorem}

\begin{figure}[tb]%
\centering \includegraphics[width=0.7\textwidth, keepaspectratio, viewport=55 100 400 400, clip]{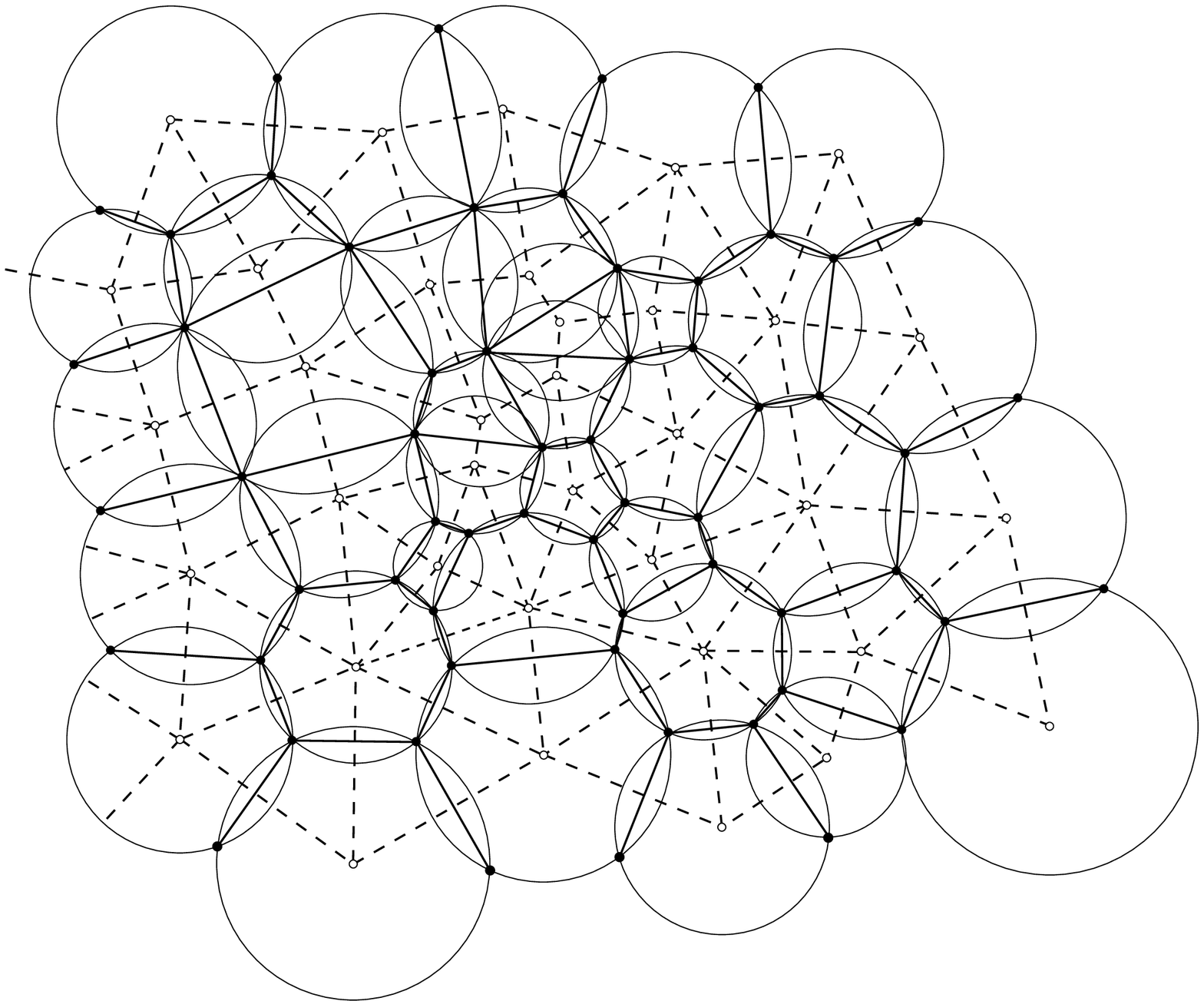}%
\caption{A Delaunay{} type circle pattern with Delaunay{} and Dirichlet{} 
  cell  decompositions.}%
\label{fig:pattern_altogether}%
\end{figure}

The two circle packings of theorem~\ref{thm:ortho} form a pattern of
orthogonally intersecting circles. By way of a further generalization, one
may consider circle patterns with circles intersecting at arbitrary angles.
\begin{quasidefinition}
  A {\em circle pattern}\/ is a configuration of circles in a constant
  curvature surface which corresponds in some way to a cell decomposition of
  the surface. According as the constant curvature is positive, zero, or
  negative, we speak of {\em{}spherical}, {\em{}euclidean}, or
  {\em{}hyperbolic circle patterns}.
\end{quasidefinition}
To obtain a real definition, the correspondence between circle pattern and
cell decomposition has to be specified. We will only be concerned with a
special class of circle patterns which are connected to Delaunay{} cell
decompositions. To be precise, we call them `Delaunay{} type circle patterns',
but we will usually refer to them simply as `circle patterns'.
Figure~\ref{fig:pattern_altogether} shows an example.
\begin{definition}
  A {\em Delaunay{} decomposition}\/ of a constant curvature surface is a
  cellular decomposition such that the boundary of each face is a geodesic
  polygon which is inscribed in a circular disc, and these discs have no
  vertices in their interior. The Poincar\'e-dual decomposition of a
  Delaunay{} decomposition with the centers of the circles as vertices and
  geodesic edges is a {\em{}Dirichlet{}
    decomposition}\/\index{Dirichlet{}!decomposition} (or {\em Voronoi{}
    diagram}\/\index{Voronoi{} diagram}).  A {\em Delaunay{} type circle
    pattern}\/ is the circle pattern formed by the circles of a Delaunay{}
  decomposition.  More generally, we allow the constant curvature surface to
  have cone-like singularities in the vertices and in the centers of the
  circles.
\end{definition}

Figure~\ref{fig:pattern_altogether} shows a Delaunay{} decomposition (black
vertices and solid lines), the dual Dirichlet{}
decomposition\index{Dirichlet{}!decomposition} (white vertices and dashed
lines) and the corresponding circle pattern. The faces of the Delaunay{}
decomposition correspond to circles. The vertices are intersection points of
circles.

\begin{figure}[tb]%
\centering%
\input{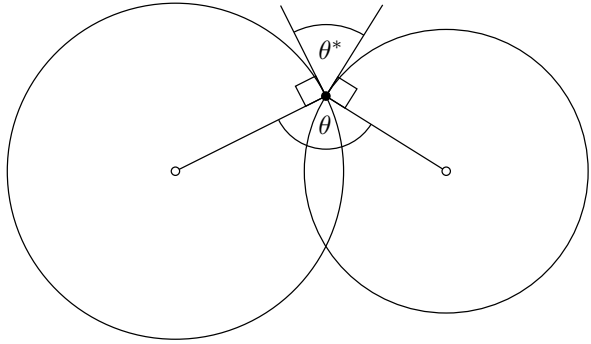}%
\caption{Interior intersection angle $\theta^*$ and exterior intersection
  angle $\theta$ of two circles. $\theta+\theta^*=\pi$.}%
\label{fig:intersection_angle}%
\end{figure}

\begin{remark}
  In section~\ref{sec:quad_graphs}, we will give an alternative and slightly
  more general definition for Delaunay{} type circle
  patterns. 
  
  A different class of circle patterns, {\em{}Thurston{} type circle
    patterns}, has been introduced by Thurston{}~\cite{thurston}. Here, the
  circles correspond to the faces of a cell decomposition, but the vertices
  do not correspond to intersection points. All vertices have degree $3$. Two
  circles corresponding to faces which share a common edge intersect (or
  touch) with an interior intersection angle $\theta^*$ (see
  figure~\ref{fig:intersection_angle}) satisfying $0\leq \theta^*\leq\pi/2$.

  Those Thurston type circle patterns in the sphere with the property that
  the sum of the three angles $\theta$\/ around each vertex is less or equal
  to $2\pi$ correspond to hyperbolic polyhedra of finite volume with dihedral
  angles at most $\pi/2$. Such polyhedra have been studied by
  Andreev{}~\cite{andreev70a}, \cite{andreev70b}. Delaunay type circle
  patterns in the sphere (without cone-like singularities) correspond to
  convex hyperbolic polyhedra with finite volume and all vertices on the
  infinite boundary of hyperbolic space.
\end{remark}

From a Delaunay{} type circle pattern, one may obtain the following data. (For
the definition of {\em interior}\/ and {\em exterior intersection angle}\/ of
two circles see figure~\ref{fig:intersection_angle}. We will always denote
the exterior intersection angle by $\theta$ and the interior intersection
angle by $\theta^*$. Note that $\theta+\theta^*=\pi$.)

\begin{itemize}
\item A cell decomposition $\Sigma$ of a 2-dimensional manifold.
\item For each edge $e$ of $\Sigma$ the exterior (or interior) intersection
  angle $\theta_e$ (or $\theta^*_e$). It satisfies
  $0<\theta_e<\pi$.
\item For each face $f$ of $\Sigma$ the cone angle $\Phi_f>0$ of the
  cone-like singularity at the center of the circle corresponding to $f$. (If
  there is no cone-like singularity at the center, then $\Phi_f=2\pi$.)
\end{itemize}

Note that the cone angle $\Theta_v$ at a vertex $v$ of $\Sigma$ is already
determined by the intersection angles $\theta_e$:
\begin{equation}
  \label{eq:capital_theta}
  \Theta_v = {\textstyle\sum}\, \theta_e,
\end{equation}
where the sum is taken over all edges around $v$. (See
figure~\ref{fig:intersection_angle}.) If $\Theta_v=2\pi$, there is no
singularity at $v$. The curvature in a vertex $v$ is 
\begin{equation}
  \label{eq:curvature_v}
  K_v=2\pi-\Theta_v,
\end{equation}
and the curvature in the center of a face $f$ is 
\begin{equation}
  \label{eq:curvature_f}
  K_f=2\pi-\Phi_f.
\end{equation}
The following theorems assume that the surface is closed, and that there are
no cone-like singularities. Only the last theorem~\ref{thm:fundamental} deals
with the general case.

Consider the following problem: Given such data as listed above, is there a
corresponding circle pattern? If so, is it unique?  The following theorem of
Rivin{}~\cite{rivin96} gives an answer for circle patterns in the sphere
without cone-like singularities. 

\begin{theorem}[Rivin{}, circle pattern version]
  \label{thm:rivin_circ}
  Let $\Sigma$ be a strongly regular cell decomposition of the sphere and let
  an angle $\theta_e$ with $0<\theta_e<\pi$ be given for every edge $e$ of
  $\Sigma$. Let $\Sigma^*$ be the Poincar\'e dual decomposition of $\Sigma$,
  and for each edge $e$ of $\Sigma$ denote the dual edge of $\Sigma^*$ by
  $e^*$.
  
  A Delaunay{} pattern corresponding to $\Sigma$ with exterior intersection
  angles $\theta_e$ exists if and only if the conditions {\it (i)}\/ and {\it
    (ii)}\/ are satisfied.
  
  \smallskip
  \makebox[1.5em][l]{\it (i)} If some edges $e^*_1,\ldots,e^*_n$ form the
  boundary of a face of $\Sigma^*$, then
  \begin{equation*}
    \sum\theta_{e_j} = 2\pi.
  \end{equation*}
  
  \makebox[1.5em][l]{\it (ii)} If some edges $e^*_1,\ldots,e^*_n$ form a
  closed path of $\Sigma^*$ which is not the boundary of a face, then
  \begin{equation*}
    \sum\theta_{e_j} > 2\pi.
  \end{equation*}
  
  \smallskip
  If it exists, the circle pattern is unique up to M\"obius
  transformations of the sphere.
\end{theorem}

The paper cited above deals with polyhedra in hyperbolic $3$-space with
vertices on the infinite boundary. The above theorem has the
following, equivalent, form.

\begin{theorem}[Rivin{}, ideal hyperbolic polyhedra version]
  \label{thm:rivin_poly}
  Let $\Sigma$ be a strongly regular cell decomposition of the sphere and let
  an angle $\theta_e$ with $0<\theta_e<\pi$ be given for every edge $e$ of
  $\Sigma$. Let $\Sigma^*$ be the Poincar\'e dual decomposition of $\Sigma$,
  and for each edge $e$ of $\Sigma$ denote the dual edge of $\Sigma^*$ by
  $e^*$.
  
  A polyhedron in hyperbolic $3$-space with vertices on the infinite boundary
  which realizes $\Sigma$ and has exterior dihedral angles $\theta_e$ exists
  if and only if the conditions~{\it (i)}\/ and~{\it (ii)}\/ are satisfied.
  
  \smallskip \makebox[1.5em][l]{\it (i)} If the edges $e^*_1,\ldots,e^*_n$
  form the boundary of a face of $\Sigma^*$, then
  \begin{equation*}
    \sum\theta_{e_j} = 2\pi.
  \end{equation*}
    
  \makebox[1.5em][l]{\it (ii)} If the edges $e^*_1,\ldots,e^*_n$ form a
  closed path of $\Sigma^*$ which is not the boundary of a face, then
  \begin{equation*}
    \sum\theta_{e_j} > 2\pi.
  \end{equation*}
  
  \smallskip If it exists, the polyhedron is unique up to an isometry of
  hyperbolic space.
\end{theorem}

Theorem~\ref{thm:higher_genus} below is a generalization of Rivin{}'s
theorem~\ref{thm:rivin_circ} to higher genus surfaces. Only the case of
oriented surfaces has to be considered: Because of the uniqueness claim,
non-orientable surfaces may be treated by applying the theorem to the
orientable double cover.

\begin{figure}[tb]%
\hfill%
\includegraphics[width=0.4\textwidth,keepaspectratio,clip]{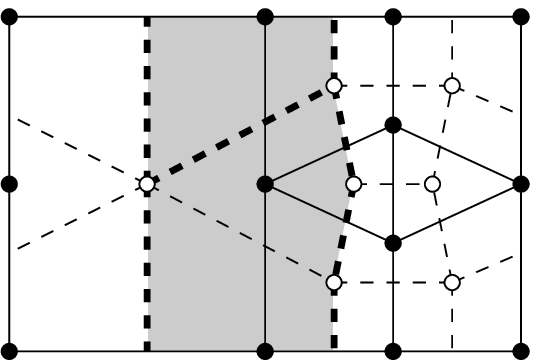}%
\hfill%
\includegraphics[width=0.44\textwidth,keepaspectratio,clip]{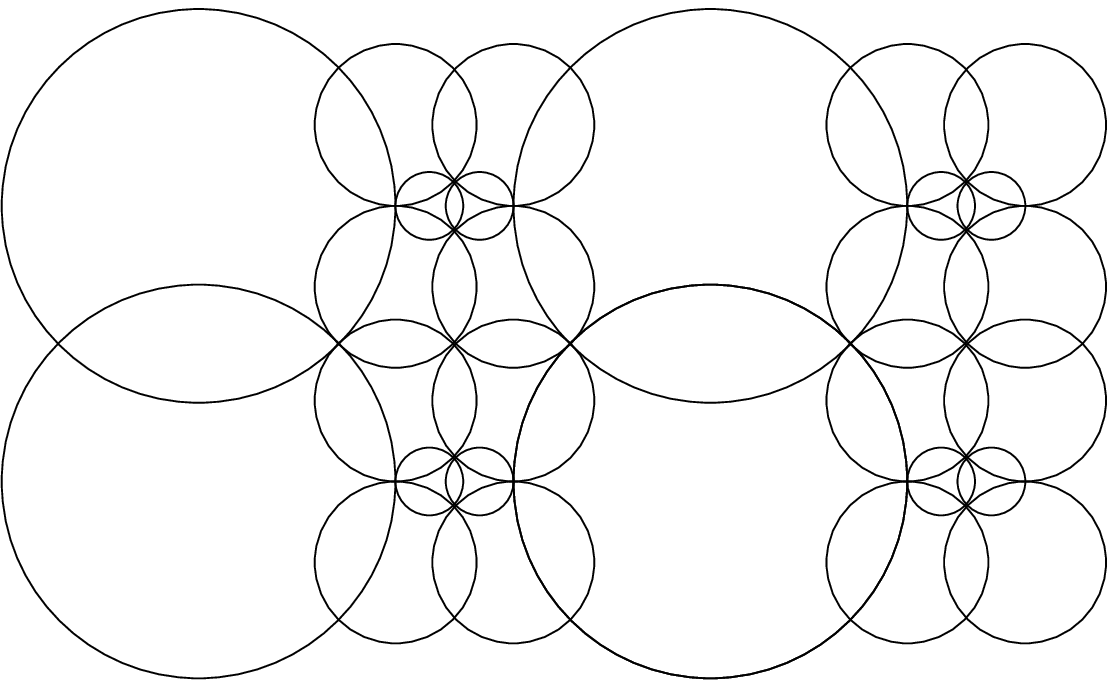}%
\hfill{}%
\caption{{\em{}Left:}\/ A cell decomposition of a torus (solid lines and
  black dots) and its Poincar\'e{} dual (dashed lines and white dots). The
  top and bottom side, and the left and right side of the rectangle are
  identified.  {\em{}Right:}\/ The corresponding double periodic circle
  pattern with orthogonally intersecting circles.}%
\label{fig:torus}%
\end{figure}

In the higher genus case, it is too restrictive to allow only strongly
regular cell decompositions. For example, figure~\ref{fig:torus} shows a cell
decomposition of the torus which is not even regular and the corresponding
circle pattern with orthogonally intersecting circles. With appropriate
cone-like singularities in the vertices and centers of circles, even faces
whose boundary is a loop and vertices of degree one are possible. As a
consequence of this tolerance regarding which cell decompositions are
acceptable, theorems~\ref{thm:higher_genus} and~\ref{thm:fundamental} are
really only true for Delaunay{} type circle patterns in the sense of the
slightly more general definition in section~\ref{sec:quad_graphs}.

The condition of theorem~\ref{thm:higher_genus} does not involve cycles in
the Poincar\'e{} dual decomposition like theorem~\ref{thm:rivin_circ}.
Instead, one has to consider cellular immersions of discs, whose interior is
embedded, into the Poincar\'e{} dual decomposition. The image of the boundary
path of the disc is a closed path in which some edges might appear twice. For
example, the shaded area in figure~\ref{fig:torus} {\em{}(left)}\/ is the
image of such a disc, where the heavy dashed lines are the image of its
boundary. One edge appears twice in the image of the boundary path. The
condition of theorem~\ref{thm:higher_genus} is that the sum of $\theta$\/
over such a boundary, where the edges are counted with appropriate
multiplicities, is at least $2\pi$.

\begin{theorem}\label{thm:higher_genus}
  Let $\Sigma$ be a cell decomposition of a closed compact oriented surface
  of genus $g>0$. Suppose exterior intersection angles are prescribed by a
  function $\theta\in(0,\pi)^E$ on the set $E$ of edges. Then there exists a
  circle pattern corresponding to these data in a surface of constant
  curvature (equal to\/ $0$ if $g=1$ and equal to $-1$ if $g>1$), if and only
  if the following condition is satisfied.
  
  \smallskip{} Suppose $\Delta\rightarrow\Sigma^*$ is any cellular immersion of
  a cell decomposition of the closed disc $\Delta$ into the Poincar\'e{} dual
  $\Sigma^*$ of $\Sigma$ which embeds the interior of $\Delta$.  Let ${\hat
    e}_1,\ldots,{\hat e}_k$ be the boundary edges of $\Delta$, let be
  $e^*_1,\ldots,e^*_k$ their images in $\Sigma^*$, and let $e_1,\ldots,e_k$
  be their dual edges in $\Sigma$. (An edge of $\Sigma$ may appear
  twice among the $e_j$.) Then
  \begin{equation}
    \label{eq:theta_sum_inequality}
    \sum_{j=1}^k \theta_{e_j} \geq 2\pi,
  \end{equation}
  where equality holds if and only if $\Delta$ has only one face.
  
  \smallskip{} In the case they exists, the constant mean curvature surface
  and the circle pattern are unique up to similarity, if $g=1$, and unique
  up to isometry, if $g>1$.
\end{theorem}

Schlenker~\cite{schlenker02} independently proves an existence and uniqueness
result for hyperbolic manifolds with polyhedral boundary, all vertices at
infinity, and prescribed dihedral angles.  Theorem~\ref{thm:higher_genus}
follows as a special case.  Interestingly, to obtain the general result,
Schlenker needs to first prove this special case separately.

\begin{figure}[tb]%
\hfill{}%
\includegraphics[width=0.45\textwidth]{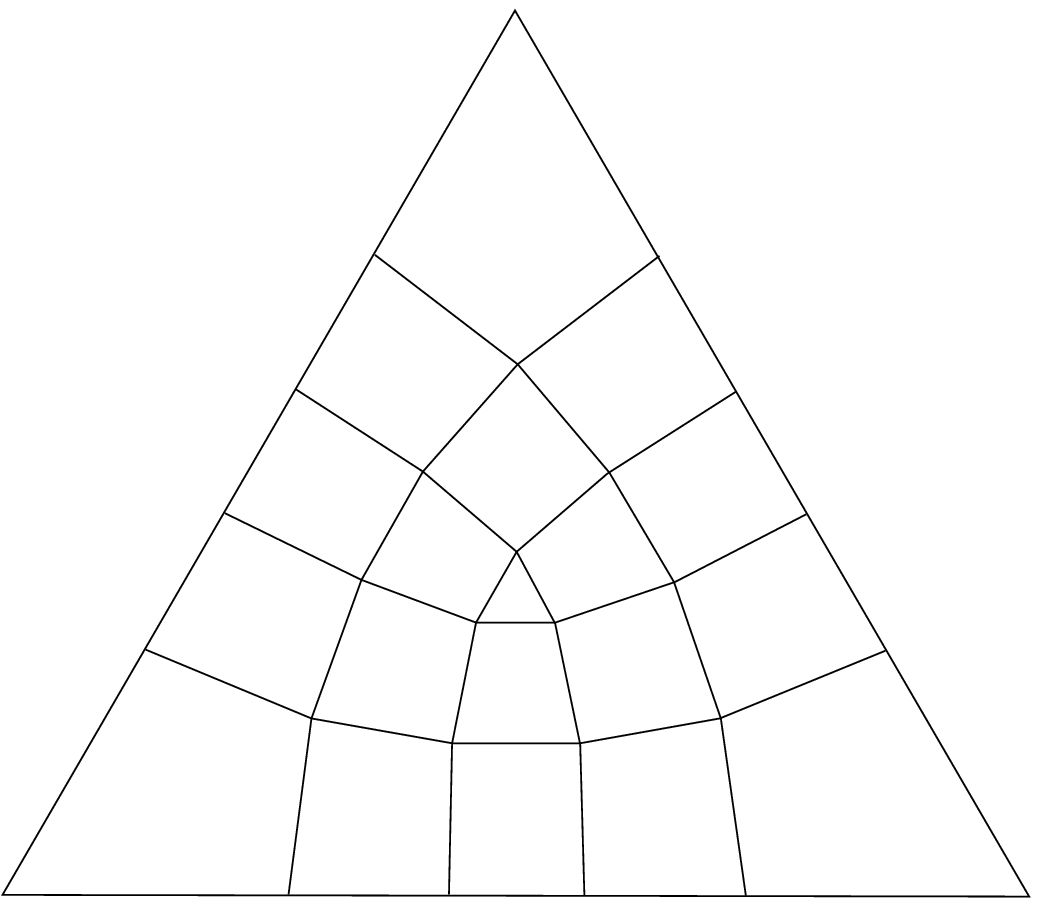}%
\hfill{}%
\includegraphics[width=0.45\textwidth]{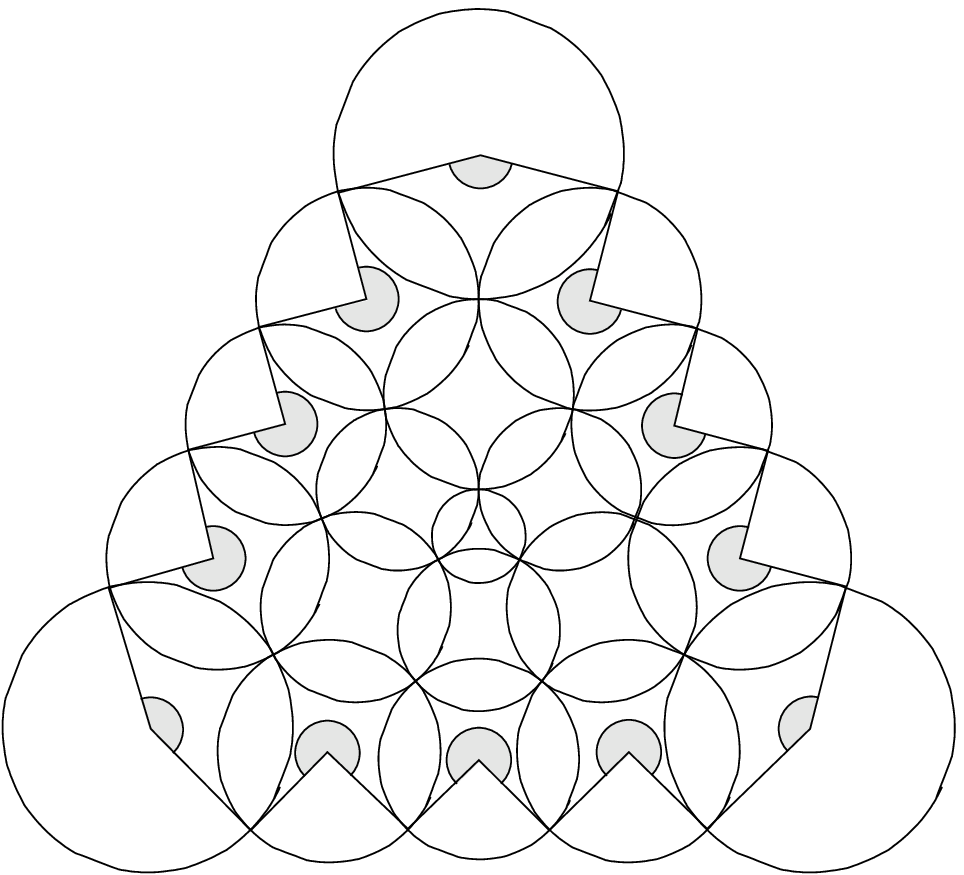}%
\hfill{}%
\caption{{\em{}Left:}\/ A cell decomposition of the disc. {\em{}Right:}\/ 
  A corresponding circle pattern with orthogonally intersecting circles.
  The shaded angles are the angles $\Phi$ for boundary faces. Here, they
  are $5\pi/6$ for the three corner circles and $3\pi/2$ for the other
  boundary circles.}
\label{fig:bounded_pattern}%
\end{figure}

We deduce theorem~\ref{thm:higher_genus} from the following more technical,
but also more general theorem~\ref{thm:fundamental}. It is not assumed that
$\theta$\/ sums to $2\pi$ around each vertex. Hence there may be cone-like
singularities in the vertices with cone angle $\theta$\/ given by
equation~\eqref{eq:capital_theta}. Also, cone-like singularities with
prescribed cone angle $\Phi_f$ are allowed in the centers of circles.
Furthermore, it applies also to surfaces with boundary.  For a boundary face
$f$, the angle $\Phi_f$ does not prescribe a cone angle, but the angle
covered by the neighboring circles, as shown in
figure~\ref{fig:bounded_pattern} {\em(right)}. These angles on the boundary
constitute {\em{}Neumann{} boundary conditions}.  Alternatively, one might
also prescribe the radii of the boundary circles.  This would constitute
{\em{}Dirichlet{} boundary conditions}\index{Dirichlet{}!boundary
  conditions}. We consider only the Neumann problem.

\begin{theorem}\label{thm:fundamental}
  Let $\Sigma$ be a cell decomposition of a compact oriented surface, with or
  without boundary. Suppose (interior) intersection angles are prescribed by
  a function $\theta^*\in(0,\pi)^{E_0}$ on the set $E_0$ of interior edges. Let
  $\Phi\in(0,\infty)^F$ be a function on the set $F$ of faces, which
  prescribes, for interior faces, the cone angle and, for boundary faces, the
  Neumann boundary conditions.
  
  (i) A euclidean circle pattern corresponding to these data exists if
  and only if the following condition is satisfied:
  
  \begin{list}{}{\setlength{\leftmargin}{\parindent}\setlength{\rightmargin}{\parindent}}
  \item %
    If $F'\subseteq F$ is any nonempty set of faces and $E'\subseteq E_0$ is
    the set of all interior edges which are incident with any face $f\in F'$,
    then
    \begin{equation}
      \label{eq:phi_theta_sum_inequality}
      \sum_{f\in F'} \Phi_f \leq \sum_{e\in E'} 2\theta^*_e,
    \end{equation}
    where equality holds if and only if $F'=F$.
  \end{list}

  \noindent{}
  If it exists, the circle pattern is unique up to similarity.
  
  (ii) A hyperbolic circle pattern corresponding to
  these data exists, if and only if, in the above condition, strict
  inequality holds also in the case $F'=F$.
  If it exists, the circle pattern is unique up to hyperbolic isometry.
\end{theorem}

Similar results were obtained by Bowditch~\cite{bowditch91},
Garrett~\cite{garrett92}, Rivin~\cite{rivin99}, and Leibon~\cite{leibon02}.
Bowditch treats the euclidean case for closed surfaces without cone-like
singularities in the centers of the circles. His proof is topological in
nature: It hinges on the fact that a certain function (essentially the
gradient of our functional $\Seuc$) is injective and proper.  Rivin extends
this result to surfaces with boundary.  Also, he considers not only singular
euclidean structures, but singular similarity structures. That is, he admits
not only cone-like singularities (with rotational holonomy) but singularities
with dilatational and rotational holonomy.  The proof uses his variational
principle~\cite{rivin94}. Leibon treats the hyperbolic case for closed
surfaces and without cone-like singularities. The proof uses his variational
principle (see section~\ref{sec:leibon_func}).  Garrett~\cite{garrett92}
obtains a similar theorem for euclidean and hyperbolic circle {\em
  packings}\/ with prescribed cone-like singularities. He considers the
Dirichlet boundary value problem (prescribed radii at the boundary). His
proof uses the relaxation method developed by Thurston~\cite{thurston}.

\section{The method of proof}
\label{sec:method_of_proof}

Chapter~\ref{cha:functionals} contains our proofs of the
theorems~\ref{thm:ortho}, \ref{thm:strong_steinitz}, \ref{thm:rivin_circ},
\ref{thm:rivin_poly}, \ref{thm:higher_genus}, and \ref{thm:fundamental}.
Here, we give a brief outline of the proofs of
theorems~\ref{thm:fundamental}, \ref{thm:higher_genus},
and~\ref{thm:rivin_circ}, which are the most involved. Most of the effort
goes into proving the fundamental theorem~\ref{thm:fundamental}, from which
the others are deduced. We extend methods introduced by Colin~de~Verdi\`ere
for circle packings~\cite{colin91}.

First, the geometric problem of constructing a
circle pattern is transformed into the analytic problem of finding the
(euclidean or hyperbolic) radii of the circles, which have to satisfy some
non-linear equations (closure conditions).  These non-linear closure
conditions turn out to be variational: The functionals $\Seuc$ and $\Shyp$
(defined in sections~\ref{sec:euc_func} and~\ref{sec:hyp_func}) are functions
of the radii, and their critical points are the solutions of the closure
conditions.  The functionals $\Seuc$, $\Shyp$ are convex (except for
scale-invariance in the euclidean case). This implies the uniqueness claims
of theorem~\ref{thm:fundamental} (section~\ref{sec:convexity}). The existence
claim is more difficult. We have to show that the functionals tend to
infinity if some radii go to zero or infinity. In
section~\ref{sec:coherent_angle_systems}, we show that this is the case if a
`coherent angle system' exists. This is a function on the oriented edges,
which satisfies a system of linear equations and inequalities. The existence
problem for circle patterns is thus reduced to the feasibility problem of a
linear program. In section~\ref{sec:proof_thm:fundamental}, we prove the
existence of a coherent angle system if the conditions of
theorem~\ref{thm:fundamental} are satisfied. This is done by interpreting a
coherent angle system as a feasible flow in a network (with capacity bounds
on the branches) and invoking the feasible flow theorem.

In section~\ref{sec:proof_thm_higher_genus}, theorem~\ref{thm:higher_genus}
for circle patterns in hyperbolic surfaces is deduced from
theorem~\ref{thm:fundamental} using methods of combinatorial topology. The
basis of this deduction is lemma~\ref{lem:general_euler}. We present a
self-contained proof of it in appendix~\ref{app:combi_top}.

In section~\ref{sec:proof_thm_rivin_circ} the analogous
theorem~\ref{thm:rivin_circ} for circle patterns in the sphere is deduced
from theorem~\ref{thm:higher_genus}. First, the problem is transferred to the
plane by stereographic projection. Then the proof proceeds in a similar way
as in the hyperbolic case. 

\section{Variational principles}
\label{sec:var_princ}

In chapter~\ref{cha:functionals}, we define the functionals $\Seuc$, $\Shyp$,
and $\Ssph$ for euclidean, hyperbolic, and spherical circle patterns. The
functional $\Ssph$ is not convex. Thus, we cannot use it to prove existence
and uniqueness theorems. The variables are (up to a coordinate
transformation) the (euclidean, hyperbolic, or spherical) radii of the
circles.  A Legendre transformation of these functionals
(section~\ref{sec:legendre}) leads to a new variational principle involving
{\em one}\/ new functional $\Shat$ for all geometries (euclidean, hyperbolic,
and spherical). The variables of $\Shat$ are certain angles; and the
variation is constrained to coherent angle systems. Depending on whether the
constraint involves euclidean, hyperbolic, or spherical coherent angle
systems (section~\ref{sec:coherent_angle_systems}), the critical points
correspond to euclidean, hyperbolic, or spherical circle patterns.

Colin~de~Verdi\`ere first used a variational principle to prove existence and
uniqueness for circle packings~\cite{colin91}. He constructs two functionals,
one for the euclidean case, one for the hyperbolic case. The variables are
the radii of the circles. Critical points correspond to circle packings.
Explicit formulas are given only for the derivatives of the functionals, not
for the functionals themselves. In section~\ref{sec:colin_func}, we derive
Colin~de~Verdi\`ere's functionals from our functionals~$\Seuc$ and $\Shyp$.
In particular, this effects the integration of Colin~de~Verdi\`ere's
differential formulas.

Apparently, Br\"agger~\cite{bragger92} had already tried to integrate
Colin~de~Verdi\`ere's formulas. He derives a new variational
principle for euclidean circle packings. The variables of his functional are
certain angles, and the variation is constrained to coherent angle
systems. This functional turns out to be related to Colin~de~Verdi\`ere's
functional by a Legendre transformation. In section~\ref{sec:bragger_func},
we derive it from our functional~$\Shat$. 

Rivin's functional for euclidean circle patterns~\cite{rivin94} is also
derived from $\Shat$ (section~\ref{sec:rivin_func}). It is less general,
because the cell decomposition is assumed to be a triangulation, and there
can be no curvature in the centers of circles.

Leibon~\cite{leibon99}, \cite{leibon02} derived a functional for hyperbolic
circle patterns which can be seen as a hyperbolic version of Rivin's
functional (section~\ref{sec:leibon_func}). It is therefore natural to expect
that Leibon's functional can be derived from $\Shat$ as well. However, this
is not the case. The Legendre dual of Leibon's functional
(section~\ref{sec:leibon_dual_func}) is not $\Shyp$, but a new
functional. Unfortunately, we cannot present an explicit formula for this
functional. In section~\ref{sec:common_ancestor}, we derive yet another
functional, from which both $\Shat$ and Leibon's functional can be derived.

At least since Br\"agger~\cite{bragger92}, there was an awareness of the fact
that the circle packing functionals have something to do with the volume of
hyperbolic $3$-simplices. Chapter~\ref{cha:volume} deals with the connection
between the circle pattern functionals and the volumes of hyperbolic
polyhedra. Schl\"afli's differential volume formula turns out to be the
unifying principle behind all circle pattern functionals. This geometric
approach is essential for the construction of the common ancestor of $\Shat$
and Leibon's functional.

Thurston's method to construct circle patterns~\cite{thurston} (implemented
in Stephenson's program {\tt CirclePack}~\cite{collins_stephenson03})
involves iteratively adjusting the radius of each circle so that the
neighboring circles fit around. This is equivalent to minimizing our
functionals ($\Seuc$, $\Shyp$) iteratively in each coordinate direction.

\section{Open questions}
\label{sec:open_questions}

There is a functional for Thurston type circle patterns (at least in the
euclidean case), its variables being the radii of the circles. In fact, Chow
and Luo~\cite{chow_luo02} show that the corresponding closure conditions are
variational. Can an explicit formula be derived? (For Thurston type circle
patterns with ``holes", we derive a functional in
section~\ref{sec:Thurston_with_holes}.)

One may also consider Thurston type circle patterns with non-intersecting
circles. Instead of the intersection angle, one prescribes the inversive
distance (an imaginary intersection angle) between neighboring
circles (see Bowers and Hurdal~\cite{bowers_hurdal03}). Is there a functional
for such circle patterns, and can one write an explicit formula? Can this
approach be used to prove existence and uniqueness theorems? These questions
are interesting, because inversive-distance circle patterns
may be the key to `conformally parametrized' polyhedral surfaces in $\R^3$.

Even though the spherical functional is not convex, may it be used to prove
existence and (M\"obius-)uniqueness theorems for branched circle patterns in
the sphere? (See also Bowers and
Stephenson~\cite{bowers_stephenson96}.) Can the representation
theorem~\ref{thm:strong_steinitz} be generalized for star-polyhedra? This
question is interesting because branched circle patterns in the
sphere can be used to construct `discrete minimal
surfaces'~\cite{bobenko_hoffmann_springborn03}.

Rodin and Sullivan showed that circle packings approximate conformal
mappings~\cite{rodin_sullivan87}. Schramm proved a similar result for `circle
patterns with the combinatorics of the square grid'~\cite{schramm97}. There
have been numerous refinements (for example, the proof of
`$C^{\infty}$-convergence' by He and Schramm~\cite{he_schramm98}), but all
convergence results deal with circle patterns with the topology of the disc
and regular combinatorics (square grid or hexagonal). Can the variational
approach help in proving convergence results for circle patterns with
non-trivial topology and irregular combinatorics?

\section{Acknowledgments}
\label{sec:acknowledge}

I would like to thank my academic advisor, Alexander I.\ Bobenko, for being a
great teacher, for his judicious guidance, and for letting me work at my own
pace. I also thank Ulrich Pinkall, not only but in particular for help with
the proof of the strong Steinitz theorem. I thank G\"unter M.\ Ziegler for
his kind and active interest in my work. His expert advice on discrete and
combinatorial matters has been extremely helpful.  

I am also indebted to my parents, but that is beyond the scope of this
thesis.

While I was working on this thesis, I was supported by the DFG's
Sonderforschungsbereich 288. Some of the material is also contained in a
previous article by the author~\cite{springborn03} and in joint articles with
Bobenko~\cite{bobenko_springborn02} and with Bobenko and
Hoffmann~\cite{bobenko_hoffmann_springborn03}.

\chapter[The functionals. Proof of existence and uniqueness theorems]%
{The functionals. Proof of the existence and uniqueness theorems}
\label{cha:functionals}

\section[Quad graphs]{Quad graphs and an alternative definition for Delaunay{} circle
  patterns}
\label{sec:quad_graphs}

\begin{figure}[tb]%
\hfill%
\includegraphics[keepaspectratio,clip]{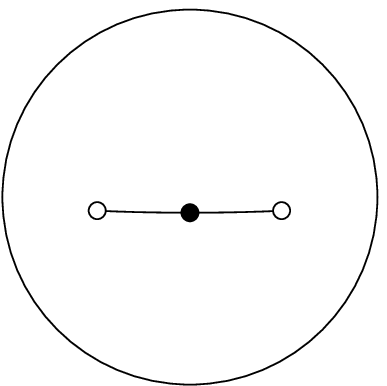}%
\hfill%
\includegraphics[width=0.5\textwidth{}]{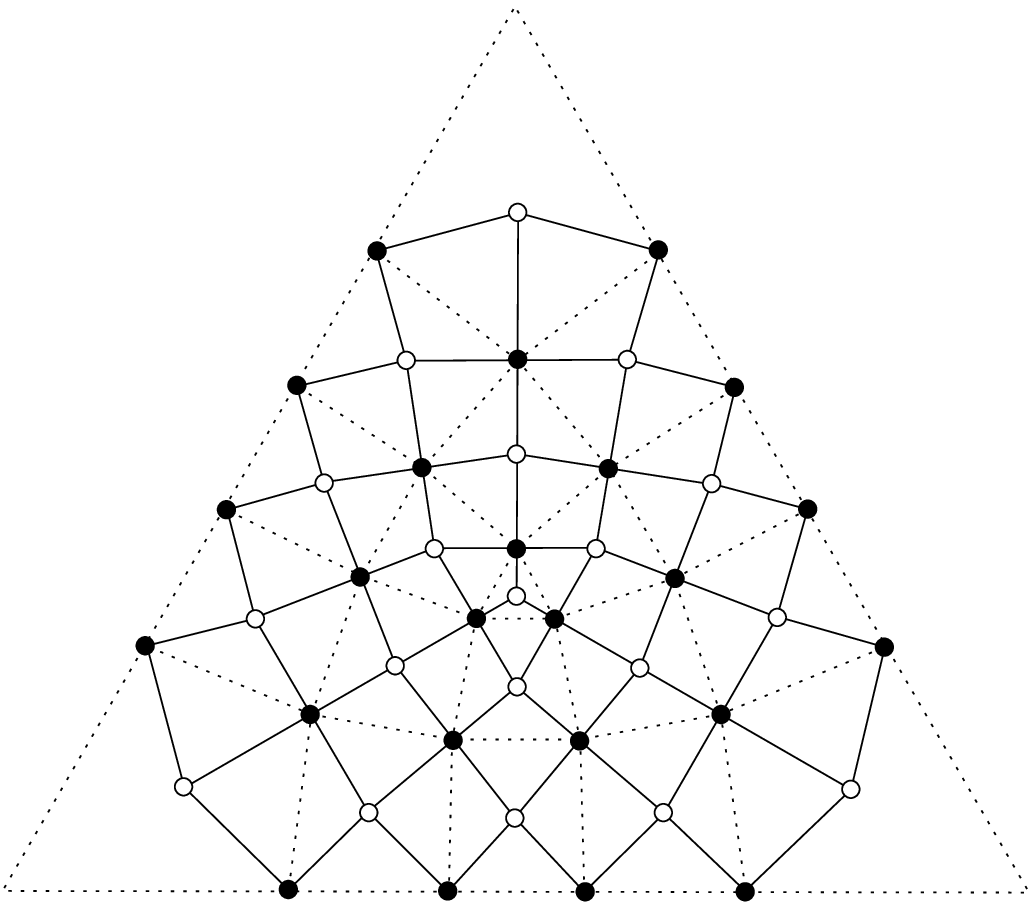}%
\hspace*{\fill}%
\caption{{\em Left:}~The smallest quad graph decomposition of the
  sphere. {\em Right:}~The quad graph corresponding to the cell
  decomposition of figure~\ref{fig:bounded_pattern} (dotted). Note that
  boundary edges do not correspond to quadrilaterals.}
\label{fig:quad_graph}
\end{figure}

A `quad graph' (the term was coined by Bobenko and
Suris~\cite{bobenko_suris02}) is a cell decomposition of a surface such
that the faces are quadrilaterals. We also demand that the vertices are
bicolored. On the other hand, we allow identifications on the boundary of a
face. For example, figure~\ref{fig:quad_graph} {\em (left)}\/ shows a quad
graph decomposition of the sphere with only one `quadrilateral'.
To put is more precisely:

\begin{definition}
  A {\em{}quad graph}\/ is a cell decomposition of a surface, such that each
  closed face is the image of a quadrilateral under a cellular map which
  immerses the open cells, and the vertices are colored black and white so
  that each edge is incident with one white and one black vertex.
\end{definition}

From any cell decomposition $\Sigma$ of a surface one obtains a quad graph
$\mathcal Q$ such that the correspondence between elements of $\Sigma$ and
elements of $\mathcal Q$ is as follows:

\smallskip{}
\begin{center}
  \begin{tabular}{c | c}
    $\Sigma$ & $\mathcal Q$ \\
    \hline{}
    vertices & black vertices \\
    faces & white vertices \\
    interior edges & quadrilaterals
  \end{tabular}  
\end{center}

\smallskip{}\noindent{}%
Figure~\ref{fig:quad_graph} {\em (right)}\/ shows an example which 
should make the construction clear. This construction may be reversed, such
that from every quad graph one obtains a cell decomposition. (The reverse
construction is not quite unique in the case of surfaces with boundary,
because one is free to insert any number boundary edges. Fortunately,
boundary edges play no role in our treatment of circle patterns.)

The following simple definition of Delaunay{} type circle patterns in terms of
quad graphs is a bit more general than the one in
section~\ref{sec:ex_and_uni}.

\begin{definition}
  A {\em{}(generalized) Delaunay{} type circle pattern}\/ is a quad graph in a
  constant curvature surface, possibly with cone-like singularities in the
  vertices, such that the edges are geodesic and all edges incident with the
  same white vertex have the same length.
\end{definition}

\begin{figure}[tb]%
\centering%
\includegraphics[keepaspectratio,clip]{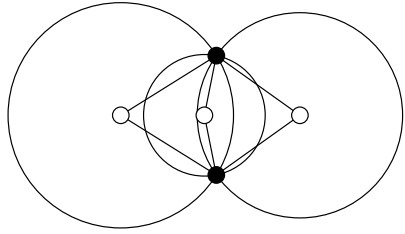}%
\caption{Not a Delaunay{} type circle pattern by the definition of 
  section~\ref{sec:ex_and_uni}.}%
\label{fig:digon_pattern}%
\end{figure}

This definitions allows for configurations as shown in
figure~\ref{fig:digon_pattern}. The corresponding cell decomposition has a
digon corresponding to the white vertex in the middle.

\section[Analytic formulation; euclidean case]%
{Analytic formulation of the circle pattern problem; euclidean case}
\label{sec:analytic}

Consider the following euclidean circle pattern problem: For a given finite
cell decomposition $\Sigma$ of a compact surface with or without boundary, a
given angle $\theta_e$ with $0<\theta_e<\pi$ for each interior edge $e$, and
a given angle $\Phi_f$ for each face $f$, construct a euclidean Delaunay type
circle pattern (as defined in section~\ref{sec:ex_and_uni}) with cell
decomposition $\Sigma$, intersection angles given by $\theta$, and cone
angles and Neumann{} boundary conditions given by $\Phi$.

We will reduce this problem to solving a set of nonlinear equations for the
radii of the circles.  The following lemma is the basis for this.

\begin{lemma}\label{lem:given_theta_and_r}
  Let $\Sigma$ be a cell decomposition of a compact surface, possibly with
  boundary. Let $\theta:\Eint\rightarrow(0,\pi)$ be a function on the set
  $\Eint$ of interior edges, and $r:F\rightarrow(0,\infty)$ be a function on
  the set $F$ of faces. Then there exists a unique euclidean circle pattern
  with cone-like singularities in the vertices and in the centers of circles
  such that the corresponding cell decomposition is $\Sigma$, the
  intersection angles are given by $\theta$ and the radii are given by $r$.
  
  The cone angle $\Theta_v$ in a vertex $v$ is given by
  $\Theta_v=\sum\theta_e$, where the sum is taken over all edges $e$ around
  $v$.  The cone angle in the center of an interior face $f_j$ (or the
  boundary angle for a boundary face) is
  \begin{equation}
    \label{eq:capital_phi_of_r}
    \Phi_{f_j} = 2\;
    \sum_{\makebox[0pt][r]{\scriptsize$f_j\circ$}\edge\makebox[0pt][l]{\scriptsize$\circ f_k$}}\;
    \frac{1}{2i}
    \log\frac{r_{f_j}-r_{f_k} e^{-i\theta_e}}{r_{f_j}-r_{f_k} e^{i\theta_e}}\,,
  \end{equation}
  where the sum is taken over all interior edges $e$ between the face $f_j$
  and its neighbors~$f_k$.
\end{lemma}

\begin{figure}[tb]
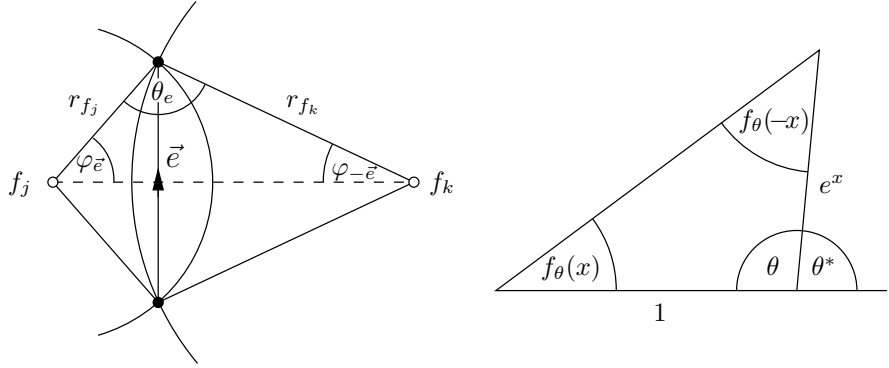
%
\hfill%
\input{kite.tex}%
\hfill%
\raisebox{0.5cm}{\input{f_theta.tex}}%
\hspace*{\fill}%
\caption{{\em Left:} A kite shaped quadrilateral of the quad graph. The
  oriented edge $\vece$ has the face $f_j$ on its left side and the face
  $f_k$ on its right
  side. {\em Right:} The function $f_{\theta}(x)$.}%
\label{fig:kite_f_theta}%
\end{figure}

\begin{proof}
  Given the cell decomposition, intersection angles, and radii, one
  constructs the corresponding circle pattern as follows. For each interior
  edge $e$ with faces $f_j$ and $f_k$ on either side, construct a euclidean
  kite shaped quadrilateral with side lengths $r_{f_j}$ and $r_{f_k}$ and
  angle $\theta_e$ as in figure~\ref{fig:kite_f_theta}~{\em (left)}.  Glue
  these quadrilaterals together to obtain a flat surface with cone-like
  singularities, and, in fact, the desired circle pattern. The uniqueness
  claim is obvious, as is the claim about $\Theta_v$.  For each oriented edge
  $\vece$, let $\varphi_{\vece}$ be half the angle covered by $\vece$ as seen
  from the center of the circle on the left side of $\vece$. See
  figure~\ref{fig:kite_f_theta}~{\em (left)}. Now,
  \begin{equation*}
    \Phi_f=2{\textstyle \sum}\varphi_{\vece},
  \end{equation*}
  where the sum is taken over all oriented edges in oriented the boundary of
  $f$, and
  \begin{equation}
    \label{eq:phi_euc}
    \varphi_{\vece}=\frac{1}{2i}
    \log\frac{r_{f_j}-r_{f_k} e^{-i\theta_e}}{r_{f_j}-r_{f_k} e^{i\theta_e}}.
  \end{equation}
  (The argument of a non-zero complex number $z$ is
  $\arg z=\frac{1}{2i}\log\frac{z}{\bar{z}}$.)
  Equation~\eqref{eq:capital_phi_of_r} follows.
\end{proof}

It is convenient to introduce the logarithmic radii 
\begin{equation}
  \label{eq:rho_of_r_euc}
  \rho=\log r
\end{equation}
as variables. Then, equation~\eqref{eq:phi_euc} may be rewritten as
\begin{equation}
  \label{eq:phi_e_of_rho_f_euc}
  \varphi_{\vece}=f_{\theta_e}(\rho_{f_k}-\rho_{f_j}),
\end{equation}
where, for $0<\theta<\pi$, the function $f_{\theta}:\R\rightarrow\R$ is
defined by
\begin{equation}
  \label{eq:f}
  f_\theta(x) := \frac{1}{2i}
  \log\frac{1-e^{x-i\theta}}{1-e^{x+i\theta}},
\end{equation}
and the branch of the logarithm is chosen such that 
\begin{equation*}
  0<f_\theta(x)<\pi.
\end{equation*}
In the following lemma, we list a few properties of the function
$f_\theta(x)$ for reference.

\begin{lemma}\label{lem:f_theta}\

  \smallskip{}
  (i) The remaining angles of a triangle with sides $1$ and $e^x$ and with
  enclosed angle $\theta$ are $f_{\theta}(x)$ and $f_{\theta}(-x)$, as shown
  in figure~\ref{fig:kite_f_theta} {\it (right)}.

  \smallskip
  (ii) The derivative of $f_\theta(x)$ is
  \begin{equation}
    \label{eq:f_prime}
    f_\theta'(x)=
    \frac{\sin\theta}{2(\cosh x-\cosh\theta)}>0,
  \end{equation}
  so $f_\theta(x)$ is strictly increasing. 
  
  \smallskip
  (iii) The function $f_\theta(x)$ satisfies the functional equation
  \begin{equation}
    \label{eq:f_symmetry}
    f_\theta(x)+f_\theta(-x)=\pi-\theta.
  \end{equation}

  \smallskip
  (iv) The limiting values of $f_\theta(x)$ are 
  \begin{equation}
    \label{eq:f_limits}
    \lim_{x\rightarrow-\infty}f_\theta(x)=0 \quad\text{and}\quad
    \lim_{x\rightarrow\infty}f_\theta(x)=\pi-\theta,
  \end{equation}

  \smallskip
  (v) For $0<y<\pi-\theta$, the inverse function is
  \begin{equation}
    \label{eq:f_inverse}
    f_\theta^{-1}(y)=\log\frac{\sin y}{\sin(y+\theta)}.
  \end{equation}

  \smallskip
  (vi) The integral of $f_{\theta}(x)$ is
  \begin{equation}
    \label{eq:f_integral}
    F_\theta(x):=\int_{-\infty}^x f_{\theta}(\xi)\,d\xi=\im\Li(e^{x+i\theta}),
  \end{equation}
  where $\Li(z)$ is the dilogarithm function; see appendix~\ref{app:dilog}.
\end{lemma}

By lemma~\ref{lem:given_theta_and_r}, the euclidean circle pattern problem is
equivalent to the non-linear equations~\eqref{eq:euler-lagrange_euc} below.

\begin{lemma}
  Given a cell decomposition $\Sigma$ of a compact surface with or without
  boundary, an angle $\theta_e$ with $0<\theta_e<\pi$ for each interior edge
  $e$, and an angle $\Phi_f$ for each face $f$. Suppose $r\in\R_{+}^{F}$ and
  $\rho\in\R^{F}$ are related by equation~\eqref{eq:rho_of_r_euc}. Then the
  following statements (i) and (ii) are equivalent:
  
  \smallskip
  (i) There is a euclidean circle pattern with radii $r_f$, intersection
  angles $\theta_e$ and cone/boundary angles $\Phi_f$.

  \smallskip
  (ii) For each face $f\in F$,
  \begin{equation}
    \label{eq:euler-lagrange_euc}
    \Phi_f-2\;
    \sum_{\makebox[0pt][r]{\scriptsize$f\,\circ$}\edge\makebox[0pt][l]{\scriptsize$\circ f_k$}}\;f_{\theta_e} (\rho_{f_k}-\rho_{f})=0,
  \end{equation}
  where the sum is taken over all interior edges $e$ between the face $f$ and
  its neighbors $f_k$.
\end{lemma}

\section[The euclidean functional]{The euclidean circle pattern functional}
\label{sec:euc_func}

The euclidean circle pattern functional defined below is a function of the
logarithmic radii $\rho_f$. Equations~\eqref{eq:euler-lagrange_euc} are the
conditions for a critical point. 

\begin{definition}
  \label{def:euc_functional}
  The {\em euclidean circle pattern functional}\/ is the function
  \begin{equation*}
    \nonumber
    \Seuc: \R^F\longrightarrow \R
  \end{equation*}
  \begin{multline}\label{eq:SEuc}
      \Seuc(\rho)= \\
      \sum_{f_j\circ\edge\circ f_k} \Big(
      \im\Li\big(e^{\rho_{f_k}-\rho_{f_j}+i\theta_e}\big)
      +\im\Li\big(e^{\rho_{f_j}-\rho_{f_k}+i\theta_e}\big)
      -(\pi-\theta_e)\big(\rho_{f_j}+\rho_{f_k}\big)\Big)\\ 
      +\sum_{\circ f}\Phi_f
      \rho_f.
  \end{multline}
  The first sum is taken over all interior edges $e$, and $f_j$ and $f_k$ are
  the faces on either side of $e$. (The terms are symmetric in $f_j$ and
  $f_k$, so it does not matter which face is considered as $f_j$ and which as
  $f_k$.) The second sum is taken over all faces $f$.
\end{definition}

\begin{lemma}
  \label{lem:func_euc}
  A function $\rho\in\R^F$ is a critical point of the euclidean circle pattern
  functional $\Seuc$, if and only if it satisfies
  equations~\eqref{eq:euler-lagrange_euc}. The critical points of $\Seuc$ are
  therefore in one-to-one correspondence with the solutions of the euclidean
  circle pattern problem.
\end{lemma}

\begin{proof}
  Using equations~\eqref{eq:f_integral} and~\eqref{eq:f_symmetry}, one
  obtains
  \begin{equation*}
    \frac{\partial\Seuc}{\partial\rho_{f}}=
    \Phi_f-2\;
    \sum_{\makebox[0pt][r]{\scriptsize$f\,\circ$}
      \edge\makebox[0pt][l]{\scriptsize$\circ f_k$}}\;
    f_{\theta_e} (\rho_{f_k}-\rho_{f})\,,
  \end{equation*}
  where the sum is taken over all edges $e$ between the face $f$ and its
  neighbors $f_k$.
\end{proof}

\section[The hyperbolic functional]{The hyperbolic circle pattern functional}
\label{sec:hyp_func}

This case is treated in the same fashion as the euclidean case. Of course,
the trigonometric relations are different:

\begin{lemma}
  \label{lem:phi_hyp}
  Suppose $r_1$ and $r_2$ are two sides of a hyperbolic triangle, the
  enclosed angle between them is $\theta$, and the remaining angles are
  $\varphi_1$ and $\varphi_2$, as shown in figure~\ref{fig:phiRhoTriangle}.
  Then
  \begin{equation}
    \label{eq:phi_of_rho_hyp}
    \varphi_{1}=f_{\theta}(\rho_{2}-\rho_{1})-f_{\theta}(\rho_{2}+\rho_{1}),
  \end{equation}
  where $f_\theta(x)$ is defined by equation~\eqref{eq:f}, and
  \begin{equation}
    \label{eq:rho_hyp}
    \rho = \log\tanh\frac{r}{2}.
  \end{equation}
  The inverse relation is
  \begin{equation}\label{eq:rho_of_phi_hyp}
    \rho_{1}=\frac{1}{2}\log
    \frac{\sin\Big(
      \frac{\displaystyle \theta^*-\varphi_{1}-\varphi_{2}}
      {\displaystyle 2}\Big)
      \sin\Big(
      \frac{\displaystyle \theta^*-\varphi_{1}+\varphi_{2}}
      {\displaystyle 2}\Big)}
    {\sin\Big(
      \frac{\displaystyle \theta^*+\varphi_{1}+\varphi_{2}}
      {\displaystyle 2}\Big)
      \sin\Big(
      \frac{\displaystyle \theta^*+\varphi_{1}-\varphi_{2}}
      {\displaystyle 2}\Big)},
  \end{equation}
  where $\theta^*=\pi-\theta$,\, $\varphi_{1,2}>0$, and $\varphi_1+\varphi_2<\theta^*$.
\end{lemma}

\begin{figure}
  \input{phiRhoTriangle.tex}
  \caption{} \label{fig:phiRhoTriangle}
\end{figure} 

Equation~\eqref{eq:phi_of_rho_hyp} is derived in appendix~\ref{app:trig}.
Equation~\eqref{eq:rho_of_phi_hyp} follows from \eqref{eq:phi_of_rho_hyp}
(and the corresponding equation for $\varphi_2$) by a straightforward
calculation using equations~\eqref{eq:f_inverse} and \eqref{eq:f_symmetry}.
Note that positive radii $r$ correspond to negative $\rho$.

In the hyperbolic case, the new variables $\rho$ are given by
equation~\eqref{eq:rho_hyp}. Instead of
equation~\eqref{eq:phi_e_of_rho_f_euc}, one has
\begin{equation}
  \label{eq:phi_e_of_rho_f_hyp}
  \varphi_{\vec e}
  =f_{\theta_e}(\rho_{f_k}-\rho_{f_j})-f_{\theta_e}(\rho_{f_k}+\rho_{f_j}).
\end{equation}
and the nonlinear equations for the variables $\rho_f$ are
\begin{equation}
  \label{eq:euler-lag_hyp}
  \Phi_f-2\;
  \sum_{\makebox[0pt][r]{\scriptsize$f\circ$}
    \edge\makebox[0pt][l]{\scriptsize$\circ f_k$}}\;
  \big(f_{\theta_e} (\rho_{f_k}-\rho_{f})
  -f_{\theta_e}(\rho_{f_k}+\rho_f)\big)=0,
\end{equation}
where the sum is taken over all interior edges $e$ between the face $f$ and
its neighbors $f_k$.

\begin{lemma}
  \label{lem:eul_lag_eqns_hyp}
  If $\rho\in\R^F$ is a solution of the equations \eqref{eq:euler-lag_hyp},
  then $\rho_f<0$ for all $f\in F$. The solutions of the equations
  \eqref{eq:euler-lag_hyp} are therefore in one-to-one correspondence with
  the solutions of the hyperbolic circle pattern problem.
\end{lemma}
\begin{proof}
  Since, by equation~\eqref{eq:f_prime}, the function $f_{\theta}(x)$ is
  strictly increasing for $0<\theta<\pi$,
  \begin{equation*}
    f_{\theta_e} (\rho_{f_k}+\rho_{f})
    -f_{\theta_e}(\rho_{f_k}-\rho_f) \geq 0,\quad \text{if}\quad\rho_f\geq 0.
  \end{equation*}
  Since $\Phi_f>0$ by assumption, the left hand side of equation
  \eqref{eq:euler-lag_hyp} is positive if $\rho_f\geq 0$.
\end{proof}

\begin{definition}
  \label{def:hyp_functional}
  The {\em hyperbolic circle pattern functional}\/ is the function
  \begin{equation*}
    \nonumber
    \Shyp: \R^F\longrightarrow \R
  \end{equation*}
  \begin{equation}
    \begin{split}\label{eq:Shyp}
      \Shyp(\rho)= \sum_{f_j\circ\edge\circ f_k} \Big( &
      \im\Li\big(e^{\rho_{f_k}-\rho_{f_j}+i\theta_e}\big)
      +\im\Li\big(e^{\rho_{f_j}-\rho_{f_k}+i\theta_e}\big)\\
      &+\im\Li\big(e^{\rho_{f_j}+\rho_{f_k}+i\theta_e}\big)
      +\im\Li\big(e^{-\rho_{f_j}-\rho_{f_k}+i\theta_e}\big)
      \Big)\\
      +\sum_{\circ f} &\Phi_f\rho_f.
    \end{split}
  \end{equation}
  The first sum is taken over all interior edges $e$, and $f_j$ and $f_k$ are
  the faces on either side of $e$. (The terms are symmetric in $f_j$ and
  $f_k$, so it does not matter which face is considered as $f_j$ and which as
  $f_k$.) The second sum is taken over all faces $f$.
\end{definition}

\begin{lemma}
  \label{lem:func_hyp}
  A function $\rho\in\R^F$ is a critical point of the hyperbolic circle pattern
  functional $\Shyp$, if and only if $\rho$ satisfies
  equations~\eqref{eq:euler-lag_hyp}. In that case, $\rho$ is negative. The
  critical points of $\Shyp$ are therefore in one-to-one correspondence with
  the solutions of the hyperbolic circle pattern problem.
\end{lemma}

\begin{proof}
  Similarly as in the euclidean case, one finds that 
  \begin{equation}
    \label{eq:partial_Shyp}
    \frac{\partial\Shyp}{\partial\rho_f}=
    \Phi_f-2\;
    \sum_{\makebox[0pt][r]{\scriptsize$f\,\circ$}
      \edge\makebox[0pt][l]{\scriptsize$\circ f_k$}}\;
    \big(f_{\theta_e} (\rho_{f_k}-\rho_{f})
    -f_{\theta_e}(\rho_k+\rho_f)\big)\,,
  \end{equation}
  such that $d\Shyp=0$, if and only if equations~\eqref{eq:euler-lag_hyp} are
  satisfied. By lemma~\ref{lem:eul_lag_eqns_hyp}, $\rho<0$ follows.
\end{proof}
For future reference, we note that equation~\eqref{eq:partial_Shyp} and the
proof of lemma~\ref{lem:eul_lag_eqns_hyp} imply
\begin{equation}
  \label{eq:rho_geq_0}
  \frac{\partial\Shyp}{\partial\rho_f}>0,\quad\text{if}\quad\rho_f\geq0.
\end{equation}

\section[Convexity and uniqueness]%
{Convexity of the euclidean and hyperbolic functionals. Proof of
  the uniqueness claims of theorem~$\text{\ref{thm:fundamental}}$}
\label{sec:convexity}

\begin{lemma}
  \label{lem:scale_inv}
  If a euclidean circle pattern with data $\Sigma$, $\theta$, $\Phi$
  exists, then the euclidean functional is\/ {\em{}scale-invariant}:
  Multiplying all radii $r$ with the same positive factor (equivalently,
  adding the same constant to all $\rho$) does not change its value.
\end{lemma}

\begin{proof}
  Let $1_F\in\R^F$ be the function which is $1$ on every face $f\in
  F$. Equation~\eqref{eq:SEuc} implies 
  \begin{equation*}
    \Seuc(\rho+h\,1_F)=\Seuc(\rho)
    +h \Big( \sum_{f\in F} \Phi_f - 2 \sum_{e\in\Eint} (\pi-\theta_e) \Big),
  \end{equation*}
  where $\Eint$ is the set of interior edges.
  Clearly, the functional can have a critical point only if the coefficient of
  $h$ vanishes. In this case, the functional is scale invariant.
\end{proof}

If the euclidean functional is scale invariant, one may restrict the
search for critical points to the subspace
\begin{equation}\label{eq:U}
U=\{\rho\in\mathbb{R}^F|\sum_{f\in F}\rho_f=0\}.
\end{equation}

\begin{lemma}
  The euclidean functional $\Seuc$ is strictly convex on the subspace $U$. The
  hyperbolic functional $\Shyp$ is strictly convex on $\R^F$.
\end{lemma}
\begin{proof}
By a straightforward calculation, one finds that the second derivative of the
euclidean functional is the quadratic form
\begin{equation*}
\Seuc''=\sum_{f_j\circ\edge\circ f_k}
\frac{\sin\theta_e}
     {\cosh(\rho_{f_k}-\rho_{f_j})-\cos\theta_e}
\,(d\rho_{f_k}-d\rho_{f_j})^2,
\end{equation*}
where the sum is taken over all interior edges $e$, and $f_j$ and $f_k$ are
the faces on either side.
Since it is (quietly) assumed that the surface is connected, the second
derivative is positive unless $d\rho_{f_j}=d\rho_{f_k}$ for all $f_j, f_k\in
F$. Hence it is positive definite on $U$.

For the hyperbolic functional, one obtains
\begin{multline*}
\Shyp''=\sum_{f_j\circ\edge\circ f_k}\left(
\frac{\sin\theta_e}
     {\cosh(\rho_{f_k}-\rho_{f_j})-\cos\theta_e}
\,(d\rho_{f_j}-d\rho_{f_k}\big)^2 +
\right.
\\
\left.
\frac{\sin\theta_e}
     {\cosh(\rho_{f_j}+\rho_{f_k})-\cos\theta_e}
\,(d\rho_{f_j}+d\rho_{f_k}\big)^2
\right),
\end{multline*}
which is positive definite on $\mathbb{R}^F$.
\end{proof}

This proves the uniqueness claims of theorem~\ref{thm:fundamental}.

\section[The spherical functional]{The spherical circle pattern functional}
\label{sec:sph_func}

Like in the euclidean and hyperbolic cases, there is a functional for
spherical circle patterns whose critical points correspond to solutions of
the circle pattern problem.

\begin{lemma}
  \label{lem:phi_sph}
  Suppose $r_1$ and $r_2$ are two sides of a spherical triangle (with
  $0<r_{1,2}<\pi$), the included angle between them is $\theta$, and the
  remaining angles are $\varphi_1$ and $\varphi_2$, as shown in
  figure~\ref{fig:phiRhoTriangle}.  Then
  \begin{equation}
    \label{eq:phi_of_rho_sph}
    \varphi_{1}=f_{\theta}(\rho_{2}-\rho_{1})
    +f_{\pi-\theta}(\rho_{2}+\rho_{1}),
  \end{equation}
  where $f_\theta(x)$ is defined by equation~\eqref{eq:f}, and
  \begin{equation}
    \label{eq:rho_sph}
    \rho = \log\tan\frac{r}{2}.
  \end{equation}
  The inverse relation is
  \begin{equation}\label{eq:rho_of_phi_sph}
    \rho_{1}=\frac{1}{2}\log
    \frac{\sin\Big(
      \frac{\displaystyle -\theta^*+\varphi_{1}+\varphi_{2}}
      {\displaystyle 2}\Big)
      \sin\Big(
      \frac{\displaystyle \theta^*-\varphi_{1}+\varphi_{2}}
      {\displaystyle 2}\Big)}
    {\sin\Big(
      \frac{\displaystyle \theta^*+\varphi_{1}+\varphi_{2}}
      {\displaystyle 2}\Big)
      \sin\Big(
      \frac{\displaystyle \theta^*+\varphi_{1}-\varphi_{2}}
      {\displaystyle 2}\Big)},
  \end{equation}
  where $\theta^*=\pi-\theta$,\, $\varphi_{1,2}>0$, and\,
  $\theta^*<\varphi_1+\varphi_2<2\pi-\theta^*$.
\end{lemma}

Equation~\eqref{eq:phi_of_rho_sph} is derived in appendix~\ref{app:trig}.
Equation~\eqref{eq:rho_of_phi_sph} follows from~\eqref{eq:phi_of_rho_sph}
(and the corresponding equation for $\varphi_2$) by a straightforward
calculation using equations~\eqref{eq:f_inverse} and \eqref{eq:f_symmetry}.

To construct the spherical circle pattern functional, one proceeds like in
the euclidean and hyperbolic cases
(sections~\ref{sec:analytic}--\ref{sec:hyp_func}). In this case, the new
variables $\rho_f$ are given by equation~\eqref{eq:rho_sph}. There is a
one-to-one correspondence between radii $r$ with $0<r<\pi$, and $\rho\in\R$.
Instead of equations~\eqref{eq:phi_e_of_rho_f_euc}
or~\eqref{eq:phi_e_of_rho_f_hyp}, one has
\begin{equation}
  \label{eq:phi_e_of_rho_f_sph}
  \varphi_{\vece}=f_{\theta_e}(\rho_{f_k}-\rho_{f_j})
  +f_{\theta^*_e}(\rho_{f_k}+\rho_{f_j}),
\end{equation}
where $\theta^*=\pi-\theta$. Consequently, the nonlinear equations for the
variables $\rho_f$ are
\begin{equation}
  \label{eq:euler-lag_sph}
  \Phi_f-2\;
    \sum_{\makebox[0pt][r]{\scriptsize$f\,\circ$}
      \edge\makebox[0pt][l]{\scriptsize$\circ f_k$}}\;
    \big(f_{\theta_e} (\rho_{f_k}-\rho_{f})
    +f_{\theta^*_e}(\rho_{f_k}+\rho_f)\big)=0,
\end{equation}
where the sum is taken over all interior edges $e$ around $f$, and $f_k$ is
the face on the other side of $e$.

\begin{definition}
  \label{def:sph_functional}
  The {\em spherical circle pattern functional}\/ is the function
  \begin{equation*}
    \nonumber
    \Ssph: \R^F\longrightarrow \R
  \end{equation*}
  \begin{equation}
    \begin{split}\label{eq:Ssph}
      \Ssph(\rho)= \sum_{f_j\circ\edge\circ f_k} \Big( &
      \im\Li\big(e^{\rho_{f_k}-\rho_{f_j}+i\theta_e}\big)
      +\im\Li\big(e^{\rho_{f_j}-\rho_{f_k}+i\theta_e}\big)     \\
      &-\im\Li\big(e^{\rho_{f_j}+\rho_{f_k}+i(\pi-\theta_e)}\big)
      -\im\Li\big(e^{-\rho_{f_j}-\rho_{f_k}+i(\pi-\theta_e)}\big)  \\
      & - \pi(\rho_{f_j} + \rho_{f_k})
      \Big)\\
      +\sum_{\circ f} &\Phi_f\rho_f.
    \end{split}
  \end{equation}
  The first sum is taken over all interior edges $e$, and $f_j$ and $f_k$ are
  the faces on either side of $e$. (The terms are symmetric in $f_j$ and
  $f_k$, so it does not matter which face is considered as $f_j$ and which as
  $f_k$.) The second sum is taken over all faces $f$.
\end{definition}

\begin{lemma}
  \label{lem:func_sph}
  A function $\rho\in\R^F$ is a critical point of the spherical circle
  pattern functional $\Ssph$, if and only if $\rho$ satisfies
  equations~\eqref{eq:euler-lag_sph}. The critical points of $\Ssph$ are
  therefore in one-to-one correspondence with the solutions of the spherical
  circle pattern problem.
\end{lemma}

This proposition is proved in the same way as the lemmas~\ref{lem:func_euc}
and~\ref{lem:func_hyp}.

The spherical circle pattern functional is not convex. By a
straightforward calculation, one finds that
\begin{multline}
    \Ssph''= \sum_{f_j\circ\edge\circ f_k}\left(
\frac{\sin\theta_e}
     {\cosh(\rho_{f_k}-\rho_{f_j})-\cos\theta_e}
\,(d\rho_{f_j}-d\rho_{f_k}\big)^2
\right.
\\
\left.
-\frac{\sin\theta_e}
     {\cosh(\rho_{f_j}+\rho_{f_k})+\cos\theta_e}
\,(d\rho_{f_j}+d\rho_{f_k}\big)^2
\right).
\end{multline}
This quadratic form is negative for the tangent vector $1_F\in\R^F$,
which has a $1$ in every component. Hence, the negative index is at least
one. 

\begin{remark}
  This has a geometric explanation. Consider a circle pattern in the sphere.
  Focus on one flower: a central circle with its neighbors. The neighbors
  nicely fit around the central circle. Now decrease the radii of all circles
  by the same factor. The effect is the same as increasing the radius of the
  sphere by that factor. This makes the sphere flatter. The neighbors will
  not fit around the central circle anymore, but there will be a gap. To
  adjust the radius of the central circle so that the neighbors fit around,
  one would make it even smaller.
\end{remark}

Nevertheless, in numerical experiments, the spherical functional has been
used with amazing success to construct circle patterns in the sphere. This
following method was used. Consider the reduced functional
\begin{equation}
  \label{eq:SsphTilde}
  \SsphTilde(\rho)=\max_{t\in\R}\, \Ssph(\rho + t\, 1_F).
\end{equation}
Clearly, $\SsphTilde(\rho)$ is invariant under a shift
$\rho\mapsto\rho+t 1_F$. To solve a spherical circle pattern problem,
minimize $\SsphTilde(\rho)$ on the set $U$ defined by equation~\eqref{eq:U}.
Equivalently, minimize $\SsphTilde(\rho)$ under the constraint $\rho_f=0$ for
some fixed $f\in F$. 

The numerical evidence suggests that this method works whenever the circle
pattern in question exists. In particular, this method can be used to
construct branched circle patterns in the sphere. In general, such patterns
cannot be constructed using the euclidean functional after a stereographic
projection.

We return to the definition of the reduced functional $\SsphTilde(\rho)$,
equation~\eqref{eq:SsphTilde}.  The following proposition provides a
geometric interpretation for the condition
\begin{equation*}
  \frac{d}{dt}\,d\Ssph(\rho+t\,{1}_F)\Big|_{t=0}=0.
\end{equation*}
Consider the spherical circle pattern problem with data $\Sigma$, $\theta\in
(0,\pi)^E$, and $\Phi\in (0,\infty)^F$. For the sake of simplicity, assume
that $\Sigma$ is a decomposition of a surface without boundary. Suppose there
exists a corresponding circle pattern in a surface $M$ with constant
curvature $1$ and cone-like singularities. Then, due to the Gauss-Bonnet{}
theorem, the total area $A$ of $M$ is determined by the formula: 
\begin{equation}
  \label{eq:gauss_bonnet_sph}
  A + K_v + K_f = 2\pi(|F|-|E|+|V|),
\end{equation}
where $K_v$ and $K_f$ are as in equations~\eqref{eq:curvature_v}
and~\eqref{eq:curvature_f}, and $|F|$, $|E|$, and $|V|$ are the numbers of
faces, edges and vertices of $\Sigma$, respectively.

\begin{proposition}
  Let $\rho\in \R^F$ be arbitrary. Given a cell decomposition $\Sigma$ of a
  closed compact surface, the intersection angles $\theta$, and radii
  $r=2\arctan e^{\rho}$, construct the surface $M^{(\rho)}$ with constant
  curvature $1$ and cone-like singularities by gluing together spherical
  kite-shaped quadrilaterals, as was described (for the euclidean case) in
  the proof of lemma~\ref{lem:given_theta_and_r}.  Let $A^{(\rho)}$ be the
  total area of $M^{(\rho)}$. Then
  \begin{equation}
    \label{eq:Ssph_scale}
    \frac{d}{dt}\,d\Ssph(\rho+t\,{1}_F)\Big|_{t=0} = A - A^{(\rho)}
  \end{equation}
  where $A$ is determined by equation \eqref{eq:gauss_bonnet_sph}.
\end{proposition}

\begin{proof}
  Let $A^{(\rho)}_e$ be the area of the quadrilateral corresponding to an
  edge $e$ as in figure~\ref{fig:kite_f_theta}~{\em (left)}. Then
  \begin{equation*}
    \begin{split}
      A^{(\rho)}_e &\;\;= 2\varphi + 2\varphi' + 2\theta_e - 2\pi \\
      &\overset{\eqref{eq:phi_e_of_rho_f_sph}}{=}
      2f_{\theta_e}(\rho_k-\rho_j)+2f_{\theta_e}(\rho_j-\rho_k)
      +4f_{(\pi-\theta_e)}(\rho_j+\rho_k)+2\theta_e - 2\pi \\
      &\;\overset{\eqref{eq:f_symmetry}}{=} 
      4 f_{(\pi-\theta_e)}(\rho_j+\rho_k).
    \end{split}
  \end{equation*}
  Also, 
  \begin{equation*}
    \sum_{e\in E}2\theta_e=2\pi|V|-\sum_{v\in V}K_v
  \end{equation*}
  and
  \begin{equation*}
    \sum_{f\in F}\Phi_f=2\pi|F|-\sum_{f\in F}K_f.
  \end{equation*}
  Finally,
  \begin{equation*}
    \begin{split}
      &\frac{d}{dt}\,d\Ssph(\rho+t\mathbf{1}_F)
      = \Big(\sum_{f\in F}\frac{d}{d\rho_f}\Big)\,\Ssph(\rho)\\
      &\quad
      \overset{\genfrac{}{}{0pt}{}{\eqref{eq:Ssph}}{\eqref{eq:f_integral}}}{=}
      \sum_{f_j\circ\edge\circ f_k}\big(-2f_{(\pi-\theta_e)}(\rho_j+\rho_k)+
      2f_{(\pi-\theta_e)}(-\rho_j-\rho_k) -2\pi\big)+\sum_{\circ f} \Phi_f \\
      &\quad\,\overset{\eqref{eq:f_symmetry}}{=}\, \sum_{f_j\circ\edge\circ
        f_k} \big(-4f_{(\pi-\theta_e)}(\rho_j+\rho_k)+2\theta_e -2\pi\big)
      +\sum_{\circ f} \Phi_f
    \end{split}
  \end{equation*}
  Equation~\eqref{eq:gauss_bonnet_sph} follows. 
\end{proof}

\section[Coherent angle systems. Existence of circle patterns]
{Coherent angle systems. The existence of circle patterns}
\label{sec:coherent_angle_systems}

In section~\ref{sec:convexity}, the uniqueness of a circle pattern was
deduced from the convexity of the euclidean and hyperbolic functionals. This
section and the next one are devoted to the existence part of
theorem~\ref{thm:fundamental}.  To establish that the euclidean functional
attains a minimum, we will show that 
\begin{equation*}
  \Seuc(\rho)\longrightarrow\infty\quad\text{if}\quad
  \|\rho\|\longrightarrow\infty\text{ in }U,
\end{equation*}
where $U\in\R^F$ is the subspace defined by equation~\eqref{eq:U}. 
To establish that the hyperbolic functional attains a minimum, we will to
show that 
\begin{equation*}
  \Shyp(\rho)\longrightarrow\infty\quad\text{if}\quad
  \|\rho\|\longrightarrow\infty\text{ with }\rho<0.
\end{equation*}
Because of the inequality~\eqref{eq:rho_geq_0}, this suffices.

To estimate the functionals from below, one has to compare the sum over
interior edges with the sum over faces in equations (\ref{eq:SEuc}) and
(\ref{eq:Shyp}).  This is achieved with the help of a so called `coherent
angle system'. In this section, we prove that the functionals have minima if
and only if coherent angle systems exist. In
section~\ref{sec:proof_thm:fundamental}, we will show that the conditions of
theorem~\ref{thm:fundamental} are necessary and sufficient for the existence
of a coherent angle system.

Coherent angle systems also play an important role in
chapter~\ref{cha:other}, when we derive other variational principles by
Legendre transformations. Spherical coherent angle systems are defined below,
but not used until chapter~\ref{cha:other}.

Let $\vecE$ be the set of oriented edges. For an oriented edge
${\vece}\in\vecE$, denote by $-\vece\in\vecE$ the edge with the opposite
orientation, and by $e$ the corresponding non-oriented edge.

\begin{definition}
  A {\em euclidean coherent angle system}\/ is a function
  $\varphi\in\mathbb{R}^{\vecEint}$ on the set $\vecEint$ of interior
  oriented edges which satisfies the following two conditions.

  {\it (i)}\/ For all oriented edges $\vece\in\vecEint$,
    \begin{equation*}
      \varphi_{\vece}>0\quad\text{and}\quad\varphi_{\vece}+\varphi_{-\vece}=\pi-\theta_e.
    \end{equation*}

    {\it (ii)}\/ For all faces $f\in F$,
    $$
    \sum 2\varphi_{\vece}=\Phi_f, $$
    where the sum is taken over all
    oriented interior edges $\vece$\/ in the oriented boundary of $f$.

    A {\em hyperbolic coherent angle system}\/ satisfies 

    {\it ($\text{i}\,'$)}\/ For all oriented edges $e\in\vecE$,
    \begin{equation*}
      \varphi_{\vece}>0\quad\text{and}\quad
      \varphi_{\vece}+\varphi_{-\vece}<\pi-\theta_e.
    \end{equation*}
    
    \noindent{}
    and condition {\it (ii)}\/ above.

  A {\em spherical coherent angle system}\/ satisfies 

    {\it ($\text{i}\,''$)}\/ For all oriented edges $e\in\vecE$,
    \begin{equation*}
      0<\varphi_{\vece}<\pi,\quad
      \pi-\theta_e<\varphi_{\vece}+\varphi_{-\vece}<\pi+\theta_e,
      \quad\text{and}\quad 
      \big|\varphi_{\vece}-\varphi_{-\vece}\big|<\pi-\theta_e,
    \end{equation*}

    \noindent{}
    and condition {\it (ii)}\/ above.
\end{definition}

(Note that the exterior angles of a spherical triangle satisfy the triangle
inequalities.)  The following lemma reduces the question of existence of a
(euclidean or hyperbolic) circle pattern to the question of existence of a
coherent angle system.

\begin{lemma}\label{lem:coh_angles}
The functional $\Seuc$ ($\Shyp$) has a critical point, if and only if a
euclidean (hyperbolic) coherent angle system exists.
\end{lemma}

\begin{proof}
  If the functional $\Seuc$ ($\Shyp$) has a critical point $\rho$, then
  equation \eqref{eq:phi_e_of_rho_f_euc} (equation
  \eqref{eq:phi_e_of_rho_f_hyp}) yields a coherent angle system. It is left
  to show that, conversely, the existence of a coherent angle system implies
  the existence of a critical point.
  
  Consider the euclidean case. Suppose a euclidean coherent angle system
  $\varphi$ exists. This implies
\begin{equation*}
  \sum_{f\in F}\Phi_f=2\sum_{e\in \Eint}(\pi-\theta_e).
\end{equation*}
Hence, the functional $\Seuc$ is scale invariant. (See the proof of
lemma~\ref{lem:scale_inv}.) We will show that $\Seuc(\rho)\rightarrow\infty$
if $\rho\rightarrow\infty$ in the subspace $U$ defined in
equation~\ref{eq:U}. More precisely, we will show that for $\rho \in U$,
\begin{equation}\label{eq:SEucEstimate}
\Seuc(\rho)>2\min_{\vece\in\vecEint}\varphi_{\vece}\;\max_{f\in F}\,|\rho_f|.
\end{equation}
The functional $\Seuc$ must therefore attain a minimum, which is a critical
point.

For $x\in\mathbb{R}$ and $0<\theta<\pi$,
\begin{equation*}
\im\Li(e^{x+i\theta})+\im\Li(e^{-x+i\theta})>(\pi-\theta)\,|x|,
\end{equation*}
and hence, by equation~\eqref{eq:SEuc}, 
\begin{equation*}
\Seuc(\rho)>-2\sum_{e\in E}
(\pi-\theta_e)\min\big(\rho_{f_k},\rho_{f_j}\big) 
+ \sum_{f\in F}\Phi_f\rho_f,
\end{equation*}
where the sum is taken over the unoriented interior edges $e$, and $f_j$ and
$f_k$ are the faces on either side of $e$. Now, we use the coherent angle
system $\varphi$ to merge the two sums. Because
\begin{equation*}
\sum_{f\in F}\Phi_f\rho_f 
= 2\sum_{e\in\Eint} (\varphi_{\vece}\,\rho_{f_j} 
+ \varphi_{-\vece}\,\rho_{f_k}), 
\end{equation*}
one obtains
\begin{equation*}
\Seuc(\rho)>2\sum_{e\in\Eint}
\min\big(\varphi_{\vece},\varphi_{-\vece}\big)\,
\big|\rho_{f_k}-\rho_{f_j}\big|.
\end{equation*}
Since we assume the cellular surface to be connected, we get
\begin{equation*}
\Seuc(\rho)>2\min_{\vece\in\vecEint}\varphi_{\vece}\,
\big(\max_{f\in F}\rho_f-\min_{f\in F}\rho_f\big),
\end{equation*}
and from this the estimate (\ref{eq:SEucEstimate}).

The hyperbolic case is similar. One shows that, if all $\rho_f<0$,
\begin{equation*}
\Shyp(\rho)>2
\min_{e\in\Eint}\big|\varphi_{\vece}+\varphi_{-\vece}-(\pi-\theta_e)\big|\;
\max_{f\in F}\big|\rho_f\big|.
\end{equation*}
\end{proof}

\section{Conclusion of the proof of theorem~$\text{\ref{thm:fundamental}}$}
\label{sec:proof_thm:fundamental}

With this section we complete the proof of theorem~\ref{thm:fundamental}.  In
section~\ref{sec:convexity}, we have shown the uniqueness claim. By
lemma~\ref{lem:coh_angles}, the circle patterns exist, if and only if a
coherent angle system exists.  All that is left to show is the following
lemma.
\begin{lemma}
\label{lem:exist}
A euclidean/hyperbolic coherent angle system exists if and only if the
conditions of theorem~\ref{thm:fundamental} hold.
\end{lemma}

The rest of this section is devoted to the proof of lemma~\ref{lem:exist}.
It is easy to see that these conditions are necessary. To prove that they are
sufficient, we apply the feasible flow theorem of network theory. Let $(N,X)$
be a network (i.e.\ a directed graph), where $N$ is the set of nodes and $X$
is the set of branches. For any subset $N'\subset N$ let $\outof(N')$ be
the set of branches having their initial node in $N'$ but not their terminal
node. Let $\into(N')$ be the set of branches having their terminal node in
$N'$ but not their initial node. Assume that there is a lower capacity bound
$a_x$ and an upper capacity bound $b_x$ associated with each branch $x$,
with $-\infty\leq a_x\leq b_x\leq \infty$. 

\begin{definition}
A {\em feasible flow}\/ is a function
$\varphi\in\mathbb{R}^X$, such that Kirchoff's current law is
satisfied, that is, for each $n\in N$,
\begin{equation*}
\sum_{x\in \outof(\{n\})} \varphi_x = \sum_{x\in \into(\{n\})} \varphi_x,
\end{equation*}
and $a_x\leq\varphi_x\leq b_x$ for all branches $x\in X$.
\end{definition}

\begin{flowThm}
A feasible flow exists if and only if for every nonempty subset $N'\subset
N$ of nodes with $N'\not=N$, 
\begin{equation*}
\sum_{x\in \outof(N')} b_x \geq \sum_{x\in \into(N')} a_x.
\end{equation*}
\end{flowThm}

A proof is given by Ford and Fulkerson \cite[ch.~II, \S3]{ford_fulkerson62}.
(Ford and Fulkerson assume the capacity bounds to be non-negative, but this
is not essential.)

\begin{figure}%
\begin{center}%
\input{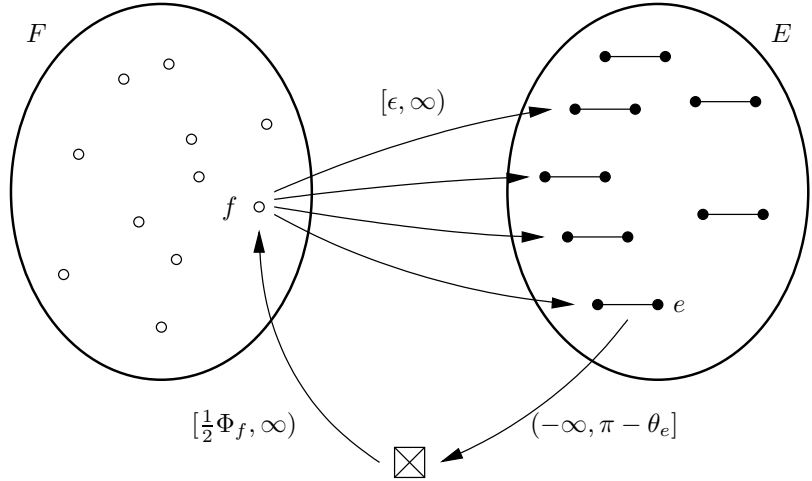}%
\caption{The network $(N,X)$. Only a few of the branches and capacity
  intervals are shown.}%
\label{fig:network}%
\end{center}%
\end{figure}

To prove lemma~\ref{lem:exist} in the euclidean case, consider the following
network; see figure~\ref{fig:network}.
The nodes are all faces and non-oriented interior edges of the cellular
surface, and one further node that we denote by $\boxtimes$: $N=F\cup E\cup
\{\boxtimes\}$. There is a branch in $X$ going from $\boxtimes$ to each face
$f\in F$ with capacity interval $[\frac{1}{2}\Phi_f,\infty)$. From each face
$f$ there is a branch in $X$ going to the non-oriented interior edges of the
boundary of $f$ with capacity interval $[\epsilon,\infty)$, where
$\epsilon>0$ will be determined later. Finally there is a branch in $X$ going
from each non-oriented edge $e \in E$ to $\boxtimes$ with capacity
$(-\infty,\pi-\theta_e]$.

Assume the conditions of theorem~\ref{thm:fundamental} are fulfilled.  A
feasible flow in the network yields a coherent angle system. Indeed, since we
have equality in~\eqref{eq:phi_theta_sum_inequality} if $F'=F$,
Kirchoff's current law at $\boxtimes$ implies that the flow into each face
$f$ is $\frac{1}{2}\Phi_f$ and the flow out of each edge $e$ is
$\pi-\theta_e$. It follows that the flow in the branches from $F$ to $E$
constitutes a coherent angle system.

We
need to show that the condition of the feasible flow theorem is
satisfied. Suppose $N'$ is a nonempty proper subset of $N$. Let $F'=N'\cap
F$ and $E'=N'\cap\Eint$.

Consider first the case that $\boxtimes\in N'$, which is the easy one. Since
$N'$ is a proper subset of $N$ there is a face $f\in F$ or an edge $e\in E$
which is not in $N'$. In the first case there is a branch out of $N'$ with
infinite upper capacity bound. In the second case there is a branch into $N'$
with negative infinite lower capacity bound. Either way, the condition of the
feasible flow theorem is trivially fulfilled.

Now consider the case that $\boxtimes\not\in N'$. We may assume that for each
face $f\in F'$, the interior edges in the boundary of $E$ are contained in
$E'$.  Otherwise, there would be branches out of $N'$ with infinite upper
capacity bound.  For subsets $A,B\subset N$ denote by $A\rightarrow B$ the
set of branches in $X$ having initial node in $A$ and terminal node in $B$.
Then the condition of the feasible flow theorem is equivalent to
\begin{equation*}
\sum_{f\in F'}\frac{1}{2}\Phi_f + \epsilon\,|F\setminus F'\rightarrow E'|
\leq 
\sum_{e\in E'} \theta^*_e.
\end{equation*}
It is fulfilled if we choose
\begin{equation}\label{eq:epsilon}
\epsilon<\,\frac{1}{2|E|}
\min_{F'} \left(\sum_{e\in E'(F')}(\pi-\theta_e)
-\sum_{f\in F'}\frac{1}{2}\Phi_f\right),
\end{equation}
where the minimum is taken over all proper nonempty subsets $F'$ of $F$ and
$E'(F')$ is the set of all non-oriented edges incident with a face in $F'$.
The minimum is greater than zero because strict inequality holds
in~\eqref{eq:phi_theta_sum_inequality} if $F'$ is a proper subset of $F$.

In the hyperbolic case, the proof is only a little bit more complicated.  In
the network, the flow in the branches going from $\boxtimes$ to a face $f$
must constrained to be exactly $\frac{1}{2}\Phi_f$; and the capacity
interval of branches going from an edge $e$ to $\boxtimes$ has to be changed
to $(-\infty,\theta^*_e-\epsilon]$.

\section{Proof of theorem~$\text{\ref{thm:higher_genus}}$}
\label{sec:proof_thm_higher_genus}

\begin{figure}%
\begin{center}%
\input{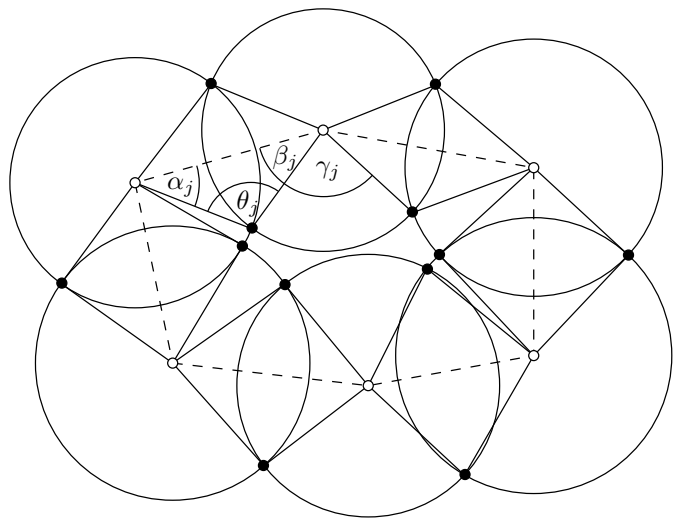}%
\end{center}%
\caption{The circles and quadrilaterals corresponding to a closed path in the
  Poincar\'e dual {\em{}(dashed)}, which cuts out a disc.}%
\label{fig:copath}%
\end{figure}
The necessity of the condition of theorem~\ref{thm:higher_genus} can be seen
from the following geometrical argument. Figure~\ref{fig:copath} shows a
closed path in the Poincar\'e dual cell decomposition $\Sigma^*$, which cuts
out a disc. By the Gauss-Bonnet theorem,
\begin{equation*}
  \sum_j (\pi-\alpha_j-\beta_j-\gamma_j)=2\pi-\epsilon A,
\end{equation*}
where $\epsilon\in\{0,-1\}$ is the curvature of the surface and $A$ is the
area enclosed by the path. But
\begin{equation*}
  \pi-\alpha_j-\beta_j = \theta_j - \epsilon A_j,
\end{equation*}
where $A_j$ is the are of triangle $j$. Hence,
\begin{equation}
  \label{eq:theta_gamma_sum}
  \sum_j\theta_j = 2\pi + \sum_j \gamma_j - \epsilon(A-\sum_j A_j).
\end{equation}
Inequality~\eqref{eq:theta_sum_inequality} follows, with equality holding if
and only if all $\gamma_j$ are zero and $\epsilon(A-\sum A_j)=0$. This is the
case if and only if the path contains only a single vertex of $\Sigma$ in its
interior.

We will now prove the sufficiency of the condition. The principle tool in the
proof is the following lemma.

\begin{lemma}
  \label{lem:general_euler}
  Suppose a nonempty $1$-dimensional cell complex (a graph)\/ $\Gamma=(V_1,E_1)$ is
  cellularly embedded in a cell decomposition $\Sigma=(V, E, F)$ of a closed
  compact surface. Suppose\/ $\Gamma$ separates\/ $\Sigma$ into $r$
  contiguous regions. Then 
  \begin{equation}
    \label{eq:general_euler}
    r-|E_1|+|V_1|=|F|-|E|+|V|+\sum_{j=1}^r h_j,
  \end{equation}
  where $h_j$ is the dimension of the first $\Z_2$-homology group of the
  $j^{\text{\it{th}}}$ region.
\end{lemma}

If all regions are simply connected, then all $h_j=0$, and
equation~\eqref{eq:general_euler} states that the Euler characteristic is
invariant under subdivisions. A multiply connected region can be turned into
a simply connected one by adding $h_j$ cuts to $\Gamma$ which do not separate
the region. Each cut reduces the left hand side of
equation~\eqref{eq:general_euler} by one.  For a more formal proof of
lemma~\ref{lem:general_euler}, see Giblin \cite[ch.~9]{giblin77} or
appendix~\ref{app:combi_top}.

Suppose the condition of theorem~\ref{thm:higher_genus} holds. We will deduce
the condition of theorem~\ref{thm:fundamental}. Thus, let $F'\subseteq F$\/
be a nonempty subset of the set $F$ of faces, and let $E'\subseteq E$ be the
set of edges, which are incident with any face in $F'$. We have to show
inequality~\eqref{eq:phi_theta_sum_inequality} where all $\Phi_f=2\pi$. That
is, we have to show
\begin{equation}
  \label{eq:to_show_1.7}
  \sum_{e\in E'}2\theta^*_e \geq 2\pi|F'|,
\end{equation}
where equality holds if and only if the genus $g=1$ and $F'=F$. 

Let 
\begin{equation*}
  F''=F\setminus F' \quad\text{and}\quad E''=E\setminus E'
\end{equation*}
be the complements of $F'$ and $E'$. Identify the faces, edges, and vertices
of the cell decomposition $\Sigma$ with the corresponding vertices, edges,
and faces of its Poincar\'e dual $\Sigma^*$. Consider the $1$-dimensional
subcomplex $\Gamma=(F'',E'')$ of $\Sigma^*$ which consists of the vertex set
$F''$ and edge set $E''$. (It is easy to see that it is indeed a subcomplex.)
First, suppose that $F'\not=F$. Then the graph $\Gamma$ is not empty. By
lemma~\ref{lem:general_euler}, applied to $\Gamma$ in $\Sigma^*$,
\begin{equation*}
r-|E''|+|F''|=|V|-|E|+|F|+\sum_{j=1}^{r}h_j,
\end{equation*}
or, equivalently,
\begin{equation}
  \label{eq:generalizedEulerApplied}
  |F'|-|E'|+|V|=\sum_{j=1}^{r}(1-h_j).
\end{equation}
Now,
\begin{equation*}
  \begin{split}
    \sum_{e\in E'}2\theta^*_e &= 2\pi|E'|- \sum_{e\in E'}2\theta_e. \\
    &= 2\pi|E'| - \sum_{e\in E}2\theta_e + \sum_{e\in E''}2\theta_e. \\
  \end{split}
\end{equation*}
Because, by the condition of theorem~\ref{thm:higher_genus}, $\theta$ sums to
$2\pi$ around each vertex,
\begin{equation*}
  \sum_{e\in E}2\theta_e = 2\pi|V|,
\end{equation*}
and hence
\begin{equation}
  \label{eq:sum_theta_star}
  \sum_{e\in E'}2\theta^*_e = 2\pi\big(|E'| - |V|\big) 
  + \sum_{e\in E''}2\theta_e. 
\end{equation}
With equation~\eqref{eq:generalizedEulerApplied}, we get
\begin{equation}
  \label{eq:sum2theta-2piFprime}
  \sum_{e\in E'}2\theta^*_e - 2\pi|F'| = 
  \sum_{e\in E''}2\theta_e - 2\pi\sum_{j=1}^{r}(1-h_j).
\end{equation}
Let $e^j_1,\ldots,e^j_{b_j}$ be the boundary path of the $j^{\text{\it{}th}}$
region into which $\Sigma^*$ decomposes if cut along the edges in $E''$. (The
same edge may appear twice in one boundary path.) Each edge of $E''$ appears
exactly twice in any of the boundary paths. Therefore,
\begin{equation*}
  \sum_{e\in E''}2\theta_e = \sum_{j=1}^{r}\sum_{i=1}^{b_j}\theta_{e^j_i},
\end{equation*}
and, with equation~\eqref{eq:sum2theta-2piFprime},
\begin{equation*}
    \sum_{e\in E'}2\theta^*_e - 2\pi|F'| = 
  \sum_{j=1}^{r}\Big(\sum_{i=1}^{b_j}\theta_{e^j_i}-2\pi(1-h_j)\Big).
\end{equation*}
In any case,
\begin{equation}\label{eq:sum_theta_e_ij}
  \sum_{i=1}^{b_j}\theta_{e^j_i}-2\pi(1-h_j)>0.
\end{equation}
For if the $j^{\text{\it{}th}}$ region is not a disc, then $h_j\geq 1$, and
\eqref{eq:sum_theta_e_ij} follows because $\theta>0$. If the
$j^{\text{\it{}th}}$ region is a disc, then $h_j=0$ and
\eqref{eq:sum_theta_e_ij} follows by the condition of
theorem~\ref{thm:higher_genus}. We have shown that
inequality~\eqref{eq:to_show_1.7} holds strictly, if $F'\not=F$.

Now assume $F'=F$, and hence $E'=E$. Like equation~\eqref{eq:sum_theta_star},
we get
\begin{equation*}
    \sum_{e\in E}2\theta^*_e = 2\pi\big(|E| - |V|\big).
\end{equation*}
Therefore, since $2-2g=|F|-|E|+|V|$,
\begin{equation*}
    \sum_{e\in E}2\theta^*_e = 2\pi(|F| - 2 + 2g).
\end{equation*}
Thus, inequality~\eqref{eq:to_show_1.7} holds strictly, except if $|F'|=|F|$
and $g=1$. This completes the proof of
theorem~\ref{thm:higher_genus}.

\section{Proof of theorem~$\text{\ref{thm:rivin_circ}}$}
\label{sec:proof_thm_rivin_circ}

\begin{figure}%
\hfill%
\includegraphics{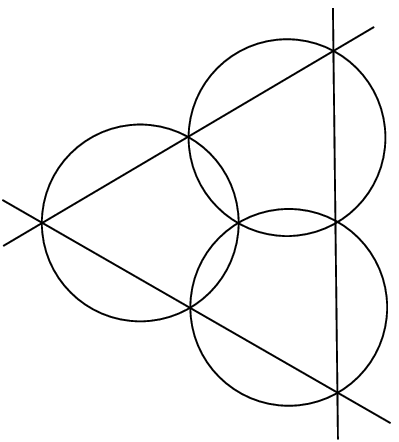}%
\hfill%
\raisebox{3mm}{\input{boundcirc.tex}}%
\hspace*{\fill}%
\caption{{\em Left:}\/ The regular cubic pattern after stereographic
  projection to the plane. {\em Right:}\/ The Neumann boundary condition
  $\Phi_f$ for the new boundary circles.}%
\label{fig:cube_boundcirc}%
\end{figure}

A circle pattern in the sphere may be projected stereographically to the
plane, choosing some vertex, $v_{\infty}$, as the center of projection. 
(For example, figure \ref{fig:cube_boundcirc} {\em{}(left)}
shows the circle pattern combinatorially equivalent to the cube and with
intersection angles $\pi/3$ after stereographic projection.) One obtains a
circle pattern in the plane, in which some circles (those corresponding to
faces incident with $v_{\infty}$) have degenerated to straight lines. Since
stereographic projection is conformal, the intersection angles are the same
as in the spherical pattern. Furthermore, M{\"o}bius-equivalent circle
patterns in the sphere correspond to patterns in the plane which differ by a
similarity transformation; that is, provided the same vertex is chosen as the
center of projection.

To prove that the condition of theorem~\ref{thm:rivin_circ} is necessary for
a circle pattern to exist, project the pattern to the plane as described
above and proceed as in section~\ref{sec:proof_thm_higher_genus}. Some
circles may have degenerated to straight lines, but
equation~\eqref{eq:theta_gamma_sum} holds nonetheless, with $\epsilon=0$.
All $\gamma_j$ are zero if the dual path has one finite vertex in its
interior, or if all circles have degenerated to straight lines. In that case
the dual path on the sphere encircles the vertex which is the center of
projection.

It is left to show that the condition of theorem~\ref{thm:rivin_circ} is
sufficient. So assume the condition holds. The idea is to show that
corresponding planar pattern exists using theorem~\ref{thm:fundamental}, and
then project it stereographically to the sphere.

To show the existence of the planar pattern, first choose a vertex
$v_{\infty}$ of the cell decomposition $\Sigma$ of the sphere. Let $F$ be the
set of faces of $\Sigma$, and let $F_{\infty}\subset F$ be the set of faces
which are incident with $v_{\infty}$. Then remove from $\Sigma$ the vertex
$v_{\infty}$, all the faces in $F_{\infty}$, and the edges between them, to
obtain a cell complex $\Sigma_0$ with face set $F_0=F\setminus F_{\infty}$.
Because $\Sigma$ is a strongly regular cell decomposition of the sphere,
$\Sigma_0$ is a cell decomposition of the closed disc.
Hence, theorem~\ref{thm:fundamental} may be applied to prove the existence
and uniqueness of the planar pattern. (Except in the trivial case when
$\Sigma_0$ has only one face.) The Neumann boundary conditions need to be
specified. For a boundary face $f$ of $\Sigma_0$, set
\begin{equation}\label{eq:PhiBoundary}
\Phi_f=2\pi-\sum 2\theta^*_e, 
\end{equation}
where the sum is taken over all boundary edges $e$ of $\Sigma_0$ which are
incident with $f$. See figure~\ref{fig:cube_boundcirc} {\em{}(right)}. For
all interior faces $f$, set
\begin{equation}\label{eq:PhiInterior} 
\Phi_f=2\pi.
\end{equation}
If the conditions of theorem \ref{thm:fundamental} are satisfied, one may
construct the corresponding planar pattern, add the lines corresponding to
the removed faces, and project to the sphere. Hence, we need to show the
statements {\em{(i)}} and {\em{(ii)}}:

\smallskip%
\makebox[1.5em][l]{\em{}(i)} For a boundary face $f$ of
$\Sigma_0$, $\Phi_f$ as defined in equation~\eqref{eq:PhiBoundary} is
positive.

\smallskip%
\makebox[1.5em][l]{\em{}(ii)} If $F_0'\in F_0$ is a nonempty set of faces of
$\Sigma_0$, and $E_0'$ is the set of of all interior edges of $\Sigma_0$
which are incident with any face in $F_0'$, then
\begin{equation}
  \label{eq:statement(ii)}
  \sum_{f\in F_0'} \Phi_f \leq \sum_{e\in E_0'} 2\theta^*_e,
\end{equation}
where equality holds if and only if $F_0'=F_0$.

\smallskip%
\begin{lemma}
  The statements {\it{}(i)}\/ and {\it (ii)}\/ above follow from the statement
  {\it{}(iii)}\/ below.
\end{lemma}

\smallskip%
\makebox[2em][l]{\em{}(iii)} If $F'$ is a nonempty subset of $F\setminus
F_{\infty}$, and $E'$ is the set of all edges of $\Sigma$ which are incident
with any face in $F'$, then
  \begin{equation}
    \label{eq:statement(iii)}
    2\pi|F'| \leq \sum_{e\in E'}2\theta^*(e),
  \end{equation}
  where equality holds if and only if $F'=F\setminus F_{\infty}$.

\smallskip%
\begin{proof}
  We deduce statement {\em{}(i)}. First let $F'=F\setminus F_{\infty}$ to
  obtain
  \begin{equation}
    \label{eq:sum_th_st_sigma_0}
    2\pi|F'| = \sum_{\genfrac{}{}{0pt}{}{\scriptstyle\text{edges}}{\scriptstyle\text{of } \Sigma_0}} 
    2\theta^*_e.
  \end{equation}
  Then let $F'=F\setminus (F_{\infty}\cup\{f\})$, where $f$ is a boundary
  face of $\Sigma_0$. Subtract the corresponding strict inequality from
  equation~\eqref{eq:sum_th_st_sigma_0} to obtain
  \begin{equation*}
    2\pi > \sum 2\theta^*_e,
  \end{equation*}
  where the sum is taken over all boundary edges of $\Sigma_0$ which are
  incident with $f$. Hence, assertion {\it (i)}\/ is true.
  
  We deduce statement {\em{}(ii)}. If $F_0'=F'$ then
  \begin{equation*}
    E_0'=E'\setminus\{\text{boundary edges of $\Sigma_0$}\}.
  \end{equation*}  
  By the definition of $\Phi_f$, inequality~\eqref{eq:statement(iii)} is
  equivalent to inequality~\eqref{eq:statement(ii)}.
\end{proof}

It is left to prove assertion {\em{}(iii)} under the assumption of
the condition of theorem~\ref{thm:rivin_circ}. We proceed in a
similar way as in section~\ref{sec:proof_thm_higher_genus}.

Suppose that $F'$ is a subset of $F\setminus F_{\infty}$, and $E'$ is the set
of all edges of $\Sigma$ which are incident with any face in $F'$.  Consider
$F''=F\setminus F'$, the complement of $F'$, and $E''=E\setminus E'$, the
complement of $E'$. Consider the Poincar\'e-dual cell decomposition
$\Sigma^*$ and, in it, the $1$-dimensional subcomplex (or graph)
$\Gamma=(F'', E'')$ with vertex set $F''$ and edge set $E''$. As for any
graph, we have
\begin{equation*}
|E''|-|F''|=c-n,
\end{equation*}
where $n$ is the number of connected components of $\Gamma$ and $c$ is the
dimension of the cycle space. Since the graph is embedded in $\Sigma^*$, a
cellular decomposition of the sphere, we have 
\begin{equation*}
c=r-1,
\end{equation*}
where $r$ is the number of regions into which $\Gamma$ separates $\Sigma^*$.
Since $E''$ contains the edges of $\Sigma$ incident with $v_{\infty}$, or,
dually, the edges of $\Sigma^*$ in the boundary of $v_{\infty}$, the number
of regions is at least two. Hence, the boundary of each region is nonempty.
By the condition of theorem~\ref{thm:rivin_circ}, the sum of $\theta$ over
the boundary of each region is at least $2\pi$. Sum over all regions to
obtain
\begin{equation*}
2\pi r\leq\sum_{e\in E''}2\theta_e.
\end{equation*}
Indeed, each edge in $E''$ appears in $0$ or $2$ boundaries. Equality holds
if and only if every edge of $E''$ is in the boundary of a region and each
boundary is the boundary of a single face of $\Sigma^*$. Thus, equality holds
if and only if $E''=E$ or $E''$ is the boundary of a single face of
$\Sigma^*$. This is the case, if and only if $F'=\emptyset$ (this is ruled
out by assumption) or $F'=F\setminus F_{\infty}$.

Thus, we have shown that
\begin{equation}\label{eq:R1}
2\pi(|E''|-|F''|)\leq \sum_{e\in E''}2\theta_e-2\pi(n+1),
\end{equation} 
with equality if and only if $F'=F\setminus F_{\infty}$. 

Equation~\eqref{eq:sum_theta_star} from
section~\ref{sec:proof_thm_higher_genus} also holds here. Thus,
inequality~\eqref{eq:R1} is equivalent to
\begin{equation*}
  -2\pi(|F|-|E|+|V|)+|F'|\leq\sum_{e\in E'}2\theta^*_e-2\pi(n+1),
\end{equation*}
and, using Euler's formula,
\begin{equation*}
|F|-|E|+|V|=2,
\end{equation*}
equivalent to
\begin{equation*}
2\pi|F'|\leq \sum_{e\in E'}2\theta^*_e-2\pi(n-1).
\end{equation*}
Since $n\geq 1$, and $n=1$ if $F'=F\setminus F_{\infty}$, we have deduced the
assertion {\em{}(iii)} above. This completes the proof of
theorem~\ref{thm:rivin_circ}.

\section{Proof of theorem~$\text{\ref{thm:rivin_poly}}$}
\label{sec:proof_thm_rivin_poly}

There is a one-to-one correspondence between Delaunay type circle patterns in
the sphere and polyhedra in hyperbolic 3-space with vertices in the infinite
boundary. In the Poincar{\'e} ball model, hyperbolic space corresponds to the
interior of the unit ball. The unit sphere corresponds to its infinite
boundary. Hyperbolic planes are represented by spheres that intersect the
unit sphere orthogonally. Hence, there is a correspondence between circles in
the unit sphere and hyperbolic planes. Furthermore, the intersection angle of
two circles equals the dihedral angle of the corresponding planes.  The
isometries of hyperbolic space correspond to the M{\"o}bius transformations
of the sphere at infinity.

\section{Proof of theorem~$\text{\ref{thm:ortho}}$}
\label{sec:proof_thm_ortho}

Theorem~\ref{thm:rivin_circ} has the following corollary.

\begin{corollary}
  Let $\Sigma$ be a strongly regular cell decomposition of the sphere, and
  suppose every vertex has $n$ edges. (Because $\Sigma$ is strongly regular,
  $n\geq 3$.) In other words, in the Poincar\'e dual decomposition
  $\Sigma^*$, every boundary of a face has $n$ edges. Suppose that every
  simple closed path in $\Sigma^*$ which is not the boundary of single face
  is more than $n$ edges long. Then there exists, uniquely up to M{\"o}bius
  transformations, a corresponding circle pattern in the sphere 
  with exterior intersection angles $2\pi/n$.
\end{corollary}

\begin{figure}
  \includegraphics[width=0.475\textwidth]{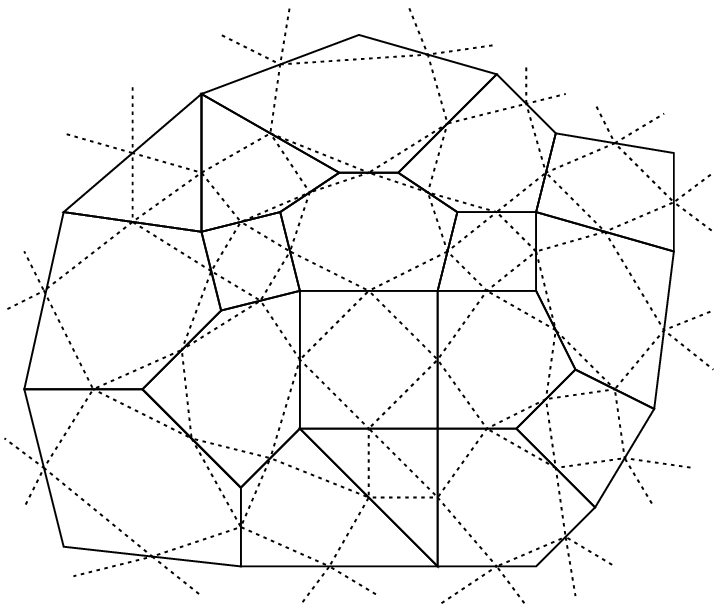} \hfill
  \includegraphics[width=0.475\textwidth]{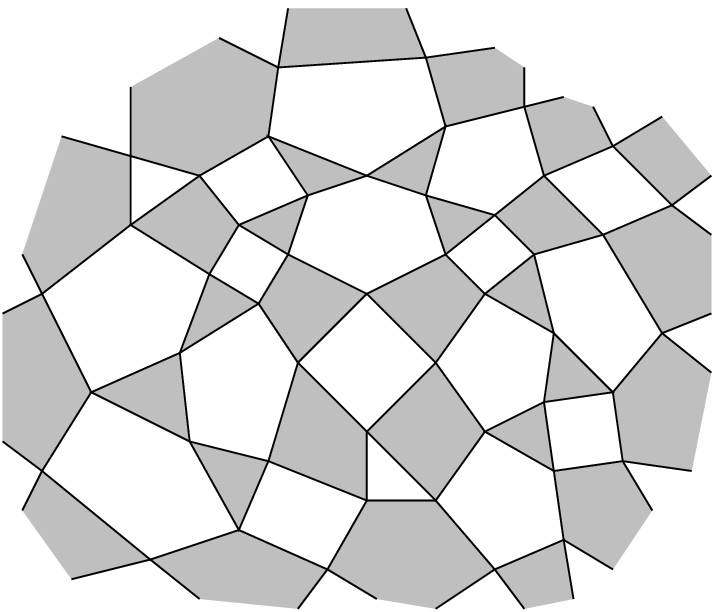}
\caption{A cellular decomposition {\em{}(left)}\/ and its medial
  decomposition {\em{}(right)}.}
\label{fig:medial}
\end{figure}

The case $n=4$ implies theorem \ref{thm:ortho}. Indeed, suppose $\Sigma$\/ is
a strongly regular cell decomposition of the sphere. In theorem
\ref{thm:ortho}, circles correspond to faces and vertices. To apply the
corollary, consider the {\em medial}\/ cell decomposition $\Sigma_m$ of
$\Sigma$. The faces of $\Sigma_m$ correspond to the faces and vertices of
$\Sigma$ and the vertices of $\Sigma_m$ correspond to edges of $\Sigma$.
Figure \ref{fig:medial} shows part of a cell decomposition of the sphere
{\em{}(left)}\/ and its medial decomposition {\em{}(right)}. The dotted lines
in the left figure represent the edges of the medial decomposition. On the
right, the faces of the medial decomposition which correspond to vertices in
the original decomposition are shaded.

The vertices of the medial decomposition $\Sigma_m$ are $4$-valent. The
assumption that $\Sigma$ is strongly regular implies that $\Sigma_m$ is
also strongly regular. It also implies that every simple closed path in the
Poincar\'e dual $\Sigma_m^*$ which is not the boundary of single face is more
than $4$ edges long. See the remark on page~\pageref{rem:regular}, point {\em
  (vi)}. Hence, the corollary implies theorem~\ref{thm:ortho}.

\section{Proof of theorem~$\text{\ref{thm:strong_steinitz}}$}
\label{sec:proof_thm:strong_steinitz}

The main part of theorem~\ref{thm:strong_steinitz} follows directly from
theorem~\ref{thm:ortho}. Given the cell decomposition $\Sigma$, by
theorem~\ref{thm:ortho}, there is a M\"obius-unique circle pattern with
orthogonally intersecting circles corresponding to the faces and vertices of
$\Sigma$. The circles which correspond to the faces lie in planes which form
the sides of the polyhedron in question. Indeed, the circles corresponding to
neighboring faces of $\Sigma$ touch, and hence the corresponding planes
intersect in a line which touches the sphere. Consider the faces of $\Sigma$
which are incident with one vertex $v$\/ of $\Sigma$. The corresponding
planes go through one point, namely, the apex of the cone touching the sphere
in the circle corresponding to the vertex $v$. Thus, one obtains a convex
polyhedron---provided, that is, that all vertex-circles are smaller than a
great circle. But one can always achieve this by a suitable M\"obius
transformation. We will show below that, by applying a suitable M\"obius
transformation, one can always get the center of gravity of the intersection
points of the orthogonal pattern into the center of the sphere. In that
position, all circles must be smaller than a great circle, because otherwise
all intersection points would lie in one hemisphere.

The construction above is reversible. Given a polyhedron with edges tangent
to the sphere, one obtains a circle pattern as in theorem~\ref{thm:ortho}.

The dual polyhedron is obtained by interchanging the role of face-circles and
vertex-circles.

Any M\"obius transformation of the sphere is the restriction of a projective
transformation of the ambient space which maps the sphere onto itself.
Conversely, a projective transformation of the ambient space which maps the
unit sphere onto itself induces a M\"obius transformation of the
sphere. Thus, the uniqueness claim of theorem~\ref{thm:strong_steinitz}
follows form the uniqueness claim of theorem~\ref{thm:ortho}.

It is left to show that, by applying a suitable M\"obius transformation to
the orthogonal circle pattern, one can get the center of gravity of the
intersection points into the center of the sphere. This follows from the
following lemma. The remainder of this section is devoted to proving it.

\begin{lemma}
  \label{lem:moebius_center}
  Let $v_1,\ldots,v_n$ be $n\geq 3$ distinct points in the $d$-dimensional
  unit sphere $S^d\subset\R^{d+1}$. There exits a M\"obius transformation $T$
  of $S^d$, such that
  \begin{equation*}
    \sum_{j=1}^n Tv_j=0.
  \end{equation*}
  If $\widetilde T$ is another such M\"obius transformation, then
  $\widetilde{T}=RT$, where $R$ is an isometry of $S^d$.
\end{lemma}

The proof relies on the close connection between the M\"obius geometry of
$S^d$ and the geometry of $(d+1)$-dimensional hyperbolic space $H^{d+1}$. In
the Poincar\'e ball model of hyperbolic space, $H^{d+1}$ is identified with
the unit ball in $\R^{d+1}$, its infinite boundary is $S^d$. The isometries
of $H^{d+1}$ extend to M\"obius transformations of $S^d$.  Conversely, every
M\"obius transformation of $S^d$ is the extension of a unique isometry of
$H^{d+1}$.

\begin{figure}
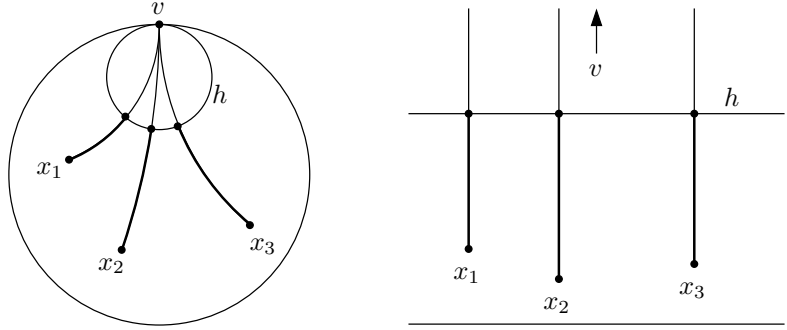

  \hfill
  \input{horosphere_ball.tex}
  \hfill
  \input{horosphere_halfspace.tex}
  \hspace*{\fill}
  \caption{The `distance' to an infinite point $v$ is measured by cutting 
    off at some horosphere through $v$. {\em Left:}\/ Poincar\'e ball model.
    {\em Right:}\/ half-space model.}
  \label{fig:horosphere}
\end{figure}

Given $n\geq 3$ points $v_1,\ldots,v_n\in S^d$, we are going show that there
is a unique point $x\in H^{d+1}$ such that the sum of the `distances' to
$v_1,\ldots,v_n$ is minimal. Of course, the distance to an infinite point is
infinite. The quantity to use is the distance to a horosphere through the
infinite point. See figure~\ref{fig:horosphere}

\begin{definition}
  For a horosphere $h$\/ in $H^{d+1}$, define
  \begin{gather*}
    \delta_h:H^{d+1}\rightarrow \R,\\
    \delta_h(x) = 
    \begin{cases}
      -\operatorname{dist}(x, h) & \text{if $x$ is inside $h$}, \\
      0 & \text{if $x\in h$}, \\
      \operatorname{dist}(x, h) & \text{if $x$ is outside $h$},
    \end{cases}
  \end{gather*}
  where $\operatorname{dist}(x, h)$ is the distance from the point $x$ to the
  horosphere $h$.
\end{definition}

Suppose $v$ is the infinite point of the horosphere $h$. Then the shortest
path from $x$ to $h$ lies on the geodesic connecting $x$ and $v$. If $h'$ is
another horosphere through $v$, then $\delta_h-\delta_{h'}$ is constant. If
$g:\R\rightarrow H^{d+1}$ is an arc-length parametrized geodesic, then
$\delta_h\circ g$ is a strictly convex function if $v$ is not an infinite
point of the geodesic $g$. Otherwise, $\delta_h\circ g(s)=\pm(s-s_0)$. These
claims are straightforward to prove using the Poincar\'e half-space model
where $v$ is the infinite point of the boundary plane. Also, one finds that 
\begin{equation*}
  \lim_{x\rightarrow\infty}\sum_{j=1}^n \delta_{h_j}(x)=\infty,
\end{equation*}
where $h_j$ are horospheres through different infinite points and
$n\geq3$. Thus, the following definition is proper.

\begin{definition}[\textsc{and Lemma}]
  Let $v_1,\ldots,v_n$ be $n$ points in the infinite boundary of $H^{d+1}$,
  where $n\geq 3$. Choose horospheres $h_1,\ldots,h_n$ through
  $v_1,\ldots,v_n$, respectively. There is a unique point $x\in H^{d+1}$
  for which $\sum_{j=1}^n \delta_{h_j}(x)$ is minimal. This point $x$ does not
  depend on the choice of horospheres. It is the {\em{}point of minimal
    distance sum}\/ from the infinite points $v_1,\ldots,v_n$.
\end{definition}

There seems to be no simpler characterization for the point of
minimal distance sum. However, it is easy to check whether it is the origin
in the Poincar\'e ball model.

\begin{lemma}
  Let $v_1,\ldots,v_n$ be $n\geq3$ different points in the infinite boundary
  of $H^{d+1}$. In the Poincar\'e ball model, $v_j\in S^d\subset\R^{d+1}$.
  The origin is the point of minimal distance sum, if and only if $\sum
  v_j=0$.
\end{lemma}

\begin{proof}
  If $h_j$ is a horosphere through $v_j$, then the gradient of $\delta_{h_j}$
  at the origin is the unit vector $-\frac{1}{2}v_j$.  (The metric is
  $ds^2=\big(\frac{2}{1-\sum x_j^2}\big)^2\sum dx_j^2$.)
\end{proof}

Lemma~\ref{lem:moebius_center} is now almost
immediate. Let $x$ be the point of minimal distance sum from the
$v_1,\ldots,v_n$ in the Poincar\'e ball model.  There is a hyperbolic
isometry $T$ which moves $x$ into the origin. If $\widetilde T$ is another
hyperbolic isometry which moves $x$ into the origin, then $\widetilde T=RT$,
with $R$ is an orthogonal transformation of $\R^{d+1}$.
Lemma~\ref{lem:moebius_center} follows.

This concludes the proof of theorem~\ref{thm:strong_steinitz}.

\chapter{Other variational principles}
\label{cha:other}

\section{Legendre transformations}
\label{sec:legendre}
 
In this section, we derive different variational principles for circle
patterns by Legendre transformations of the functionals $\Seuc$, $\Shyp$, and
$\Ssph$. We obtain the functional $\Shat$, defined below, which depends not
on the (transformed) radii $\rho\in\R^F$ but on the angles
$\varphi\in\R^{\vecEint}$.  According to whether the variation is constrained
to the space of euclidean, hyperbolic, or spherical coherent angle systems,
the critical points correspond to circle patterns of the respective geometry.
On the euclidean and hyperbolic coherent angle systems, the functional
$\Shat$ is strictly convex upwards, so that there can be only one critical
point, which is a maximum. (A function $f$ on a convex domain $D$ is called
{\em strictly convex upwards} {\em (downwards)}, if
\begin{equation*}
  f\big((1-t)x_1+tx_2\big)\genfrac{}{}{0pt}{1}{>}{(<)}(1-t)f(x_1)+tf(x_2)
\end{equation*}
for all $x_1,x_2\in D$ and $0<t<1$.) The function $\Shat(\varphi)$ is defined
in terms of Clausen's integral $\Cl(x)$, see appendix~\ref{app:dilog}.

\begin{theorem}
  \label{thm:Shat}
  Let $\vecEint$ be the set of oriented interior edges and define
  \begin{gather}
    \label{eq:Shat}
    \notag
    \Shat: \R^{\vecEint}\rightarrow \R\\
    \begin{split}
      \Shat(\varphi)= \sum_{e:\,\vece\,\uparrow\!\downarrow\,-\vece} 
      \Big(& \Cl(\theta^*_e+\varphi_{\vece}-\varphi_{-\vece})
      + \Cl(\theta^*_e-\varphi_{\vece}+\varphi_{-\vece})\\
      &+ \Cl(\theta^*_e+\varphi_{\vece}+\varphi_{-\vece}) +
      \Cl(\theta^*_e-\varphi_{\vece}-\varphi_{-\vece}) - 2\Cl(2\theta^*_e)
      \Big),
    \end{split}
  \end{gather}
  where the sum is taken over all non-oriented interior edges $e$, and $\vece$,
  $-\vece$ are the two oppositely oriented edges corresponding to it.
  
  (i) The function $\Shat$ is strictly convex upwards on the set of all
  euclidean coherent angle systems.
  
  If the function $\Seuc$ attains its minimum at $\rho\in\R^F$, then the
  restriction of $\Shat$ to the space of euclidean coherent angle systems
  attains its maximum at $\varphi\in\R^{\vecEint}$ defined by
  equations~\eqref{eq:phi_e_of_rho_f_euc}.
  
  Conversely, suppose the restriction of $\Shat$ to the space of euclidean
  coherent angle systems attains its maximum at $\varphi\in\R^{\vecE}$. Then
  the equations
  \begin{equation}
    \label{eq:rho_f_of_phi_e_euc}
    \rho_{f_k}-\rho_{f_j}=  
    \log\frac{\sin\varphi_{\vece}}{\sin(\varphi_{\vece}+\theta_e)}
    \quad\text{for}\quad f_j\;\underset{\vece}{\edgeup}\; f_k
  \end{equation}
  are compatible. (There is one equation for each oriented interior edge $\vece$ and
  $f_j$, $f_k$ are the faces to its left and right, respectively.) Therefore,
  they define $\rho\in\R^F$ uniquely up to an additive constant. The function
  $\Seuc$ attains its minimum at this~$\rho$.
  
  If they exist, the extremal values are equal:
  \begin{equation*}
    \min_{\rho\in\R^F}\Seuc(\rho)=
    \max_{\varphi \text{euc.~coherent}} \Shat(\varphi)
  \end{equation*}  
  
  (ii) The function $\Shat$ is strictly convex upwards on the set of
  hyperbolic coherent angle systems.
  
  If the function $\Shyp$ attains its minimum at $\rho\in\R^F$, then the
  restriction of $\Shat$ to the space of hyperbolic coherent angle systems
  attains its maximum at the $\varphi\in\R^{\vecEint}$ defined by
  equation~\eqref{eq:phi_e_of_rho_f_hyp}.
  
  Conversely, suppose the restriction of $\Shat$ to the space of hyperbolic
  coherent angle systems attains its maximum at $\varphi\in\R^{\vecEint}$.
  Then the equations
  \begin{equation}\label{eq:rho_f_of_phi_e_hyp}
    \rho_{f}=\frac{1}{2}\log
    \frac{\sin\Big(
      \frac{\displaystyle \theta^*-\varphi_{\vece}-\varphi_{-\vece}}
      {\displaystyle 2}\Big)
      \sin\Big(
      \frac{\displaystyle \theta^*-\varphi_{\vece}+\varphi_{-\vece}}
      {\displaystyle 2}\Big)}
    {\sin\Big(
      \frac{\displaystyle \theta^*+\varphi_{\vece}+\varphi_{-\vece}}
      {\displaystyle 2}\Big)
      \sin\Big(
      \frac{\displaystyle \theta^*+\varphi_{\vece}-\varphi_{-\vece}}
      {\displaystyle 2}\Big)}
  \end{equation}
  are compatible and define, therefore, a unique $\rho\in\R^F$.  (There is
  one equation for each oriented interior edge $\vece$, and $f$ is the face
  on its left side.) The function $\Shyp$ attains its minimum at this $\rho$.
  
  If they exist, the extremal values are equal:
  \begin{equation*}
    \min_{\rho\in\R^F}\Shyp(\rho)=
    \max_{\genfrac{}{}{0pt}{}{\scriptstyle\varphi\in\vecE,\text{ hyp.}}
      {\scriptstyle\text{coherent}}} \Shat(\varphi)
  \end{equation*}  
  
  (iii) If $\rho\in\R^F$ is a critical point of the function $\Ssph$, then
  $\varphi\in\R^{\vecEint}$ defined by equations~\eqref{eq:phi_e_of_rho_f_sph}
  is a critical point of the restriction of $\Shat$ to the space of spherical
  coherent angle systems.
  
  Conversely, suppose the restriction of $\Shat$ to the space of spherical
  coherent angle systems is critical at $\varphi\in\R^{\vecE}$. Then
  the equations
  \begin{equation}\label{eq:rho_f_of_phi_e_sph}
    \rho_{f}=\frac{1}{2}\log
    \frac{\sin\Big(
      \frac{\displaystyle -\theta^*+\varphi_{\vece}+\varphi_{-\vece}}
      {\displaystyle 2}\Big)
      \sin\Big(
      \frac{\displaystyle \theta^*-\varphi_{\vece}+\varphi_{-\vece}}
      {\displaystyle 2}\Big)}
    {\sin\Big(
      \frac{\displaystyle \theta^*+\varphi_{\vece}+\varphi_{-\vece}}
      {\displaystyle 2}\Big)
      \sin\Big(
      \frac{\displaystyle \theta^*+\varphi_{\vece}-\varphi_{-\vece}}
      {\displaystyle 2}\Big)}
  \end{equation}
  are compatible and define, therefore, a unique $\rho\in\R^F$. (There is one
  equation for each oriented interior edge $\vece$, and $f$ is the face on
  its left side.) The function $\Ssph$ is critical at this $\rho$.
  
  The values at corresponding critical points are equal:
  \begin{equation*}
    \Ssph(\rho_{\text{critical}})=\Shat(\varphi_{\text{critical}}).
  \end{equation*}
\end{theorem}

\begin{remark}
  1. If $\varphi$ is a euclidean coherent angle system, then
  equation~\eqref{eq:Shat} simplifies to 
  \begin{equation}
    \label{eq:ShatEuc}
    \Shat(\varphi)=\sum\big(\Cl(2\varphi_{\vece})+\Cl(2\varphi_{-\vece})
    -\Cl(2\theta^*_e)\big),
  \end{equation}
  where the sum is taken over all interior non-oriented edges $e$, and
  $\vece$, $-\vece$ are the corresponding oriented edges.
  
  2. Equations~\eqref{eq:rho_f_of_phi_e_hyp} and~\eqref{eq:rho_f_of_phi_e_sph}
  may be subsumed under the equation
  \begin{equation*}
    \rho_{f}=\frac{1}{2}\log\left|
      \frac{\sin\Big(
        \frac{\displaystyle \theta^*-\varphi_{\vece}-\varphi_{-\vece}}
        {\displaystyle 2}\Big)
        \sin\Big(
        \frac{\displaystyle \theta^*-\varphi_{\vece}+\varphi_{-\vece}}
        {\displaystyle 2}\Big)}
      {\sin\Big(
        \frac{\displaystyle \theta^*+\varphi_{\vece}+\varphi_{-\vece}}
        {\displaystyle 2}\Big)
        \sin\Big(
        \frac{\displaystyle \theta^*+\varphi_{\vece}-\varphi_{-\vece}}
        {\displaystyle 2}\Big)}\right|.
  \end{equation*}
\end{remark}

The rest of this section is devoted to the derivation of the functional
$\Shat(\varphi)$ and the proof of theorem~\ref{thm:Shat}.  

In classical mechanics, the Legendre transformation translates between the
Lagrangian and Hamiltonian descriptions of a mechanical system. (For a
thorough treatment, see, for example, Arnold's textbook \cite{arnold78}.)
The motion of such a system is typically described by functions $q(t)$ which
are critical for the functional $\int L(\dot q, q)\,dt$, where $L(v,q)$ is
the Lagrangian of the system. The Hamiltonian $H(p,q)$\/ is obtained by a
Legendre transformation with respect to $v$: $H(p,q)=pv-L(v,q)$, where
$p=\partial L(v,q)/\partial v$. It turns out that the motion of the system is
described by functions $p(t)$, $q(t)$ which are critical for the functional
$\int(pq-H(p,q))\,dt$. Here, $p$ and $q$ vary independently.

The construction of $\Shat$ proceeds in a similar way. The role of the
position variables $q$ is played by the variables $\rho$. Whereas $q$ depends
continuously on the time $t$, $\rho$ is a function on the finite set of
faces.  Instead of integrals, the functionals are finite sums. For each
functional, euclidean, hyperbolic, and spherical, we will define a Lagrangian
form first. From it, we will obtain a Hamiltonian form by a Legendre
transformation. A reduction leads to the functional $\Shat$ in every case.

First, we recapitulate the definition and basic properties of the Legendre
transformation of a smooth function of one variable. Then we deal with the
euclidean, hyperbolic, and spherical functionals separately.

\subsection{The Legendre transformation of a smooth function of one variable}
\label{sec:Legendre_one_variable}

The Legendre transformation is defined for convex functions. In general, they
need not be smooth and may depend on many variables. For the following
sections, it is sufficient to consider only the simplest case. The
generalization to many variables, which is needed in
section~\ref{sec:leibon_dual_func}, is straightforward.

\begin{definition}
  Suppose $F$ is a smooth real valued function on some interval in $\R$, and
  $F''>0$ (or $F''<0$).  Then
  \begin{equation}\label{eq:yOfx}
    y=F'(x)
  \end{equation}
  defines a smooth coordinate transformation. The {\em Legendre transform}\/
  of $F$ is
  \begin{equation}
    \label{eq:legendre_F}
    \widehat{F}(y) = xy - F(x),
  \end{equation}
  where $x$ is related to $y$ by (\ref{eq:yOfx}).
\end{definition}

The Legendre transformation is an involution: If $\widehat F$ is the Legendre
transform of $F$ and, say, $F''>0$, then $\widehat{F}''>0$ and $F$ is the
Legendre transform of $\widehat{F}$.  Indeed,
\begin{equation*}
  \widehat{F}'(y)\overset{\eqref{eq:legendre_F}}{=}
  \frac{dx}{dy}\,y + x - \frac{dx}{dy}F'(x)\overset{\eqref{eq:yOfx}}{=}x
\end{equation*}
and
\begin{equation*}
  \widehat{F}''(y) = \frac{dx}{dy} = \frac{1}{F''(x)}.
\end{equation*}

\subsection{The euclidean functional}
\label{sec:legendre_euc}

Equation~\eqref{eq:SEuc}, defining the functional $\Seuc$, may be rewritten
as
\begin{equation*}
  \Seuc(\rho) = \sum_{f_j\circ\edgeup\circ f_k} \big(
  F_{\theta_e}(\rho_{f_k}-\rho_{f_j})-(\pi-\theta_e)\rho_{f_k}\big)
  + \sum_{\circ f} \Phi_f\rho_{f}.
\end{equation*}
Here, the first sum is taken over all interior {\em{}oriented}\/ edges
$\vece$, and $f_j$ and $f_k$ are the faces on the left and right of $\vece$,
and $e$ is the corresponding non-oriented edge. The second sum is taken over
all faces $f$. The function $F_{\theta}(x)$ is defined in
equation~\eqref{eq:f_integral}. Define the Lagrangian form of the functional
by
\begin{gather*}
  \SeucL:\mathbb{R}^{\vecEint}\times\mathbb{R}^{F}\rightarrow\mathbb{R},\\
  \SeucL(v, \rho) = \sum_{f_j\circ\edgeup\circ f_k} \big(
  F_{\theta_e}(v_{\vece})-(\pi-\theta_e)\rho_{f_k}\big)
  + \sum_{\circ f} \Phi_f \rho_f.
\end{gather*}
If, for all oriented interior edges $\vece\in\vecE$,
\begin{equation*}
  v_{\vece}=\rho_{f_k}-\rho_{f_j},
\end{equation*}
where $f_j$ and $f_k$ to are the faces to the the left and to the right of
$\vece$, respectively, then
\begin{equation*}
  \SeucL(v,\rho)=\Seuc(\rho). 
\end{equation*}
In general, however, we do not assume that $v$ comes from a $\rho$ in this
way, or even that $v_{\vece}=-v_{-\vece}$.

By equations~\eqref{eq:f_prime} and~\eqref{eq:f_integral},
$F_{\theta}''(x)>0$. The Legendre transform of $F_{\theta}(x)$ is
\begin{equation*}
  \widehat{F}_{\theta}(y) = -\frac{1}{2}\Cl(2y) + \frac{1}{2}\Cl(2y+2\theta)
  - \frac{1}{2}\Cl(2\theta)\,,
\end{equation*}
where
\begin{equation*}
y=F_{\theta}'(x)=f_{\theta}(x).
\end{equation*}
This follows from equation~\eqref{eq:LiCl} of appendix~\ref{app:dilog}. By
equation~\eqref{eq:f_limits}, the domain of $\widehat{F}_{\theta}(y)$ is
\begin{equation*}
  0<y<\pi-\theta.
\end{equation*}
Define the `Hamiltonian' form of the functional by
\begin{gather*}
  \SeucH:D\times\R^F\longrightarrow\R,\\
  \SeucH(\varphi, \rho) = \sum_{f_j\circ\edgeup\circ f_k} 
  \big(\varphi_{\vece}(\rho_{f_k}-\rho_{f_j})-
  \widehat{F}_{\theta_e}(\varphi_{\vece})-(\pi-\theta_e)\rho_{f_k}\big)
  + \sum_{\circ j} \Phi_f\rho_{f},
\end{gather*}
where the $\varphi$-domain $D$\/ is
\begin{equation}
  \label{eq:domainD}
  D=\big\{\varphi\in\R^{\vecEint}\,\big|\;
  0<\varphi_{\vece}<\pi-\theta_{e}\,
  \text{\it{} for all }\,\vece\in\vecEint\big\}.
\end{equation}
(This is a superset of the space of euclidean coherent angle systems.) Writing
it all out, one has
\begin{multline*}
  \SeucH(\varphi, \rho)=
  \sum_{f_j\circ\edgeup\circ f_k}
  \Big(\varphi_{\vece}\,(\rho_{f_k}-\rho_{f_j})
  +\frac{1}{2}\big(\Cl(2\varphi_{\vece})-\Cl(2\varphi_{\vece}+2\theta_e)
  +\Cl(2\theta_e) \big)\\
  -(\pi-\theta_e)\rho_{f_k}\Big)+ \sum_{\circ f}
  \Phi_f  \rho_f.
\end{multline*}
If $\varphi$ and $\rho$ are
related by equations~\eqref{eq:phi_e_of_rho_f_euc}, then
\begin{equation*}
  \SeucH(\varphi, \rho)=\Seuc(\rho).
\end{equation*}
Now, $\widehat{F}_{\theta}''>0$. This follows from $F_{\theta}''>0$ and the
properties of the Legendre transformation. (It can also be verified by a direct
calculation.) Therefore, $\SeucH(\varphi,\rho)$ is strictly convex upwards
with respect to $\varphi$. 

The motivation for this construction of $\SeucH$ is, that the critical points
of $\Seuc$ correspond to critical points of $\SeucH$ where $\rho$ and
$\varphi$ may vary independently:
\begin{lemma}\label{lem:crit_points_euc}
  If $(\varphi, \rho)$ is a critical point of $\SeucH$, then $\rho$ is a
  critical point of $\Seuc$. Conversely, suppose $\rho$ is a critical point
  of $\Seuc$. Define $\varphi$ by equations~\eqref{eq:phi_e_of_rho_f_euc}. Then
  $(\varphi, \rho)$ is a critical point of $\SeucH$. At a critical point,
  $\Seuc(\rho)=\SeucH(\varphi,\rho)$.
\end{lemma}
\begin{proof}
  For any interior oriented edge $\vece\in\vecE$,
  \begin{equation*}
    \frac{\partial \SeucH}{\partial\varphi_{\vece}}(\varphi,\rho)=
    \rho_{f_k}-\rho_{f_j}-\widehat{F}_{\theta_e}'(\varphi_{\vece}),
  \end{equation*}
  Hence, all partial derivatives with respect to the variables $\varphi_{\vece}$ vanish, if and only if $\varphi$ and $\rho$ are related by
  $\rho_{f_k}-\rho_{f_j}=\widehat{F}_{\theta_e}'(\varphi_{\vece})$, or,
  equivalently, by $\varphi_{\vece}=F_{\theta_e}'(\rho_{f_k}-\rho_{f_k})
  =f_{\theta_e}(\rho_{f_k}-\rho_{f_j})$.  In this case, $\SeucH(\varphi,
  \rho)=\Seuc(\rho)$.
\end{proof}

Note that $\SeucH(\varphi,\rho)$ depends linearly on $\rho$. Collecting the
coefficients of each $\rho_f$, one finds that, if $\varphi$ is a euclidean
coherent angle system, then $\SeucH(\varphi, \rho)$ does not depend on $\rho$
at all.  In fact, in that case,
\begin{equation*}
  \SeucH(\varphi, \rho)=\Shat(\varphi).
\end{equation*}
To see this, note that if $\varphi$ is a coherent angle system, then
$\varphi_{\vece}+\theta_e=\pi-\varphi_{-\vece}$. Hence $\Cl(2\varphi_{\vece}+2\theta_e)=-\Cl(2\varphi_{-\vece})$. Also,
$\Cl(2\theta)=-\Cl(2\theta^*)$.

The convexity claim of theorem~\ref{thm:Shat}{\em{}(i)}\/ follows from the
convexity of $\SeucH(\varphi, \rho)$ with respect to $\varphi$.
Equations~\eqref{eq:rho_f_of_phi_e_euc} are the inverse relations to
equations~\eqref{eq:phi_e_of_rho_f_euc}. It is left to prove that they are
compatible if $\varphi$ is a critical point of $\Shat$ under variations in
the space of coherent angle systems. The rest follows from
lemma~\ref{lem:crit_points_euc}. Suppose $\varphi$ is a critical point of
$\Shat$ under variations in the space of coherent angle systems. Since we are
only interested in the restriction of $\Shat$ to the space of euclidean
coherent angle systems, we may use equation~\eqref{eq:ShatEuc} for $\Shat$.
The partial derivatives are then
\begin{equation*}
  \frac{\partial\Shat}{\partial\varphi_{\vece}}
  =-2\log(2\sin\varphi_{\vece}).
\end{equation*}
The tangent space to the space of euclidean coherent angle systems is spanned
by vectors of the form
\begin{equation*}
  \sum_{\vece\in\gamma}\Big(\frac{\partial}{\partial\varphi_{\vece}}
  -\frac{\partial}{\partial\varphi_{-\vece}}\Big),
\end{equation*}
where $\gamma$ is some oriented closed path in the Poincar\'e dual cell
decomposition.  (If we add to $\varphi_{\vece}$ we have to subtract the same
amount from $\varphi_{-\vece}$. But the sum of $\varphi$ around each face has
to remain constant, so we step from one face to the next, adding and
subtracting along some closed path of the Poincar\'e dual.) Since $\varphi$
is a critical point under variations in the space of coherent angle systems,
\begin{equation*}
  \sum_{\vece\in\gamma}
  \log\frac{\sin\varphi_{\vece}}{\sin(\varphi_{\vece}+\theta_e)}
  =0
\end{equation*}
for all closed paths $\gamma$ in the Poincar\'e dual. Hence,
equations~\eqref{eq:rho_f_of_phi_e_euc} are compatible. This completes the
proof of part {\em{}(i)}\/ of theorem~\ref{thm:Shat}.

\subsection{The hyperbolic functional}
\label{sec:legendre_hyp}

The construction in this case is similar to the euclidean case, but a bit more
involved.  This is due to the fact that the nonlinear terms in
equation~\eqref{eq:Shyp} (the definition of $\Shyp(\rho)$) depend not only on
the differences $\rho_k-\rho_j$, but also on the sums $\rho_k+\rho_j$. In
fact, equation~\eqref{eq:Shyp} may be written as
\begin{equation*}
  \Shyp(\rho) = \sum_{f_j\circ\edge\circ f_k} 
  \big(G_{\theta_e}(\rho_{f_k}-\rho_{f_j})
  + G_{\theta_e}(\rho_{f_j}+\rho_{f_k})\big)
  + \sum_{\circ f} \Phi_f\rho_f,
\end{equation*}
where
\begin{equation*}
  G_{\theta}(x) = F_{\theta}(x)+F_{\theta}(-x)=\im\Li\big(e^{x+i\theta}) + \im\Li\big(e^{-x+i\theta}).
\end{equation*}
The first sum is taken over all non-oriented edges $e$, and $f_j$ and $f_k$
are the faces on either side of $e$. The function $G_{\theta}(x)$ is even, so
that it does not matter which face is taken to be $f_j$ and which $f_k$. The
second sum is taken over all faces $f$.

We introduce some notation before we define the Lagrangian form of the
functional. The space $\mathbb{R}^{\vecEint}$ of functions on the oriented
interior edges splits into the direct sum of antisymmetric and symmetric
functions:
\begin{equation*}
  \mathbb{R}^{\vecEint}=\Alt(\vecEint)\oplus\Sym(\vecEint),
\end{equation*}
where 
\begin{equation}
  \label{eq:AltE}
  \Alt(\vecEint)=
  \big\{v\in\mathbb{R}^{\vecEint}\,\big|\,v_{-\vece}=-v_{\vece}\big\}
\end{equation}
and
\begin{equation*}
  \Sym(\vecEint)=
  \big\{w\in\mathbb{R}^{\vecEint}\,\big|\,w_{-\vece}=w_{\vece}\big\}.
\end{equation*}
Define the function 
\begin{gather*}
  \ShypL:\Alt(\vecEint)\times\Sym(\vecEint)\times\mathbb{R}^F
  \rightarrow\mathbb{R},\\
  \ShypL(v, w, \rho) = \sum_{f_j\circ\edge\circ f_k} 
  \big(G_{\theta_e}(v_{\vece}) + G_{\theta_e}(w_{\vece})\big) 
  + \sum_{\circ f} \Phi_f\rho_f.
\end{gather*}
The first sum is taken over all non-oriented interior edges $e$, and $\vece$
is one of the oriented representatives $e$. The second sum is taken over
all faces $f$. 

If, for all oriented edges $\vece\in\vecE$,
\begin{equation*}
  v_{\vece}=\rho_{f_k}-\rho_{f_j} \quad\text{and}\quad
  w_{\vece}=\rho_{f_k}+\rho_{f_j},    
\end{equation*}
where $f_j$ and $f_k$ are the faces to the the left and to the right of
$\vece$, respectively, then
\begin{equation*}
  \ShypL(v,w,\rho)=\Shyp(\rho). 
\end{equation*}
By a straightforward calculation,
\begin{equation*}
  G_{\theta}'(x)=f_{\theta}(x)-f_{\theta}(-x)=
  2\arctan\Big(\tan\Big(\frac{\theta^*_e}{2}\Big)
  \tanh\Big(\frac{y}{2}\Big)\Big),
\end{equation*}
and
\begin{equation*}
  G_{\theta}''(x)=f_{\theta}'(x)+f_{\theta}'(-x)=2f_{\theta}'(x)>0.
\end{equation*}
The Legendre transform of $G_{\theta}(x)$ is
\begin{equation*}
  \widehat{G}_{\theta}(y)=-\Cl(\theta^*+y)-\Cl(\theta^*-y)+\Cl(2\theta^*),
\end{equation*}
where $y=G_{\theta}'(x)$. This follows from equation~\eqref{eq:LiClp} of
appendix~\ref{app:dilog}. The domain of $\widehat{G}_{\theta}(y)$ is
$-\theta^*<y<\theta^*$.

The `Hamiltonian' form of the functional is 
\begin{gather*}
  \big\{p\in\Alt(\vecEint)\,\big|\;|p_{\vece}|<\theta^*_e\big\}
  \times
  \big\{s\in\Sym(\vecEint)\,\big|\;|s_{\vece}|<\theta^*_e\big\}
  \times\R^F\longrightarrow\R, \\
\begin{split}
  \ShypH(p, s, \rho) = 
  \sum_{f_j\circ\edge\circ f_k} &
  \big(p_{\vece} (\rho_{f_k}-\rho_{f_j})-\widehat{G}_{\theta_e}(p_{\vece}) 
     + s_{\vece} (\rho_{f_k}+\rho_{f_j})-\widehat{G}_{\theta_e}(s_{\vece})\\
  + \sum_{\circ f} & \Phi_f\rho_f,
\end{split}
\end{gather*}
where the first sum is taken over all non-oriented interior edges $e$, and
$\vece$ is one of the oriented representatives $e$. The second sum is taken
over all faces $f$. Writing it all out, we have
\begin{equation}\label{eq:ShypH}
  \begin{split}
    \ShypH(p,s,\rho)=
    \sum_{f_j\circ\edge\circ f_k} \Big(& 
    p_{\vece}\big(\rho_{f_k}-\rho_{f_j}\big)
    +\Cl\big(\theta^*_e+p_{\vece}\big)+\Cl\big(\theta^*_e-p_{\vece}\big)\\
    +&s_{\vece}\big(\rho_{f_j}+\rho_{f_k}\big)
    +\Cl\big(\theta^*_e+s_{\vece}\big)+\Cl\big(\theta^*_e-s_{\vece}\big)\\
    -&2\Cl\big(2\theta^*_e\big) \Big)
    +\sum_{\circ f} \Phi_f\rho_f.
  \end{split}
\end{equation}
The functional $\ShypH(p,s,\rho)$ is convex upwards in the variables $p$,
$s$.  If $p$ and $s$ are related to $\rho$ by
\begin{equation}
  \label{eq:p_of_rho_hyp}
  p_{\vece}=G_{\theta_e}'(\rho_k-\rho_j)=2\arctan\left(
    \tan\left(\frac{\theta^*_e}{2}\right)
    \tanh\left(\frac{\rho_k-\rho_j}{2}\right)
  \right)
\end{equation}
and
\begin{equation}
  \label{eq:s_of_rho_hyp}
  s_{\vece}=G_{\theta_e}'(\rho_k+\rho_j)=2\arctan\left(
    \tan\left(\frac{\theta^*_e}{2}\right)
    \tanh\left(\frac{\rho_k+\rho_j}{2}\right)
  \right)
\end{equation}
where $f_j$ and $f_k$ are the faces to the left and right of $\vece$,
respectively, then 
\begin{equation*}
  \ShypH(p,s,\rho)=S(\rho).
\end{equation*}
Again, the critical points of $\Shyp(\rho)$ correspond to the critical points
of $\ShypH(p,s,\rho)$ where $p$, $s$, and $\rho$ may vary independently.
(However, the domain of $\ShypH(p,s,\rho)$ has to be respected, of course. In
particular, the constraints $p\in\Alt(\vecEint)$ and $s\in\Sym(\vecEint)$
apply.)  The following lemma is proved in the same way as
lemma~\ref{lem:crit_points_euc}.

\begin{lemma}
  \label{lem:crit_points_hyp}
  If $(p, s, \rho)$ is a critical point of $\ShypH$, then $\rho$ is a
  critical point of $\Shyp$.  Conversely, suppose $\rho$ is a critical point
  of $\Shyp$. Define $p$ and $s$ by equations~\eqref{eq:p_of_rho_hyp}
  and~\eqref{eq:s_of_rho_hyp}. Then $(p,s,\rho)$ is a critical point of
  $\ShypH$.  At a critical point, $\Shyp(\rho)=\ShypH(p,s,\rho)$.
\end{lemma}

Now, we introduce the variables $\varphi\in\R^{\vecEint}$,
\begin{equation*}
\varphi_{\vece}={\textstyle\frac{1}{2}}(p_{\vece}-s_{\vece})
\end{equation*}
instead of $(p,s)$:
\begin{equation*}
\ShypH(\varphi,\rho)=\ShypH(p,s,\rho).
\end{equation*}
Since
\begin{equation*}
  \varphi_{-\vece}={\textstyle\frac{1}{2}}(-p_{\vece}-s_{\vece}),
\end{equation*}
we have 
\begin{align*}
  p_{\vece}&=\varphi_{\vece}-\varphi_{-\vece},\\
  s_{\vece}&=-\varphi_{\vece}-\varphi_{-\vece},
\end{align*}
and therefore
\begin{equation*}
  \begin{split}
\ShypH(\varphi,\rho)=&  \\ 
  \sum_{f_j\circ\edge\circ f_k} \Big(& 
  -2\big(\varphi_{\vece}\rho_{f_j}+\varphi_{-\vece}\rho_{f_k}\big)
  +\Cl\big(\theta^*_e+\varphi_{\vece}-\varphi_{-\vece}\big)
  +\Cl\big(\theta^*_e-\varphi_{\vece}+\varphi_{-\vece}\big)\\
  &+\Cl\big(\theta^*_e-\varphi_{\vece}-\varphi_{-\vece}\big)
  +\Cl\big(\theta^*_e+\varphi_{\vece}+\varphi_{-\vece}\big)
  -2\Cl\big(2\theta^*_e\big) \Big)\\
  +\sum_{\circ f} \Phi_f&\rho_f\,.
  \end{split}
\end{equation*}
The domain of $\ShypH(\varphi,\rho)$ is the set of all
$(\varphi,\rho)\in\R^{\vecEint}\times\R^F$ such that $\varphi$ satisfies the
inequalities
\begin{equation*}
  |\varphi_{\vece}+\varphi_{-\vece}|<\theta^*_e
  \quad\text{and}\quad
  |\varphi_{\vece}-\varphi_{-\vece}|<\theta^*_e.
\end{equation*}
This is a superset of the space of hyperbolic coherent angle systems.  The
relations~\eqref{eq:p_of_rho_hyp} and~\eqref{eq:s_of_rho_hyp} between $p$,
$s$ and $\rho$ are equivalent to the relation~\eqref{eq:phi_e_of_rho_f_hyp}
between $\varphi$ and $\rho$. Hence, if $\varphi$ and $\rho$ are related by
equation~\eqref{eq:phi_e_of_rho_f_hyp}, then
$\ShypH(\varphi,\rho)=\Shyp(\rho)$.  The functional $\ShypH(\varphi,\rho)$ is
linear in $\rho$. As in the euclidean case, collecting the coefficients of
each $\rho_f$, one observes that if $\varphi$ is a hyperbolic coherent angle
system, then $\ShypH$ does not depend on $\rho$.  In fact, in that case,
$\ShypH(\varphi,\rho)=\Shat(\varphi)$.

The convexity claim of theorem~\ref{thm:Shat}{\em{} (ii)}\/ follows from the
convexity of $\ShypH(\varphi, \rho)$ with respect to $\varphi$.
Equations~\eqref{eq:rho_f_of_phi_e_hyp} are the inverse relations to
equations~\eqref{eq:phi_e_of_rho_f_hyp}; see lemma~\ref{lem:phi_hyp}. It is
left to prove that they are compatible if $\varphi$ is a critical point of
$\Shat$ under variations in the space of hyperbolic coherent angle systems.
The rest follows from lemma~\ref{lem:crit_points_hyp}. Suppose the oriented
edges $\vece_1$ and $\vece_2$ are both in the oriented boundary of a face
$f$. We need to show that the value for $\rho_f$ obtained from the equation~\eqref{eq:rho_f_of_phi_e_hyp}
involving $\varphi_{\pm \vece_1}$ is the same as the value obtained from the
equation involving $\varphi_{\pm \vece_2}$. The tangent space to the space of
coherent angle systems is spanned by vectors of the form
\begin{equation*}
  \frac{\partial}{\partial\varphi_{\vece_2}}
  -\frac{\partial}{\partial\varphi_{\vece_1}}
\end{equation*}
where $\vece_1$ and $\vece_2$ are two oriented edges in the boundary of a
face. But $\partial\Shat/\partial\varphi_{\vece}$ is twice the right hand
side of equation (\ref{eq:rho_f_of_phi_e_hyp}). Thus, they are consistent if
$\varphi$ is a critical point of $\Shyp$ under variations in the space of
coherent angle systems. This completes the proof of part {\em{}(ii)}\/ of
theorem~\ref{thm:Shat}.

\subsection{The spherical functional}
\label{sec:legendre_sph}

This case is similar to the hyperbolic case. Define the `Lagrangian' form of
the spherical circle pattern functional $\Ssph$ (defined in
equation~\eqref{eq:Ssph}) to be
\begin{gather*}
  \SsphL:\Alt(\vecEint)\times\Sym(\vecEint)\times\mathbb{R}^F
  \rightarrow\mathbb{R},\\
  \SsphL(v, w, \rho) = \sum_{f_j\circ\edge\circ f_k} \big( 
  G_{\theta_e}(v_{\vece}) 
  - G_{\pi-\theta_e}(w_{\vece})
  -\pi(\rho_{f_j}+\rho_{f_k})\big) 
  + \sum_{\circ f} \Phi_f\rho_f.
\end{gather*}
If, for all oriented edges $\vece\in\vecE$,
\begin{equation*}
  v_{\vece}=\rho_{f_k}-\rho_{f_j} \quad\text{and}\quad
  w_{\vece}=\rho_{f_k}+\rho_{f_j},    
\end{equation*}
where $f_j$ and $f_k$ to are the faces to the the left and to the right of
$\vece$, respectively, then
\begin{equation*}
  \SsphL(v,w,\rho)=\Ssph(\rho). 
\end{equation*}
The `Hamiltonian' form is defined as
\begin{gather*}
  \big\{p\in\Alt(\vecEint)\,\big|\;|p_{\vece}|<\theta^*_e\big\}
  \times
  \big\{s\in\Sym(\vecEint)\,\big|\;|s_{\vece}|<\theta_e\big\}
  \times\R^F\longrightarrow\R, \\
\begin{split}
  \SsphH(p, s, \rho) = 
  \sum_{f_j\circ\edge\circ f_k} &
  \big(
  p_{\vece} (\rho_{f_k}-\rho_{f_j})
  - \widehat{G}_{\theta_e}(p_{\vece}) 
  - s_{\vece} (\rho_{f_k}+\rho_{f_j})
  + \widehat{G}_{\pi-\theta_e}(s_{\vece})\\
  &-\pi(\rho_{f_j}+\rho_{f_k})
  \big)
  + \sum_{\circ f}  \Phi_f\rho_f.
\end{split}
\end{gather*}
Writing it all out, we have
\begin{equation}\label{eq:SsphH}
  \begin{split}
    \SsphH(p,s,\rho)=
    \sum_{f_j\circ\edge\circ f_k} \Big(& 
    p_{\vece}\big(\rho_{f_k}-\rho_{f_j}\big)
    + \Cl\big(\theta^*_e+p_{\vece}\big)
    + \Cl\big(\theta^*_e-p_{\vece}\big)\\
    -&s_{\vece}\big(\rho_{f_j}+\rho_{f_k}\big)
    - \Cl\big(\theta_e+s_{\vece}\big)
    - \Cl\big(\theta_e-s_{\vece}\big)\\
    -&2\Cl\big(2\theta^*_e\big)-\pi\big(\rho_{f_j}+\rho_{f_k}\big)\Big)\\
    + \sum_{\circ f} \Phi_f&\rho_f.
  \end{split}
\end{equation}
(We have used that $\Cl(2\theta)=-\Cl(2\theta^*)$.) The functional
$\SsphH(p,s,\rho)$ is convex upwards in the variables $p$ and convex
downwards in the variables $s$. If $p$ and $s$ are related to $\rho$ by
\begin{equation}
  \label{eq:p_of_rho_sph}
  p_{\vece}=G_{\theta_e}'(\rho_k-\rho_j)=2\arctan\left(
    \tan\left(\frac{\theta^*_e}{2}\right)
    \tanh\left(\frac{\rho_k-\rho_j}{2}\right)
  \right)
\end{equation}
and
\begin{equation}
  \label{eq:s_of_rho_sph}
  s_{\vece}=G_{\theta^*_e}'(\rho_k+\rho_j)=2\arctan\left(
    \tan\left(\frac{\theta_e}{2}\right)
    \tanh\left(\frac{\rho_k+\rho_j}{2}\right)
  \right)
\end{equation}
where $f_j$ and $f_k$ are the faces to the left and right of $\vece$,
respectively, then 
\begin{equation*}
  \SsphH(p,s,\rho)=S(\rho).
\end{equation*}
Again, the critical points of $\Ssph(\rho)$ correspond to the critical points
of $\SsphH(p,s,\rho)$ where $p$, $s$, and $\rho$ may vary independently.
The following lemma is proved in the same way as
lemmas~\ref{lem:crit_points_euc} and~\ref{lem:crit_points_hyp}.

\begin{lemma}
  \label{lem:crit_points_sph}
  If $(p, s, \rho)$ is a critical point of $\SsphH$, then $\rho$ is a
  critical point of $\Ssph$.  Conversely, suppose $\rho$ is a critical point
  of $\Ssph$. Define $p$ and $s$ by equations~\eqref{eq:p_of_rho_sph}
  and~\eqref{eq:s_of_rho_sph}. Then $(p,s,\rho)$ is a critical point of
  $\SsphH$.  At a critical point, $\Ssph(\rho)=\SsphH(p,s,\rho)$.
\end{lemma}

Now we introduce the variables $\varphi\in\R^{\vecE}$,
\begin{equation*}
\varphi_{\vece}={\textstyle\frac{1}{2}}(p_{\vece}+s_{\vece}+\pi)
\end{equation*}
instead of $(p,s)$:
\begin{equation*}
\SsphH(\varphi,\rho)=\SsphH(p,s,\rho).
\end{equation*}
We have
\begin{align*}
  p_{\vece}&=\varphi_{\vece}-\varphi_{-\vece},\\
  s_{\vece}&=\varphi_{\vece}+\varphi_{-\vece}-\pi,
\end{align*}
and hence
\begin{equation*}
  \begin{split}
\SsphH(\varphi,\rho)=&  \\ 
  \sum_{f_j\circ\edge\circ f_k} \Big(& 
  -2\big(\varphi_{\vece}\rho_{f_j}+\varphi_{-\vece}\rho_{f_k}\big)
  +\Cl\big(\theta^*_e+\varphi_{\vece}-\varphi_{-\vece}\big)
  +\Cl\big(\theta^*_e-\varphi_{\vece}+\varphi_{-\vece}\big)\\
  &+\Cl\big(\theta^*_e-\varphi_{\vece}-\varphi_{-\vece}\big)
  +\Cl\big(\theta^*_e+\varphi_{\vece}+\varphi_{-\vece}\big)
  -2\Cl\big(2\theta^*_e\big)\Big)\\
  +\sum_{\circ f} \Phi_f&\rho_f\,.
  \end{split}
\end{equation*}
(We have used that $\Cl(x)$ is odd and $2\pi$-periodic.) The domain of
$\SsphH(\varphi,\rho)$ is the set of all
$(\varphi,\rho)\in\R^{\vecEint}\times\R^F$ such that $\varphi$ satisfies the
inequalities
\begin{equation*}
  |\varphi_{\vece}+\varphi_{-\vece}-\pi|<\theta_e
  \quad\text{and}\quad
  |\varphi_{\vece}-\varphi_{-\vece}|<\theta^*_e.
\end{equation*}
If $\varphi$ and $\rho$ are related by
equation~\eqref{eq:phi_e_of_rho_f_sph}, then
$\SsphH(\varphi,\rho)=\Ssph(\rho)$.  The functional $\SsphH(\varphi,\rho)$ is
linear in $\rho$. Collecting the coefficients of each $\rho_f$, one observes
that if $\varphi$ is a spherical coherent angle system, then $\SsphH$ does
not depend on $\rho$.  In fact, in that case,
$\SsphH(\varphi,\rho)=\Shat(\varphi)$.

Equations~\eqref{eq:rho_f_of_phi_e_sph} are the inverse relations to
equations~\eqref{eq:phi_e_of_rho_f_sph}; see lemma~\ref{lem:phi_sph}. It is
left to prove that they are compatible if $\varphi$ is a critical point of
$\Shat$ under variations in the space of spherical coherent angle systems.
This follows in exactly the same way as in the hyperbolic case. The rest
follows from lemma~\ref{lem:crit_points_sph}. This completes the proof of
part {\em{}(iii)}\/ of theorem~\ref{thm:Shat}.

\section{{Colin~de~Verdi\`ere}{}'s functionals}
\label{sec:colin_func}

\begin{figure}
\hfill
\input{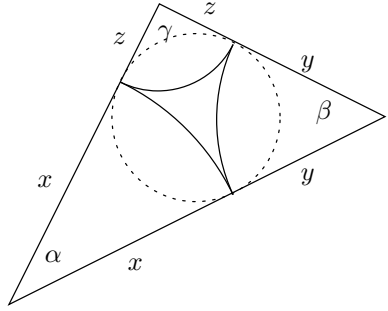}
\hspace*{\fill}
\caption{A euclidean triangle}
\label{fig:eucTri}
\end{figure}

Colin~de~Verdi{\`e}re \cite{colin91} considers circle packings in which the
circles correspond to the vertices of a triangulation. He considers the
$1$-form
\begin{equation*}
\omega = \alpha\,du+\beta\,dv+\gamma\,dw
\end{equation*}
on the space of euclidean triangles, where 
\begin{equation*}
  u=\log x, \quad v=\log y, \quad\text{and}\quad w=\log z
\end{equation*}
and $x$, $y$, $z$ and $\alpha$, $\beta$, $\gamma$ are as shown in
figure~\ref{fig:eucTri}.

It turns out that $d\omega=0$, hence one may integrate. Define the
function $f_{\alpha_0, \beta_0, \gamma_0}$ on $\mathbb{R}^3$ by
\begin{equation*}
  f_{\alpha_0, \beta_0, \gamma_0}(u,v,w) = 
  \int^{(u,v,w)} (\alpha_0 - \alpha)\,du + (\beta_0 - \beta)\,dv
  + (\gamma_0 - \gamma)\,dw.
\end{equation*}
The initial point of the integration does not matter. 

Suppose we are given a triangulation and a coherent angle system for it.
Here, a coherent angle system is a positive function on the set of angles of
the triangles, such that the sum in each triangle is $\pi$, and the sum
around each vertex is $2\pi$. For a function $\rho$ on the vertices of the
triangulation, Colin~de~Verdi{\`e}re's functional for euclidean circle
packings is
\begin{equation*}
  S_{\text{\it CdV}}(\rho)=\sum f_{\alpha_0, \beta_0, \gamma_0}(u, v, w),
\end{equation*}
where the sum is taken over all triangles, $u$, $v$, and $w$ are the values
of $\rho$ on the vertices of each triangle, and $\alpha_0$, $\beta_0$ and
$\gamma_0$ are the corresponding angles of the coherent angle system.

The critical points of this functional correspond to the logarithmic radii of
a circle packing.

The hyperbolic case is treated in the same way, except that now 
\begin{equation*}
    u=\log \tanh(x/2), \quad v=\log \tanh(y/2), 
    \quad\text{and}\quad w=\log \tanh(z/2),
\end{equation*}
and in the definition of a coherent angle system, it is required that the
sum of the angles in a triangle is less than $\pi$.

To treat circle packings with our functionals, we consider circle patterns
with orthogonally intersecting circles on the medial decomposition of the
triangulation (see figure~\ref{fig:medial}). 

Suppose $\Sigma_0$ is a triangulation. Let $\Sigma$ be the medial
decomposition. Set $\theta_e=\pi/2$ for all edges $e$. For simplicity, we
will consider only the case of closed surfaces, so let $\Phi_f=2\pi$ for all
faces $f$. The faces of $\Sigma$ are of two types: those that correspond to
faces of $\Sigma_0$ and those that correspond to vertices of $\Sigma_0$.  Let
$F=F_1\cup F_2$, where $F_1$ contains the faces of the first type and $F_2$
contains the faces of the second type.

Consider $\Seuc$ as function on $\mathbb{R}^{F_1}\times\mathbb{R}^{F_2}$ and
define $S_2:\mathbb{R}^{F_2}\rightarrow \mathbb{R}$,
\begin{equation*}
S_2(\rho_2)=\min_{\rho_1\in\mathbb{R}^{F_1}} 
\Seuc(\rho_1,\rho_2).
\end{equation*}

Hence, $S_2(\rho_2)=\Seuc(\rho_1,\rho_2)$, where $\rho_1\in\mathbb{R}^{F_1}$
is determined as follows.  Suppose $f\in F_1$ and $f_a, f_b, f_c\in F_2$ are
the neighboring faces of $f$. Construct the euclidean triangle whose sides
are
\begin{equation*}
e^{{\rho_2}_{f_a}}+e^{{\rho_2}_{f_b}},\quad
e^{{\rho_2}_{f_b}}+e^{{\rho_2}_{f_c}},\quad\text{and}\quad
e^{{\rho_2}_{f_c}}+e^{{\rho_2}_{f_a}}.
\end{equation*}
Let $\rho_1(f)$ be the logarithmic radius of the inscribed circle.

For $\rho_2\in\mathbb{R}^{F_2}$, and $\rho_1$ the corresponding
point in $\mathbb{R}^{F_1}$,
\begin{equation*}
  \frac{\partial \Seuc}{\partial \rho_1}(\rho_1, \rho_2)  = 0,\quad
  \frac{\partial \Seuc}{\partial \rho_2}(\rho_1, \rho_2) =
  \frac{\partial S_2}{\partial \rho_2}(\rho_2).
\end{equation*}
Consider $S_{\text{\it CdV}}$ as function on $\mathbb{R}^{F_2}$. It is not
hard to see that $dS_2=\,dS_{\text{\it CdV}}$.  This implies that
Colin~de~Verdi{\`e}re's euclidean functional is, up to an additive constant,
equal to $S_2$.

Colin~de~Verdi{\`e}re's functional for hyperbolic circle packings
can be derived from $\Shyp$ in the same way.

\section[Thurston type circle patterns with ``holes'']{Digression: Thurston type circle patterns with ``holes''}
\label{sec:Thurston_with_holes}

\begin{figure}%
\hfill%
\includegraphics[height=0.24\textwidth]{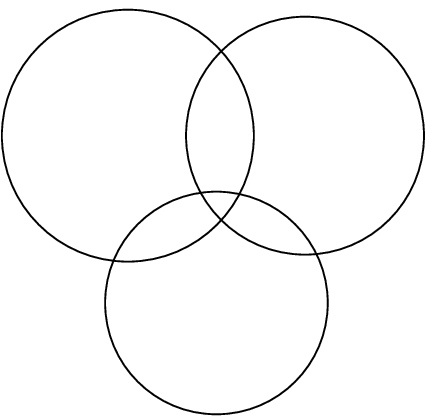}%
\hfill\hfill%
\includegraphics[height=0.24\textwidth]{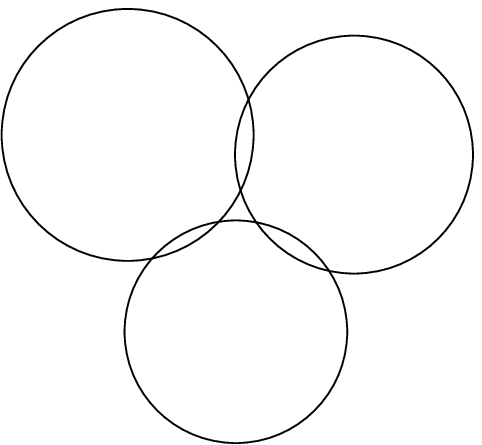}%
\raisebox{0.1\textwidth}{ $\longrightarrow$ }%
\includegraphics[height=0.24\textwidth]{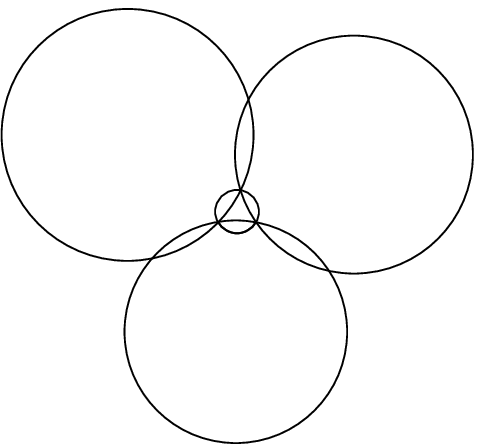}%
\hspace*{\fill}
\caption{{\it Left:} Thurston type circle pattern with overlap. {\it
    Middle:} Thurston type circle pattern with hole. {\it Right:} Add
  a circle to fill the hole.}%
\label{fig:addcirc}%
\end{figure}

\begin{figure}%
\centering%
\input{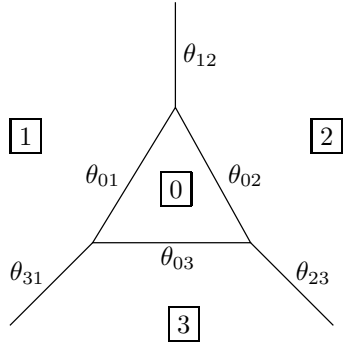}%
\caption{A ``truncated vertex". The old faces are labeled $1$, $2$, $3$; the
  new face is labeled $0$.}%
\label{fig:hole_fill_angles}%
\end{figure}

A similar reduction leads to functionals for Thurston type circle patterns
with ``holes"~\cite{springborn03}. Thurston type circle patterns were
mentioned in the introduction, section~\ref{sec:ex_and_uni}. Here, the
vertices do not correspond to intersection points. All vertices have three
edges, but the intersection angles $\theta$\/ do not have to sum to $2\pi$.
If the sum is less than $2\pi$, the circles will overlap as in
figure~\ref{fig:addcirc} {\it(left)}.
If the sum is greater than $2\pi$, there will be a hole as in
figure~\ref{fig:addcirc} {\it(middle)}. If the sum is greater than $2\pi$
around {\em{}all}\/ vertices, one obtains a Delaunay type circle pattern by
adding circles as shown in figure~\ref{fig:addcirc} {\it(right)}.
Combinatorially, this corresponds to ``truncating'' all vertices. The
intersection angles for the new edges (see figure~\ref{fig:hole_fill_angles})
are determined by the equations
\begin{equation*}
    \begin{split}
    \theta_{12} + \theta_{20} + \theta_{01} &= 2\pi \\
    \theta_{23} + \theta_{30} + \theta_{02} &= 2\pi \\
    \theta_{31} + \theta_{10} + \theta_{03} &= 2\pi,
  \end{split}
\end{equation*}
so that 
\begin{equation*}
  \theta_{01} = \pi - \frac{1}{2}(\theta_{12} - \theta_{23} + \theta_{31}),
  \quad\text{etc.}
\end{equation*}

\section{Br\"agger{}'s functional}
\label{sec:bragger_func}

Like Colin~de~Verdi{\`e}re, Br{\"a}gger \cite{bragger92} considers circle
packings in which the circles correspond to the vertices of a triangulation.
A coherent angle system in the sense of section~\ref{sec:colin_func} is
clearly equivalent to a coherent angle system in the sense of this article.
Thus, Br{\"a}gger's functional is seen to be equal to $\Shat/2$ up to an
additive constant.

\section{Rivin's functional}
\label{sec:rivin_func}

Rivin \cite{rivin94} considers euclidean circle patterns with arbitrary
prescribed intersection angles. The pattern of intersection is determined by
a triangulation. Circles correspond to faces of the triangulation. However,
since Rivin allows intersection angles $\theta=0$, adjacent triangles may in
effect belong to the same circle. Thus, cellular decompositions with
non-triangular cases may be treated by first dissecting all faces into
triangles. Rivin treats cone-like singularities in the vertices but not in
the centers of the circles. His functional is up to an additive constant
equal to $\Shat/2$ with $\Shat$ as in equation~\eqref{eq:ShatEuc}.

\section{Leibon's functional}
\label{sec:leibon_func}

Like Rivin, Leibon \cite{leibon99}, \cite{leibon02} considers circle patterns
with arbitrary prescribed intersection angles. The circles correspond to the
faces of a triangulation and the variables are the angles of the inscribed
triangles.  Whereas Rivin treats the euclidean case, Leibon treats the
hyperbolic case. Leibon's functional cannot be derived directly from our
functionals. In section~\ref{sec:common_ancestor}, we will construct yet
another functional, from which both Leibon's functional and $\Shat(\varphi)$
can be derived. This section is an exposition of Leibon's variational
principle.

Let $\Sigma$ be a triangulation of a compact surface without boundary and let
intersection angles be prescribed by a function $\theta:E\rightarrow(0,\pi)$
on the non-oriented edges. Assume that $\theta$ sums to $2\pi$ around each
vertex. (This assumption is made only to simplify the exposition. It means
that we consider only the case without cone-like singularities. If the
assumption is violated, one obtains prescribed cone-like singularities in the
vertices. Cone-like singularities in the centers of the circles cannot be
treated by this method.) The variables of Leibon's functional are the
angles of the triangles of the triangulation.  We associate each angle of a
triangle with the opposite directed edge in the oriented boundary of that
triangle. Thus, for a geodesic triangulation of a hyperbolic surface, the
angles of the triangles define a function on the oriented edges.  (In the
figures, however, we label the angles in the usual way.)

For a function $\alpha\in\R^{\vecE}$ on the oriented edges, Leibon's
functional is
\begin{equation*}
  H(\alpha)=\sum_{t\in F} V(\alpha_1^t,\alpha_2^t,\alpha_3^t),
\end{equation*}
where the sum is taken over all triangles $t$;
$\alpha_1^t,\alpha_2^t,\alpha_3^t$ are the three angles of a triangle $t$,
and
\begin{equation}
  \label{eq:V_alpha_123}
  \begin{split}
    V(\alpha_1,\alpha_2,\alpha_3)=\frac{1}{2}\Big(&
    \Cl(2\alpha_1)+\Cl(2\alpha_2)+\Cl(2\alpha_3)\\
    &+\Cl(\pi+\alpha_1-\alpha_2-\alpha_3)
    +\Cl(\pi-\alpha_1+\alpha_2-\alpha_3)\\
    &+\Cl(\pi-\alpha_1-\alpha_2+\alpha_3)
    +\Cl(\pi-\alpha_1-\alpha_2-\alpha_3)\Big).
  \end{split}
\end{equation}

A {\em coherent angle system}\/ in this setting is a function
$\alpha\in\R^{\vecE}$ on the oriented edges, which satisfies the following
conditions:

\makebox[2em][l]{\it (i)} For all $\vece\in\vecE$, $\alpha_{\vece}>0$.

\nopagebreak
\makebox[2em][l]{\it (ii)} For all triangles $t\in F$, 
$\alpha_1^t+\alpha_2^t+\alpha_3^t<\pi$.

\begin{figure}
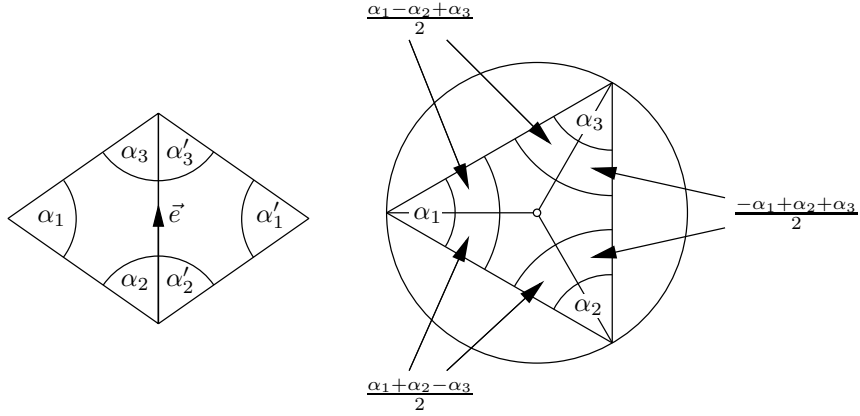
%
\hfill%
\raisebox{1.2cm}{\input{adjacent_triangles.tex}}%
\hfill%
\input{leibon_angles.tex}%
\hspace*{\fill}%
\caption{{\it Left:}\/ Adjacent triangles. {\it Right:}\/ Triangle with
  circumscribed circle. (For simplicity, euclidean triangles are drawn in the
  figure.)}%
\label{fig:leibon_angles}%
\end{figure}

\nopagebreak
\makebox[2em][l]{\it (iii)} The angles of two adjacent triangles with common
edge $e\in E$, as shown in figure~\ref{fig:leibon_angles} {\it (left)}, satisfy
\begin{equation*}
  \frac{-\alpha_1+\alpha_2+\alpha_3}{2}
  +\frac{-\alpha'_1+\alpha'_2+\alpha'_3}{2}
  =\theta_e.
\end{equation*}

Note that condition {\it (iii)}\/ implies that the angles $\alpha$ sum to
$2\pi$ around each vertex (since the angles $\theta$ do). For a Delaunay
triangulation with circles intersecting at angles given by $\theta$, the
angles of the triangles form a coherent angle system.
Figure~\ref{fig:leibon_angles} {\it (right)}\/ shows that condition {\it
  (iii)} is satisfied.

\begin{lemma}[Leibon~\cite{leibon02}]
  \label{lem:leibon_principle}
  The hyperbolic triangles with angles prescribed by $\alpha\in\R^{\vecE}$
  fit together to form a Delaunay triangulation of a hyperbolic surface with
  circles intersecting at angles given by $\theta$ if and only if
  $\alpha\in\R^{\vecE}$ is a critical point of the functional $H(\alpha)$
  under variations in the space of coherent angle systems.
\end{lemma}

Leibon also shows that the functional is strictly convex upwards and uses his
variational principle to prove an existence and uniqueness theorem for
hyperbolic circle patterns, which is a special case of
theorem~\ref{thm:fundamental}~\cite{leibon02}.

\begin{proof}[Proof of lemma~\ref{lem:leibon_principle}]
  It is clear that, if the triangles may be glued together, one obtains a
  Delaunay triangulation of a hyperbolic surface with circles intersecting at
  angles given by $\theta$. The triangles fit together if sides along which
  two triangles are supposed to be glued have the same length. Hence, the
  lemma follows from the following facts:

  {\it (a)}\/ For each oriented edge $\vece\in\vecE$, let
  \begin{equation*}
    w_{\vece}
    =\frac{\partial}{\partial\alpha_2}+\frac{\partial}{\partial\alpha_3}
    -\frac{\partial}{\partial\alpha'_2}-\frac{\partial}{\partial\alpha'_3},
  \end{equation*}
  where $\alpha_2,\alpha_3,\alpha'_2,\alpha'_3$ are as in
  figure~\ref{fig:leibon_angles} {\it (left)}.
  The tangent space to the space of coherent angle systems is spanned by
  the vectors $w_{\vece}$, $\vece\in\vecE$.

  {\it (b)}\/ For positive $\alpha_1$, $\alpha_2$, $\alpha_3$ with
  $\alpha_1+\alpha_2+\alpha_3<\pi$,
  \begin{equation}
    \label{eq:leibon_V_deriv}
    \Big(\frac{\partial}{\partial\alpha_2}+\frac{\partial}{\partial\alpha_3}
    \Big)\,V(\alpha_1,\alpha_2,\alpha_3)=2\log\sinh\frac{a_1}{2},
  \end{equation}
  where $a_1$\/ is the length of the side opposite $\alpha_1$ in a hyperbolic
  triangle with angles $\alpha_1$, $\alpha_2$, $\alpha_3$.
  
  \paragraph{Proof of {\it (a)}}
  Consider the linear functions 
  \begin{equation*}
    \begin{aligned}
      &f:   & \Alt(\vecE)&\longrightarrow\R^{\vecE}, \\
      &g:   & \R^{\vecE}&\longrightarrow\R^E, \\
      &g^*: & \R^E&\longrightarrow\R^{\vecE}
    \end{aligned}
  \end{equation*}
  defined by figure~\ref{fig:functions_fg}. 
  \begin{figure}
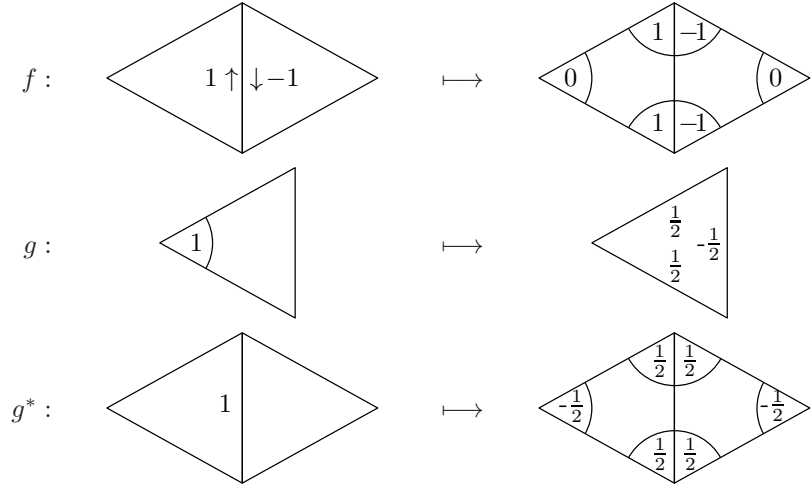
%
    \begin{equation*}%
      \begin{aligned}
        f&: & 
        &\quad\raisebox{-0.9cm}{\input{function_f_pre.tex}} &
        &\quad{\longmapsto} &
        &\quad\raisebox{-0.9cm}{\input{function_f_im.tex}} \\
        g&: &
        &\quad\quad\quad\raisebox{-0.9cm}{\input{function_g_pre.tex}} &
        &\quad{\longmapsto} &
        &
        \quad\quad\quad\raisebox{-0.9cm}{\input{function_g_im.tex}} \\
        g^*&: &
        &\quad\raisebox{-0.9cm}{\input{function_gst_pre.tex}} &
        &\quad{\longmapsto} &
        & \quad\raisebox{-0.9cm}{\input{function_gst_im.tex}}
      \end{aligned}
    \end{equation*}%
    \caption{The linear functions $f$, $g$, and $g^*$.}%
    \label{fig:functions_fg}%
  \end{figure}
  ($\Alt(\vecE)$ denotes the space
  of antisymmetric functions on the oriented edges; see
  equation~\eqref{eq:AltE}.) Suppose the vector spaces $\R^{\vecE}$, and
  $\R^E$ are equipped with the canonical scalar product. Then $g$ and $g^*$
  are adjoint linear operators. The tangent space to the space of coherent
  angle systems is
  \begin{equation*}
    (\image g^*)^{\perp} = \kernel g.
  \end{equation*}
  Hence we need to show that 
  \begin{equation}
    \label{eq:im_f_eq_ker_g}
    \image f= \kernel g.
  \end{equation}
  In fact, we will show that 
  \begin{equation*}
    0\longrightarrow\Alt(\vecE)
    \overset{f}{\longrightarrow}\R^{\vecE}
    \overset{g}{\longrightarrow}\R^E
    \longrightarrow 0
  \end{equation*}
  is an exact sequence. It is easy to see that 
  \begin{equation*}
    \kernel f=0
  \end{equation*}
  and 
  \begin{equation*}
    \image f\subseteq\kernel g
  \end{equation*}
  It is also easy to see that
  \begin{equation*}
    \kernel g^*=0
  \end{equation*}
  and hence 
  \begin{equation*}
    \image g=(\kernel g^*)^{\perp}=\R^E.
  \end{equation*}
  Since
  \begin{equation*}
    \dim\kernel g=|\vecE|-|E|=|E|=\dim\image f,
  \end{equation*}
  equation~\eqref{eq:im_f_eq_ker_g} follows. This proves {\it (a)}.

  \paragraph{Proof of {\it (b)}}
  By equation~\eqref{eq:clausen_diff},
  \begin{equation*}
    \Big(\frac{\partial}{\partial\alpha_2}+\frac{\partial}{\partial\alpha_3}
    \Big)\,V(\alpha_1,\alpha_2,\alpha_3)=
    \log\frac{
      \cos\frac{\alpha_1+\alpha_2+\alpha_3}{2}
      \cos\frac{-\alpha_1+\alpha_2+\alpha_3}{2}}{\sin\alpha_2\sin\alpha_3}.
  \end{equation*}
  From the hyperbolic angle cosine theorem (see, for example,
  Beardon~\cite{beardon83})
  \begin{equation*}
    \cosh a_1=
    \frac{\cos\alpha_2\cos\alpha_3+\cos\alpha_1}{\sin\alpha_2\sin\alpha_3}
  \end{equation*}
  one derives the half-angle formula
  \begin{equation*}
    \sinh^2\frac{a_1}{2}=\frac{\cos\frac{\alpha_1+\alpha_2+\alpha_3}{2}
    \cos\frac{-\alpha_1+\alpha_2+\alpha_3}{2}}{\sin\alpha_2\sin\alpha_3}\,.
  \end{equation*}
  This proves {\it (b)}.
\end{proof}

\section{The Legendre dual of Leibon's functional}
\label{sec:leibon_dual_func}

In this section, we construct the circle pattern functional $\widehat H$,
which is Legendre dual to Leibon's functional. The variables are
(essentially) the edge lengths of the triangles. Unfortunately, we cannot
give an explicit formula for $\widehat H$ in terms of the dilogarithm or
related functions. Whether a manageable formula exists, seems to be an
interesting question.

Introduce new variables $\beta\in\R^{\vecE}$ instead of
$\alpha\in\R^{\vecE}$: For each triangle, the angles $\beta_1$, $\beta_2$,
$\beta_3$ are related to $\alpha_1$, $\alpha_2$, $\alpha_3$ by
\begin{equation}
  \label{eq:beta_123}
  \begin{split}
  \beta_1 &= \frac{-\alpha_1+\alpha_2+\alpha_3}{2} \\
  \beta_2 &= \frac{\alpha_1-\alpha_2+\alpha_3}{2} \\
  \beta_3 &= \frac{\alpha_1+\alpha_2-\alpha_3}{2}\,.    
  \end{split}
\end{equation}
Note that
\begin{equation*}
  \frac{\partial}{\partial\beta_1} =
  \frac{\partial}{\partial\alpha_2}+\frac{\partial}{\partial\alpha_3},
  \quad \text{etc.},
\end{equation*}
such that equation~\eqref{eq:leibon_V_deriv} is equivalent to
\begin{equation*}
  \frac{\partial V}{\partial\beta_j}=2\log\sinh a_j.
\end{equation*}
Let $\widehat V$ be the Legendre transform of $V$, considered as a function
of $\beta_1$, $\beta_2$, $\beta_3$. (See
section~\ref{sec:Legendre_one_variable}.) That is, let
\begin{multline*}
  \widehat V(2\log\sinh a_1,2\log\sinh a_2,2\log\sinh a_3)= \\
  2\beta_1\log\sinh a_1 + 2\beta_2\log\sinh a_2 + 2\beta_3\log\sinh a_3
  - V(\alpha_1, \alpha_2, \alpha_3),
\end{multline*}
where $\alpha_1$, $\alpha_2$, $\alpha_3$ and $a_1$, $a_2$, $a_3$ are the
angles and side lengths of a hyperbolic triangle. Then
\begin{equation*}
  \frac{\partial\widehat V}{\partial\,(2\log\sin a_j)}=\beta_j.
\end{equation*}
For a function $a\in\R^E$ on the non-oriented edges which satisfies the
triangle inequalities, define the functional
\begin{equation*}
  \widehat H(2\log\sinh a)=\sum_{t\in F}
  \widehat V(2\log\sinh a^t_1,2\log\sinh a^t_2,2\log\sinh a^t_3) 
  - \sum_{e\in E} 2\theta_e\log\sinh a_e.
\end{equation*}
The first sum is taken over all triangles $t$, and $a^t_1$, $a^t_2$, $a^t_3$
are the sides of triangle $t$. The second sum is taken over all non-oriented
edges $e$. The functional $\widehat H$ is the Legendre dual of Leibon's
functional $H$:
\begin{proposition}
  Suppose $a\in\R^E$ satisfies the triangle inequalities. Then the hyperbolic
  triangles with edge lengths given by $a$ form a Delaunay triangulation of a
  hyperbolic surface with circles intersecting at angles given by $\theta$ if
  and only if $2\log\sinh a$ is a critical point of the functional $\widehat
  H$.
\end{proposition}

\begin{proof}
  Suppose $\vece$ and $-\vece$ are the two orientations of the unoriented
  edge $e\in E$. Then the partial derivative of $\widehat H$ with respect to
  $2\log\sinh a_e$ is
  \begin{equation*}
    \frac{\partial \widehat H(2\log\sinh a)}{\partial\,(2\log\sinh a_e)}=
    \beta_{\vece}+\beta_{-\vece}-\theta_e.
  \end{equation*}
  This partial derivative vanishes if and only if the circumscribed circles
  of the two triangles incident with $e$ intersect at the angle $\theta_e$;
  see figure~\ref{fig:leibon_angles} {\it (right)}.
\end{proof}

\chapter{Circle patterns and the volumes of hyperbolic polyhedra}
\label{cha:volume}

In chapter~\ref{cha:functionals}, we constructed circle patterns by tiling
(euclidean, hyperbolic or spherical) $2$-space with kite-shaped
quadrilaterals. In this chapter, we extend these tilings to `tilings' of
hyperbolic $3$-space with $3$-dimensional polyhedral tiles. This provides a
geometric interpretation of the functional $\Shat$ of
section~\ref{sec:legendre}. It turns out to be the sum of the volumes of the
polyhedral tiles. The functional attains its maximum (under variation
constrained to coherent angle systems) if the tiles fit together to form a
$3$-dimensional polyhedron. The maximal value is therefore the volume of this
polyhedron.

\section{Schl\"afli's differential volume formula}
\label{sec:schlafli}

In 1852, Schl\"afli discovered the remarkable formula \eqref{eq:schlafli_n}
for the differential of the volume of a spherical simplex of arbitrary
dimension \cite{schlafli1852}. In 1936, Kneser gave an elegant proof (in a
nasty journal) \cite{kneser36}. It works also for hyperbolic simplices.

\begin{theorem}[Schl\"afli's differential volume formula]
  For an $n$-dimensional simplex $S$ in spherical or hyperbolic space, let
  $V(S)$ be its $n$-dimensional volume. Denote the $(n-1)$-faces by $S_i$,
  $0\leq i\leq n$. Let $a_{jk}$ be the $(n-2)$-dimensional volume of the
  $(n-2)$-face $S_j\cap S_k$, and let $\alpha_{jk}$ be the the dihedral angle
  between $S_j$ and $S_k$. Then the differential of the volume function $V$
  is
  \begin{equation}
    \label{eq:schlafli_n}
    dV = \frac{\varepsilon}{n-1}\sum_{0\leq j<k\leq{}n} a_{jk}\,d\alpha_{jk},
  \end{equation}
  where $\varepsilon=1$ in the spherical and $\varepsilon=-1$ in the
  hyperbolic case.
\end{theorem}

For $3$-dimensional hyperbolic simplices (tetrahedra), formula
\eqref{eq:schlafli_n} may be written
\begin{equation}
  \label{eq:schlafli}
  dV = -\frac{1}{2}\sum_{\text{edges } j} a_j\,d\alpha_j,
\end{equation}
where the sum is taken over the edges $j$, and $a_j$, $\alpha_j$
are the length and dihedral angle of edge $j$. In fact, this formula holds
for arbitrary polyhedral shapes:

\begin{corollary}
  Formula \eqref{eq:schlafli} holds for arbitrary compact hyperbolic
  $3$-manifolds with polyhedral boundary (in particular, for $3$-dimensional
  hyperbolic polyhedra), when they are deformed without changing the
  combinatorial type of the polyhedral boundary.
\end{corollary}

This follows by triangulating the $3$-manifold with polyhedral boundary. At all
edges of the triangulation, except those which are contained in edges of the
polyhedral boundary, the sum of dihedral angles is constant.

Now consider hyperbolic polyhedra (or manifolds with polyhedral boundary)
with some vertices in the infinite boundary of $H^3$. Even though they are
not compact, their volume is finite. However, formula \eqref{eq:schlafli}
does not apply directly, because some edges have infinite length. The
following lemma is ascribed to Milnor
\cite[p.~576]{rivin94}.\footnote{However, the proof given by Rivin does not
  seem to be valid, because, in the ``slicing a carrot" argument, Schl\"afli's
  formula is applied to shapes bounded by planes and horospheres.}

\begin{lemma}[Milnor]
  \label{lem:milnor}
  Suppose $P$ is a $3$-dimensional hyperbolic polyhedron with some vertices
  on the infinite boundary. For each infinite vertex of $P$, choose an
  arbitrary horosphere centered at this vertex. For a finite edge $j$ of $P$,
  let $a_j$ be its length. For an edge $j$ with one (or two) vertices at
  infinity, let $a_j$ be its length up to the cut-off point(s) at the
  corresponding horosphere(s). (This may be negative if the horospheres are
  not small enough.) With this agreement, the following holds: If $P$ is
  deformed in such a way that its combinatorial type does not change and the
  infinite vertices stay infinite, then the volume differential is given by
  formula \eqref{eq:schlafli}.
\end{lemma}

\begin{proof}
  Because the sum of dihedral angles at an infinite vertex with $n$\/ edges is
  constantly $(n-2)\pi$\/ during the deformation, the right hand side of
  equation \eqref{eq:schlafli} is independent of the choice of horospheres.
  We may thus assume that the horospheres are small enough, so that each
  intersects only the planes incident with the vertex at which it is centered
  and no other planes or horospheres. By triangulating the polyhedron
  appropriately, one finds that it suffices to consider a triply orthogonal
  tetrahedron with one vertex at infinity, see figure~\ref{fig:orthoscheme}
  {\em(left)}. The volume of such a tetrahedron can be found by a
  straightforward integration without using Schl\"afli's formula, as shown in
  appendix~\ref{app:triply_ortho_tetra}.  Differentiate
  equation~\eqref{eq:vol_orthoscheme} to obtain, after a straightforward
  calculation,
  \begin{equation*}
    \frac{\partial V}{\partial\alpha}=
    -\frac{1}{2}\log\sqrt{1-\frac{\cos^2\beta}{\cos^2\alpha}}.
  \end{equation*}
  Truncate the tetrahedron by the horosphere centered at the infinite vertex
  which touches the opposite face. Then the truncated sides with dihedral
  angles $\alpha$ and $\frac{\pi}{2}-\alpha$ have length $a=0$ and
  $\tilde{a}=-\log\sqrt{1-\frac{\cos^2\beta}{\cos^2\alpha}}$, respectively.
  Hence,
  \begin{equation*}
    dV\big|_{\beta=\text{\it const.}} = 
    -\frac{1}{2}\big(a\,d\alpha + \tilde{a}\,d(\pi-\alpha)\big).
  \end{equation*}
  This proves the theorem.
\end{proof}

\section[A prototypical variational principle]{A prototypical variational principle and its Legendre dual}
\label{sec:proto_principle}

In this section, we take a little detour to illustrate the fundamental ideas
which are applied in the next sections. It may be skipped, or read after the
following sections. We derive a pair of Legendre dual variational principles
connected with triangulations of hyperbolic $3$-manifolds with polyhedral
boundary. 
The arguments below may be adapted to deal with the case where all vertices
of the triangulation are at infinity~\cite{schlenker02}. In that case, one
does obtain convex functionals.

Let $\mathcal{T}$ be a finite triangulation of a compact topological
$3$-manifold $M$\/ with boundary, and let
$\mathcal{T}_0,\ldots,\mathcal{T}_3$ be the sets of vertices, edges,
triangles, and tetrahedra. Let $\mathcal{S}$ be the manifold of shapes of
hyperbolic tetrahedra. This is diffeomorphic to a connected open subset of
$\mathds R^6$. Global coordinates are the six dihedral angles (or,
alternatively, the six edge lengths), which satisfy certain inequalities.

For a function $\sigma:\mathcal{T}_3\rightarrow\mathcal{S}$, which assigns a
shape to each combinatorial tetrahedron of $\mathcal{T}$ independently, let
\begin{equation}
  \label{eq:proto_func}
  \ShatProto(\sigma) = 2 \sum_{t\in\mathcal{T}_3} V(\sigma(t)),
\end{equation}
where $V$ is the volume function on $\mathcal{S}$. This defines a
differentiable functional
$\ShatProto:\mathcal{S}^{\mathcal{T}_3}\rightarrow\mathds R$. 

Let $\Phi:\mathcal{T}_1\cap\partial M\rightarrow\mathds R_{+}$ be a function
that assigns a positive angle to each boundary edge.  A shape assignment
$\sigma:\mathcal{T}_3\rightarrow\mathcal{S}$ is {\em coherent}, if the
dihedral angles sum to $2\pi$ around each interior edge and to
$\Phi(e)$ around each boundary edge $e$.

\begin{proposition}
  A coherent shape assignment $\sigma\in\mathcal{S}^{\mathcal{T}_3}$ is a
  critical point of the functional $\ShatProto$ under variations in the space
  of coherent shape assignments if and only if the hyperbolic tetrahedra
  $\sigma(\mathcal{T}_3)$ fit together to form a hyperbolic manifold with
  polyhedral boundary, where the boundary angles are prescribed by $\Phi$.
\end{proposition}

\begin{proof}
  Let $\sigma$ be a coherent shape assignment. Suppose
  $t_1,t_2\in\mathcal{T}_3$ share an edge $e\in\mathcal{T}_1$. Let
  $\alpha_1$, $\alpha_2$ be the dihedral angles of $\sigma(t_1),\sigma(t_2)$
  at $e$, and let $a_1, a_2$ be the side lengths of $e$ in
  $\sigma(t_1),\sigma(t_2)$. There is a variation $v$ in the space of
  coherent shape assignments with $1=d\alpha_1(v)=-d\alpha_2(v)$ which keeps
  all other dihedral angles constant. By Schl\"afli's
  formula~\eqref{eq:schlafli},
  \begin{equation*}
    d\ShatProto(v)=-a_1+a_2.
  \end{equation*}
  Hence, if $\sigma$ is a critical point of $\ShatProto$ under variations in
  the space of coherent shape assignments, then $\sigma$ assigns a unique
  length to all edges $e\in\mathcal{T}_1$. Therefore, the tetrahedra fit
  together.

  The converse statement follows, because variations like $v$ span the
  tangent space to the space of coherent shape assignments.
\end{proof}

To define the Legendre dual functional, consider a different subset of shape
assignments: the space $\mathcal{R}\subset\mathcal{S}^{\mathcal{T}_3}$ of
those $\sigma\in\mathcal{S}^{\mathcal{T}_3}$ which assign a unique length to
each edge $e\in\mathcal{T}_1$. This means, if $t_1,t_2\in\mathcal{T}_3$ are
both incident with $e$, then the length of $e$ in $\sigma(t_1)$ equals the
length of $e$ in $\sigma(t_2)$. For each such shape assignment
$\sigma\in\mathcal{R}$ and each edge $e\in\mathcal{T}_1$, let $a_e(\sigma)$
be the length assigned to $e$ by $\sigma$. The functions
$a_e:\mathcal{R}\rightarrow\mathds R_+$ (which satisfy certain inequalities)
form a system of coordinates of $\mathcal{R}$. Hence, $\mathcal{R}$ is
diffeomorphic to an open subset of $\mathds{R}^{\mathcal{T}_1}$. For a
tetrahedron $t\in\mathcal{T}_3$ with edge $e\in t\cap\mathcal{T}_1$, let
$\alpha_{t,e}(\sigma)$ be the dihedral angle at $e$ in $\sigma(t)$. Extend
the function $\Phi$, which was defined on the boundary edges, to a function
on all edges by setting $\Phi(e)=2\pi$ if $e$ is an interior edge.

Define the functional $\Sproto:\mathcal{R}\rightarrow\mathds
R$,
\begin{equation*}
  \Sproto(\sigma) = 
  -\sum_{t\in\mathcal{T}_3} 
  \Big(\sum_{e\in t\cap\mathcal{T}_1}
  a_e(\sigma)\alpha_{t,e}(\sigma) + 2V(\sigma(t))\Big)
  + \sum_{e\in\mathcal{T}_1}\Phi(e)a_e(\sigma).
\end{equation*}

\begin{proposition}
  A shape assignment $\sigma\in\mathcal{R}$ is a critical point of the
  functional $\Sproto$ (under variations in $\mathcal{R}$) if and only if the
  tetrahedra $\sigma(\mathcal{T}_3)$ fit together to form a hyperbolic
  manifold with polyhedral boundary, where the boundary angles are prescribed
  by $\Phi$.
\end{proposition}

\begin{proof} 
  By Schl\"afli's differential volume formula~\eqref{eq:schlafli},
  \begin{equation*}
    d\Sproto=\sum_{e\in\mathcal{T}_1}\big(\Phi(e)-
    \sum_{
      \begin{smallmatrix}
        t\in\mathcal{T}_3:\\
        t\ni e
      \end{smallmatrix}
    }
    \alpha_{t,e}\big)\,da_e.
  \end{equation*}
  Thus, at a critical point, the sum of all dihedral angles at an edge $e$ is
  $\Phi(e)$.
\end{proof}

\section{The euclidean functional}
\label{sec:volume_euc}

\begin{figure}[t]%
\centering%
\input{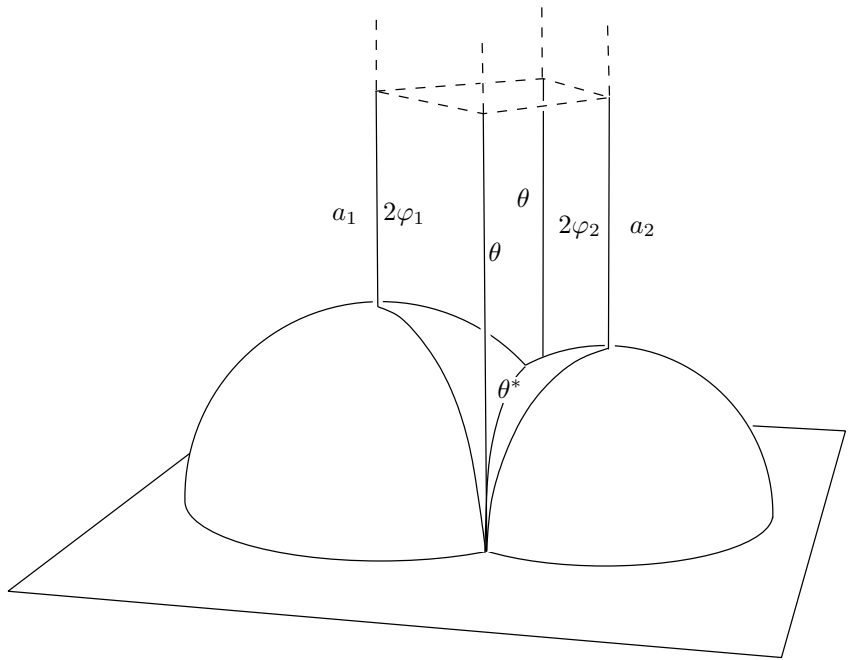}%
\caption{The hyperbolic polyhedron $P_{\theta}(\varphi_1,\varphi_2)$ for
  $\varphi_1+\varphi_2=\theta^*$ (where $\theta^*=\pi-\theta$) shown in the
  Poincar\'e half-space model. The dihedral angle at the four unmarked edges
  is $\frac{\pi}{2}$.}%
\label{fig:quadpoly_euc}
\end{figure}

Suppose that $0<\theta<\pi$ and $\varphi_1>0$, $\varphi_2>0$. If
\begin{equation*}
  \varphi_1+\varphi_2=\pi-\theta,
\end{equation*}
we define a hyperbolic polyhedron $P_{\theta}(\varphi_1,\varphi_2)$ as
follows. 

If $\varphi_1$ and $\varphi_2$ are not greater than $\frac{\pi}{2}$, let
$P_{\theta}(\varphi_1,\varphi_2)$ be the hyperbolic polyhedron shown in
figure~\ref{fig:quadpoly_euc} in the Poincar\'e half-space model of
hyperbolic $3$-space: $H^3=\{(x,y,z)\in\R^3\,|\,z>0\}$ with metric
$|ds|=\frac{1}{z}\sqrt{x^2+y^2+z^2}$. It has three infinite vertices: two lie
in the plane $z=0$ and third is the infinite point compactifying the closed
half-space. The dihedral angles are as indicated. The unmarked edges have
dihedral angle $\frac{\pi}{2}$.

\begin{figure}[t]
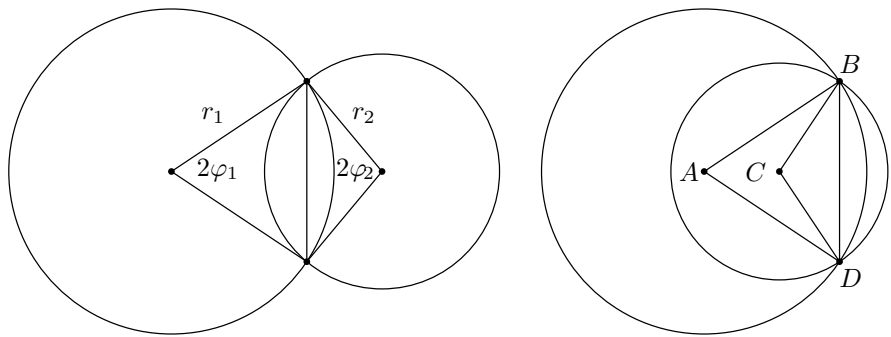
%
\hfill%
\input{quad_convex.tex}%
\hfill%
\input{quad_nonconv.tex}%
\hspace*{\fill}%
\caption{{\em Left:}\/~The polyhedron of figure~\ref{fig:quadpoly_euc},
  projected orthogonally into the $z$-plane. {\em Right:}\/~The corresponding
  picture for the case
  $\varphi_2>\frac{\pi}{2}$.}%
\label{fig:quad_convex_nonconvex}%
\end{figure}

If either $\varphi_1$ or $\varphi_2$ are greater than $\frac{\pi}{2}$, we
define $P_{\theta}(\varphi_1,\varphi_2)$ as an algebraic sum.
Figure~\ref{fig:quad_convex_nonconvex} {\em(left)}\/ shows the orthogonal
projection of figure~\ref{fig:quadpoly_euc} into the $z$-plane. The
corresponding figure for the case $\varphi_2>\frac{\pi}{2}$ is shown on the
right. In this case, let $P_1$ be the polyhedron bounded by the planes
corresponding to the circle around $A$ and the line segments $AB$, $BD$, and
$DA$; and let $P_2$ be the polyhedron bounded by the planes corresponding to
the circle around $C$ and the line segments $CB$, $BD$, and $DC$. Define
$P_{\theta}(\varphi_1,\varphi_2)$ as the algebraic sum
$P_{\theta}(\varphi_1,\varphi_2)=P_1-P_2$. Extend the volume function linearly
to such algebraic sums: $V(P_1-P_2)=V(P_1)-V(P_2)$.

\begin{proposition}
  Suppose $\varphi\in\R^{\vecEint}$ is a euclidean coherent angle system (see
  section~\ref{sec:coherent_angle_systems}). Then the functional
  $\Shat(\varphi)$ of theorem~\ref{thm:Shat} satisfies
  \begin{equation*}
    \Shat(\varphi)=\sum_{e\in E_{\text{int}}} 
    \Big(2V\big(P_{\theta_e}(\varphi_{\vece},\varphi_{-\vece})\big)
    -\Cl(2\theta^*_e)\Big),
  \end{equation*}
  where $V$ is the volume function, the sum is taken over all non-oriented
  interior edges $e$, and $\vece$, $-\vece$ are the two corresponding
  oriented edges.
\end{proposition}

\begin{proof}
Because of equation~\eqref{eq:ShatEuc} for $\Shat(\varphi)$, where $\varphi$
is a euclidean coherent angle system, the proposition follows directly from
the following lemma. 
\end{proof}

\begin{lemma}
\label{lem:VP_euc}
  The volume of $P_{\theta}(\varphi_1,\varphi_2)$ is 
  \begin{equation}
    \label{eq:VP_euc}
    V\big(P_{\theta}(\varphi_1,\varphi_2)\big)=
    \frac{1}{2}\Cl(2\varphi_1)+\frac{1}{2}\Cl(2\varphi_2).
  \end{equation}
\end{lemma}

\begin{proof}
  The polyhedron $P_{\theta}(\varphi_1,\varphi_2)$ may be decomposed into four
  triply orthogonal tetrahedra (see appendix~\ref{app:triply_ortho_tetra}),
  each having two vertices at infinity. The characteristic angle is
  $\varphi_1$ for two of the tetrahedra and $\varphi_2$ for the other two.
  Equation~\eqref{eq:VP_euc} follows from equation~\eqref{eq:Valpha}. (This
  argument works also in the case where $\varphi_1$ or $\varphi_2$ are
  greater than $\frac{\pi}{2}$.)
\end{proof}

\begin{corollary}
  If $\varphi\in\R^{\vecEint}$ is a critical point of $\Shat$ under variations
  in the space of euclidean coherent angle systems, then the polyhedra
  $P_{\theta_e}(\varphi_{\vece},\varphi_{-\vece})$ fit together to form a
  hyperbolic orbifold $M$ with polyhedral boundary. $M$ is a manifold with
  polyhedral boundary if\/ $\Phi_f=2\pi$ for all interior faces $f\in F$. The
  volume of $M$ is 
  \begin{equation*}
    V(M)
    =\frac{1}{2}\Big(\Shat(\varphi)
    +\sum_{e\in E_{\text{int}}}\Cl(2\theta^*_e)\Big)
    =\frac{1}{2}\Big(\Seuc(\rho)
    +\sum_{e\in E_{\text{int}}}\Cl(2\theta^*_e)\Big) ,
  \end{equation*}
  where $\rho\in\R^F$ is the corresponding critical point of $\Seuc$.
\end{corollary}

\section{The spherical functional}
\label{sec:volume_sph}

\begin{figure}[tbp]%
\centering%
\input{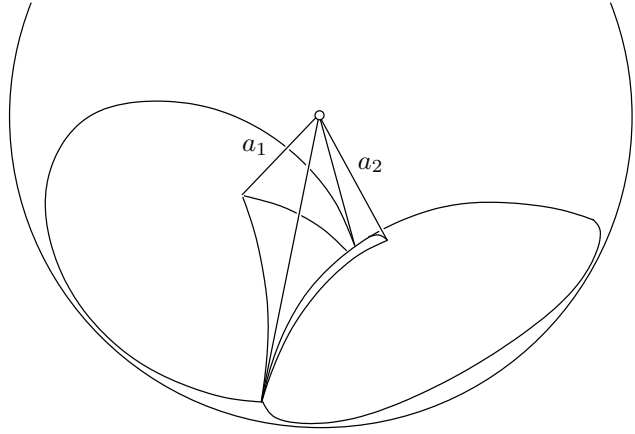}%
\caption{The hyperbolic polyhedron $P_{\theta}(\varphi_1,\varphi_2)$ for
  $\pi-\theta<\varphi_1+\varphi_2<\pi+\theta$ and
  $\big|\varphi_1-\varphi_2\big|<\pi-\theta$, shown in the Poincar\'e ball
  model. The dihedral angles at the edges marked $a_1$ and $a_2$ are
  $2\varphi_1$ and $2\varphi_2$. The other two edges of the vertex marked
  `$\circ$' (which is at $(x,y,z)=(0,0,0)$) have angle $\theta$; and the
  faces opposite vertex `$\circ$' intersect at the angle
  $\theta^*=\pi-\theta$. The four remaining edges have dihedral angle
  $\frac{\pi}{2}$.  The polyhedron
  has two vertices in the infinite boundary.}%
\label{fig:quadpoly_sph}%
\end{figure}

In the previous section, we defined $P_{\theta}(\varphi_1,\varphi_2)$ for
$\varphi_1+\varphi_2=\pi-\theta$. If
\begin{equation*}
  \pi-\theta<\varphi_1+\varphi_2<\pi+\theta
  \quad\text{and}\quad
  \big|\varphi_1-\varphi_2\big|<\pi-\theta,
\end{equation*}
let $P_{\theta}(\varphi_1,\varphi_2)$ be the hyperbolic polyhedron shown in
figure~\ref{fig:quadpoly_sph} 
in the Poincar\'e ball model of hyperbolic space:
\begin{equation*}
  H^3=\big\{(x,y,z)\in\R^3\;\big|\;x^2+y^2+z^2<1\big\}
\end{equation*}
with metric 
\begin{equation*}
  |ds|=\frac{2}{1-x^2-y^2-z^2}\sqrt{dx^2+dy^2+dz^2}.
\end{equation*}
(If $\varphi_1$ or $\varphi_2$ are greater
than $\frac{\pi}{2}$, then $P_{\theta}(\varphi_1,\varphi_2)$ is defined by an
algebraic sum, as in the previous section. We will not dwell on this point.)

\begin{proposition}
  Suppose $\varphi\in\R^{\vecEint}$ is a spherical coherent angle system (see
  section~\ref{sec:coherent_angle_systems}). Then the functional
  $\Shat(\varphi)$ of theorem~\ref{thm:Shat} satisfies
  \begin{equation*}
    \Shat(\varphi)=\sum_{e\in E_{\text{int}}} 
    2V\big(P_{\theta_e}(\varphi_{\vece},\varphi_{-\vece})\big)
  \end{equation*}
  where $V$ is the volume function, the sum is taken over all non-oriented
  interior edges $e$, and $\vece$, $-\vece$ are the two corresponding
  oriented edges.
\end{proposition}

\begin{corollary}
  If $\varphi\in\R^{\vecEint}$ is a critical point of $\Shat$ under variations
  in the space of spherical coherent angle systems, then the polyhedra
  $P_{\theta_e}(\varphi_{\vece},\varphi_{-\vece})$ fit together to form a
  hyperbolic orbifold $M$ with polyhedral boundary. $M$ is a manifold with
  polyhedral boundary if\/ $\Phi_f=2\pi$ for all interior faces $f\in F$. The
  volume of $M$ is
  \begin{equation*}
    V(M)
    ={\textstyle\frac{1}{2}}\Shat(\varphi)
    ={\textstyle\frac{1}{2}}\Ssph(\rho),
  \end{equation*}
  where $\rho\in\R^F$ is the corresponding critical point of $\Ssph$.
\end{corollary}

The proposition follows directly from the following lemma.

\begin{lemma}
\label{lem:VP_sph}
  The volume of $P_{\theta}(\varphi_1,\varphi_2)$ is 
  \begin{equation}
    \label{eq:VP_sph}
    \begin{split}
      V\big(P_{\theta}(\varphi_1,\varphi_2)\big)
      =\frac{1}{2}\Big(& 
        \Cl(\theta^* + \varphi_{1}-\varphi_{2})
      + \Cl(\theta^* - \varphi_{1}+\varphi_{2}) \\
      +&\Cl(\theta^* + \varphi_{1}+\varphi_{2}) 
      + \Cl(\theta^* - \varphi_{1}-\varphi_{2}) -2\Cl(2\theta^*) \Big).
    \end{split}
  \end{equation}
\end{lemma}

\begin{proof}[Proof of Lemma~\ref{lem:VP_sph}]
  We will show that the derivatives of both sides of
  equation~\eqref{eq:VP_sph} with respect to $\varphi_1$ and $\varphi_2$ are
  equal. Since both sides tend to $0$ as $\varphi_1\rightarrow0$ and
  $\varphi_2\rightarrow\theta^*$, equation~\eqref{eq:VP_sph} follows.
  
  Consider deformations of $P_{\theta}(\varphi_1,\varphi_2)$ during which all
  dihedral angles except $2\varphi_1$ and $2\varphi_2$ remain constant. By
  Milnor's generalization of Schl\"afli's differential volume formula
  (lemma~\ref{lem:milnor}),
   \begin{equation*}
     dV\big(P_{\theta}(\varphi_1,\varphi_2)\big)
     =-a_1\,d\varphi_1-a_2\,d\varphi_2,
   \end{equation*}
   where $a_1$, $a_2$ are the lengths of the sides with dihedral angles
   $2\varphi_1$ and $2\varphi_2$.  (And this is also true in the case where
   $P_{\theta}(\varphi_1,\varphi_2)$ is defined as an algebraic sum.) By
   equation~\eqref{eq:clausen_diff}, the differential of the right hand side
   of equation~\eqref{eq:VP_sph} is
   \begin{equation*}
     \rho_1\,d\varphi_1+\rho_2\,d\varphi_2,
   \end{equation*}
   where $\rho_1$ and $\rho_2$ are given by
   equation~\eqref{eq:rho_f_of_phi_e_sph}. Let $r_1$ and $r_2$ be the {\em
   spherical}\/ radii of the circles in the infinite boundary shown in
   figure~\ref{fig:quadpoly_sph}. Then $r_1$ and $r_2$ are related to
   $\rho_1$ and $\rho_2$ by equation~\eqref{eq:rho_sph}. This implies
   $\rho_1=-a_1$ and $\rho_2=-a_2$; see figure~\ref{fig:r_and_a_sph}.
\begin{figure}[tb]%
\centering%
\input{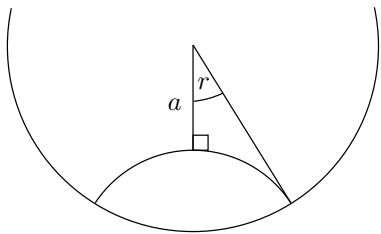}%
\caption{In a hyperbolic right-angled triangle with one vertex at 
  infinity, the length $a$ and the angle $r$ satisfy the equation 
  $a=-\log\tan\frac{r}{2}$.}%
\label{fig:r_and_a_sph}%
\end{figure}
This completes the proof.
\end{proof}

\section{The hyperbolic functional}
\label{sec:volume_hyp}

\begin{figure}[tbp]%
\centering%
\input{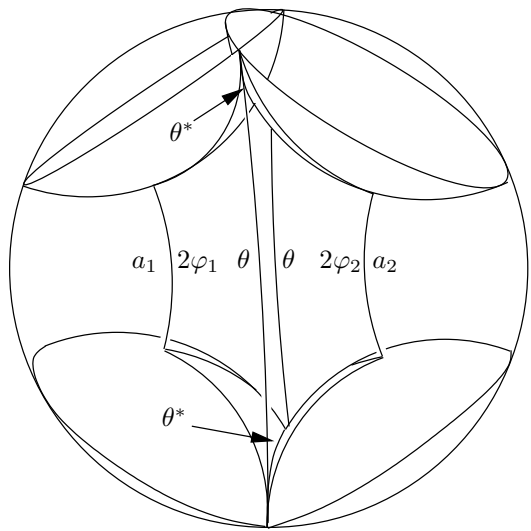}%
\caption{The hyperbolic polyhedron $P_{\theta}(\varphi_1,\varphi_2)$ for
  $\varphi_1+\varphi_2<\pi-\theta$, shown in the Poincar\'e ball model. The
  dihedral angles are as indicated, where $\theta^*=\pi-\theta$. The unmarked
  edges have dihedral angle $\frac{\pi}{2}$. The polyhedron has four vertices
  in the infinite boundary. It is symmetric with respect to reflection on the
  equatorial
plane.}%
\label{fig:quadpoly_hyp}%
\end{figure}

Now we define the hyperbolic polyhedron $P_{\theta}(\varphi_1,\varphi_2)$ for
\begin{equation*}
  \varphi_1+\varphi_2<\pi-\theta.
\end{equation*}
In this case, let $P_{\theta}(\varphi_1,\varphi_2)$ be the hyperbolic
polyhedron shown in figure~\ref{fig:quadpoly_hyp} in the Poincar\'e ball
model of hyperbolic space.
(If $\varphi_1$ or $\varphi_2$ are greater than $\frac{\pi}{2}$, then
$P_{\theta}(\varphi_1,\varphi_2)$ is defined by an algebraic sum, as
described in section~\ref{sec:volume_euc}. We will not dwell on this point.)

\begin{proposition}
  Suppose $\varphi\in\R^{\vecEint}$ is a hyperbolic coherent angle system (see
  section~\ref{sec:coherent_angle_systems}). Then the functional
  $\Shat(\varphi)$ of theorem~\ref{thm:Shat} satisfies
  \begin{equation*}
    \Shat(\varphi)=\sum_{e\in E_{\text{int}}} 
    V\big(P_{\theta_e}(\varphi_{\vece},\varphi_{-\vece})\big)
  \end{equation*}
  where $V$ is the volume function, the sum is taken over all non-oriented
  interior edges $e$, and $\vece$, $-\vece$ are the two corresponding
  oriented edges.
\end{proposition}

\begin{corollary}
  If $\varphi\in\R^{\vecEint}$ is a critical point of $\Shat$ under variations
  in the space of hyperbolic coherent angle systems, then the polyhedra
  $P_{\theta_e}(\varphi_{\vece},\varphi_{-\vece})$ fit together to form a
  hyperbolic orbifold $M$ with polyhedral boundary. $M$ is a manifold with
  polyhedral boundary if\/ $\Phi_f=2\pi$ for all interior faces $f\in F$. The
  volume of $M$ is
  \begin{equation*}
    V(M)
    =\Shat(\varphi)
    =\Shyp(\rho),
  \end{equation*}
  where $\rho\in\R^F$ is the corresponding critical point of $\Shyp$.
\end{corollary}

The proposition follows directly from the following lemma.

\begin{lemma}
\label{lem:VP_hyp}
  The volume of $P_{\theta}(\varphi_1,\varphi_2)$
  is 
  \begin{equation}
    \label{eq:VP_hyp}
    \begin{split}
      V\big(P_{\theta}(\varphi_1,\varphi_2)\big)
      =& \Cl(\theta^* + \varphi_{1}-\varphi_{2})
      + \Cl(\theta^* - \varphi_{1}+\varphi_{2}) \\
      +&\Cl(\theta^* + \varphi_{1}+\varphi_{2}) 
      + \Cl(\theta^* - \varphi_{1}-\varphi_{2}) -2\Cl(2\theta^*).
    \end{split}
  \end{equation}
\end{lemma}

Lemma~\ref{lem:VP_hyp} is proved in the same way as lemma~\ref{lem:VP_sph}:

\begin{proof}[Proof of Lemma~\ref{lem:VP_hyp}]
  We will show that the derivatives of both sides of
  equation~\eqref{eq:VP_hyp} with respect to $\varphi_1$ and $\varphi_2$ are
  equal. Since both sides tend to $0$ as $\varphi_1\rightarrow0$ and
  $\varphi_2\rightarrow\theta^*$, equation~\eqref{eq:VP_sph} follows.
  
  Under deformations of $P_{\theta}(\varphi_1,\varphi_2)$ during which all
  dihedral angles except $2\varphi_1$ and $2\varphi_2$ remain constant, the
  volume differential is
   \begin{equation*}
     dV\big(P_{\theta}(\varphi_1,\varphi_2)\big)
     =-a_1\,d\varphi_1-a_2\,d\varphi_2,
   \end{equation*}
   where $a_1$, $a_2$ are the lengths of the sides with dihedral angles
   $2\varphi_1$ and $2\varphi_2$.  (And this is also true in the case where
   $P_{\theta}(\varphi_1,\varphi_2)$ is defined as an algebraic sum.) By
   equation~\eqref{eq:clausen_diff}, the differential of the right hand side
   of equation~\eqref{eq:VP_hyp} is
   \begin{equation*}
     2\rho_1\,d\varphi_1+2\rho_2\,d\varphi_2,
   \end{equation*}
   where $\rho_1$ and $\rho_2$ are given by
   equation~\eqref{eq:rho_f_of_phi_e_hyp}. The symmetry plane intersects
   $P_{\theta}(\varphi_1,\varphi_2)$ in a hyperbolic kite-shaped
   quadrilateral with angles $\theta$, $\varphi_1$, $\varphi_2$. Let the
   hyperbolic side lengths of this quadrilateral be $r_1$ and $r_2$. Then
   $r_1$ and $r_2$ are related to $\rho_1$ and $\rho_2$ by
   equation~\eqref{eq:rho_hyp} (see lemma~\ref{lem:phi_hyp}).
\begin{figure}[tb]%
\centering%
\input{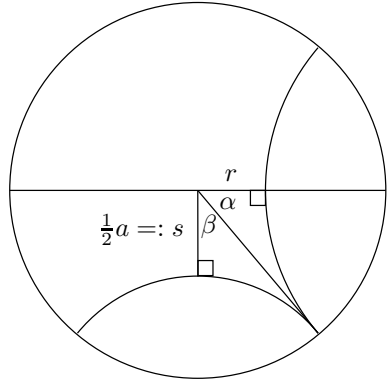}%
\caption{From $r=-\log\tan\frac{\alpha}{2}$ and
   $s=-\log\tan\frac{\beta}{2}$ with $\alpha+\beta=\frac{\pi}{2}$,
   one obtains $s=-\log\tanh\frac{r}{2}$. Hence, $s=-\rho$, with $\rho$ given
   by equation~\eqref{eq:rho_hyp}.}%
\label{fig:r_and_a_hyp}%
\end{figure}
Furthermore, $a_j=-2\rho_j$ (see figure~\ref{fig:r_and_a_hyp}). This completes
the proof.
\end{proof}

\section{Leibon's functional}
\label{sec:leibon_volume}

Schl\"afli's differential volume formula is the basis for Leibon's functional
(see section~\ref{sec:leibon_func}) as well. Lemma~\ref{lem:leibon_volume}
and Schl\"afli's differential volume formula provide another proof of Leibon's
variational principle, lemma~\ref{lem:leibon_principle}.

\begin{figure}%
\input{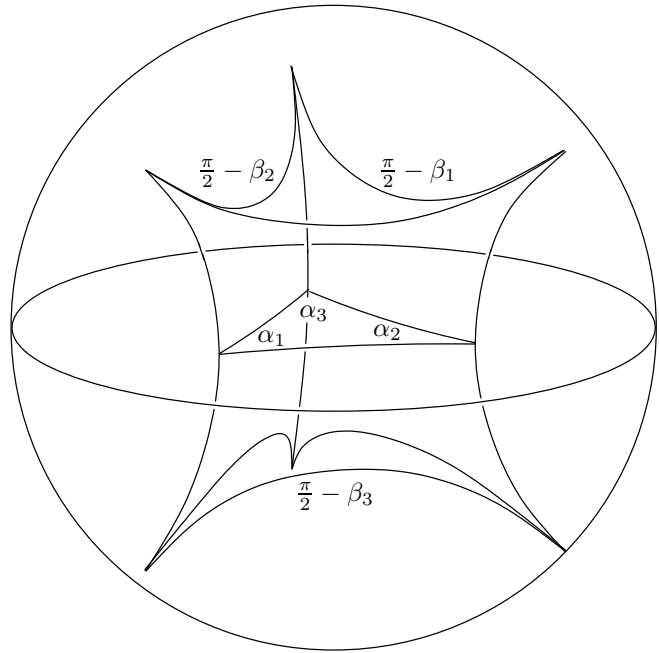}%
\caption{An ideal triangular prism, shown in the Poincar\'e ball model. The
  equatorial plane is the plane of symmetry.}%
\label{fig:ideal_prism}%
\end{figure}

An {\em ideal triangular prism}\/ is a convex hyperbolic polyhedron with six
ideal vertices, which is symmetric with respect to a hyperbolic plane (see
figure~\ref{fig:ideal_prism}).
An ideal triangular prism intersects the plane of symmetry in a hyperbolic
triangle. The angles $\alpha_1$, $\alpha_2$, $\alpha_3$ of this triangle
determine the shape of the prism. They are also the dihedral angles of the
three edges intersecting the symmetry plane. The other dihedral angles are
$\frac{\pi}{2}-\beta_1$, $\frac{\pi}{2}-\beta_2$, and
$\frac{\pi}{2}-\beta_3$, with the $\beta_j$ defined by
equations~\eqref{eq:beta_123}.

\begin{figure}%
\input{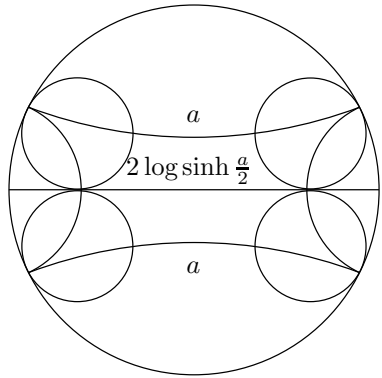}%
\caption{Truncated quadrilateral face of the ideal prism.}%
\label{fig:truncated_quadrilateral}%
\end{figure}

\begin{lemma}[Leibon~\cite{leibon99}\cite{leibon02}]
  \label{lem:leibon_volume}
  The hyperbolic volume of the ideal prism shown in
  figure~\ref{fig:ideal_prism} is $V(\alpha_1,\alpha_2,\alpha_3)$ as defined
  by equation~\eqref{eq:V_alpha_123}. 
  
  Truncate the ideal prism by horospheres which are centered at the ideal
  vertices and touch the symmetry plane. The truncated edge lengths are $0$
  for the edges which intersect the symmetry plane, and
  $2\log\sinh\frac{a_j}{2}$ for the other edges, where $a_j$ is the length of
  the corresponding side of the hyperbolic triangle with angles $\alpha_1$,
  $\alpha_2$, $\alpha_3$. (See figure~\ref{fig:truncated_quadrilateral}.)
\end{lemma}

\section[A common ancestor]{A common ancestor of Leibon's and our functionals}
\label{sec:common_ancestor}

The geometric interpretation of both Leibon's functional and the functional
$\Shat(\varphi)$ in terms of hyperbolic volume clarifies their mutual
relationship. In this section, we construct another functional,
$T(\varphi,\beta)$, from which both Leibon's functional and $\Shat(\varphi)$
can be derived. Let $\Sigma$ be a cell decomposition of a compact surface
without boundary, let intersection angles be prescribed by a function
$\theta:E\rightarrow(0,\pi)$ on the non-oriented edges, and let cone angles
in the centers of the circles be prescribed by a function
$\Phi:F\rightarrow(0,2\pi)$ on the faces.  

\begin{figure}
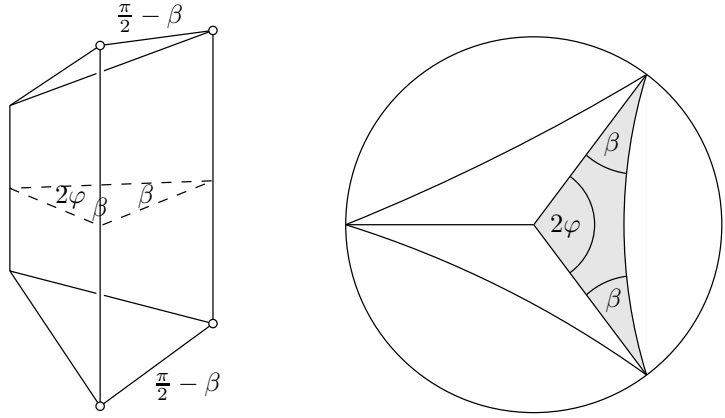
%
\hfill%
\input{decomposed_prism.tex}%
\hfill%
\input{decomposed_triangle.tex}%
\hspace*{\fill}%
\caption{{\it Left:} Schematic picture of one of the three pieces into which
  the ideal prism of figure~\ref{fig:ideal_prism} is decomposed. It has four
  ideal vertices, which are marked~`$\circ$'. The dihedral angles at four of
  the edges are $\frac{\pi}{2}$. The other dihedral angles are as indicated.
  {\it Right:} The triangular intersection with the symmetry plane is
  decomposed into three isosceles triangles.}%
\label{fig:decompose_prism}%
\end{figure}

The ideal prism shown in figure~\ref{fig:ideal_prism} can be decomposed into
three polyhedral pieces like the one shown in
figure~\ref{fig:decompose_prism}.
The hyperbolic polyhedron $P_{\theta}(\varphi_1,\varphi_2)$ (see
figure~\ref{fig:quadpoly_hyp}) can be decomposed into two pieces of this
kind. Let $W(\varphi,\beta)$ be the hyperbolic volume of such a polyhedron.

A {\it coherent angle system} in this setting is a pair
$\varphi\in\R^{\vecE}$, $\beta\in\R^{\vecE}$ of functions on the oriented
edges satisfying

\makebox[2em][l]{\it (i)} For all $\vece\in\vecE$,
\begin{equation}
  \label{eq:phi_beta_inequality}
  \varphi_{\vece}>0, \quad 
  \beta_{\vece}>0, \quad \text{and}\quad
  \varphi_{\vece}+\beta_{\vece}<\frac{\pi}{2}.
\end{equation}

\nopagebreak 
\makebox[2em][l]{\it (ii)} For all faces $f\in F$,
\begin{equation}
  \label{eq:phi_constraint}
  \sum 2\varphi_{\vece}=\Phi_f,
\end{equation}
where the sum is taken over all oriented edges in the oriented boundary of
$f$.
    
\nopagebreak
\makebox[2em][l]{\it (iii)} For all $\vece\in\vecE$,
\begin{equation}
  \label{eq:beta_constraint}
  \beta_{\vece}+\beta_{-\vece}=\theta_{e}.
\end{equation}

For $\varphi\in\R^{\vecE}$, $\beta\in\R^{\vecE}$ which satisfy {\it (i)},
define the functional
\begin{equation*}
  T(\varphi,\beta)=\sum_{\vece\in\vecE}W(\varphi_{\vece},\beta_{\vece}).
\end{equation*}
(Condition {\it (i)}\/ has to be fulfilled for the corresponding polyhedron
of figure~\ref{fig:decompose_prism} {\it (left)}\/ to exist.)

\begin{proposition}
  A coherent angle system $\varphi$, $\beta$ is a critical point of
  $T(\varphi,\beta)$ under variations in the space of coherent angle systems,
  if and only if the isosceles triangles with angles $2\varphi_{\vece}$ and
  $\beta_{\vece}$ as in figure~\ref{fig:decompose_prism}
  {\it (right)}\/ fit together to form a circle pattern with intersection
  angles $\theta$ and cone angles $\Phi$ in the centers of circles.
\end{proposition}

\begin{proof}
The tangent space to the space of of coherent angle systems is spanned by
vectors of the form 
\begin{equation*}
  \frac{\partial}{\partial\beta_{\vece}}
  -\frac{\partial}{\partial\beta_{-\vece}}
\end{equation*}
and
\begin{equation*}
  \frac{\partial}{\partial\varphi_{\vece_1}}
  -\frac{\partial}{\partial\varphi_{\vece_2}},
\end{equation*}
where $\vece_1$ and $\vece_2$ are two consecutive oriented edges in the
oriented boundary of the same face. If and only if
\begin{equation*}
  \Big(\frac{\partial}{\partial\beta_{\vece}}
  -\frac{\partial}{\partial\beta_{-\vece}}\Big)\,T(\varphi,\beta)=0,
\end{equation*}
then the polyhedral pieces corresponding to the edges $\vece$ and $-\vece$\/
fit together to form a polyhedron
$P_{\theta_e}(\varphi_{\vece_1},\varphi_{\vece_2})$ (see
figure~\ref{fig:quadpoly_hyp}). 

If and only if
\begin{equation*}
  \Big(\frac{\partial}{\partial\varphi_{\vece_j}}
  -\frac{\partial}{\partial\varphi_{\vece_k}}\Big)\,T(\varphi,\beta)=0
\end{equation*}
for all pairs of consecutive edges $\vece_j$, $\vece_k$ in the oriented
boundary of a face, then the corresponding polyhedral pieces fit together to
form a (not necessarily triangular) ideal prism as shown in
figure~\ref{fig:ideal_prism} (with a cone-like singular line if
$\Phi_f\not=2\pi$).
\end{proof}

The following two propositions explain how the functional $\Shat(\varphi)$
and Leibon's functional can be obtained as two different reductions of
$T(\varphi,\beta)$.

\begin{proposition}
  If $\varphi$, $\beta$ is a coherent angle system in the sense of this
  section, then $\varphi$ is a hyperbolic coherent angle system in the sense
  of section~\ref{sec:coherent_angle_systems}.
  
  Suppose $\varphi\in\R^{\vecE}$ is a hyperbolic coherent angle system in the
  sense of section~\ref{sec:coherent_angle_systems}. Then there exists
  exactly one $\beta(\varphi)\in\R^{\vecE}$ which
  satisfies~\eqref{eq:phi_beta_inequality} and~\eqref{eq:beta_constraint},
  such that $\varphi$, $\beta(\varphi)$ is a critical point of
  $T(\varphi,\beta)$ under variations with fixed $\varphi$, and with $\beta$
  constrained by equation~\eqref{eq:beta_constraint}. The hyperbolic circle
  pattern functional is 
  \begin{equation*}
    \Shat(\varphi)=T(\varphi,\beta(\varphi)).
  \end{equation*}
\end{proposition}

\begin{proof}[(Sketch of a) proof]
  The first claim about coherent angle systems follows directly from the
  definitions. 
  
  Now, suppose $\varphi\in\R^{\vecE}$ is a hyperbolic coherent angle system
  in the sense of section~\ref{sec:coherent_angle_systems}. For each
  non-oriented edge $e=\{\vece,-\vece\}$, construct the hyperbolic kite with
  angles $2\varphi_{\vece}$, $2\varphi_{-\vece}$, and $\theta_e$, as in
  figure~\ref{fig:kite_f_theta} {\it (left)}. The angle $\theta_e$ is divided
  by one of the diagonals into two angles $\beta_{\vece}$ and
  $\beta_{-\vece}$, such that $\varphi$, $\beta$\/ form a coherent angle
  system in the sense of this section. This $(\varphi,\beta)$\/ is the unique
  critical point of $T(\varphi,\beta)$ under variations with fixed $\varphi$,
  and with $\beta$ constrained by equations~\eqref{eq:beta_constraint}.
  Indeed, it follows from Schl\"afli's formula, that for every such critical
  point, the polyhedra of figure~\ref{fig:decompose_prism} {\it (left)}\/ fit
  together in pairs to form polyhedra as in figure~\ref{fig:quadpoly_hyp}.
  (This means that the isosceles triangles of
  figure~\ref{fig:decompose_prism} {\it (right)} fit together to form kites.)
  Equation $\Shat(\varphi)=T(\varphi,\beta(\varphi))$ follows form the
  additivity of the volume function.
\end{proof}

\begin{proposition}
  Suppose that $\Sigma$ is a triangulation, $\theta$ sums to $2\pi$ around
  each vertex, and $\Phi_f=2\pi$ for all faces $f$.
  
  If $\varphi$, $\beta$ is a coherent angle system in the sense of this
  section, then $\alpha\in\R^{\vecE}$, defined by
  equations~\eqref{eq:beta_123}, is a coherent angle system in the sense of
  section~\ref{sec:leibon_func}.
  
  Suppose $\alpha\in\R^{\vecE}$ is a coherent angle system in the sense of
  section~\ref{sec:leibon_func}, and $\beta\in\R^{\vecE}$ is defined by
  equations~\eqref{eq:beta_123}. Then there exists exactly one
  $\varphi(\beta)\in\R^{\vecE}$ which
  satisfies~\eqref{eq:phi_beta_inequality} and~\eqref{eq:phi_constraint},
  such that $\varphi$, $\beta(\varphi)$ is a critical point of
  $T(\varphi,\beta)$ under variations with fixed $\beta$, and with $\varphi$
  constrained by equation~\eqref{eq:phi_constraint}. Leibon's circle pattern
  functional is
   \begin{equation*}
    H(\beta)= T(\varphi(\beta),\beta).
  \end{equation*}    
\end{proposition}

The proof is similar to the proof of the last proposition.

\chapter{A computer implementation}
\label{cha:computer}

The variational method of constructing circle patterns described in
chapter~\ref{cha:functionals} is well suited for computer implementation.
Here, we present an implementation in Java.  The classes in the
\lstinline|cellularSurface| and \lstinline|circlePattern| package hierarchies
were written by the author. The \lstinline|mfc|,
\lstinline|numericalMethods|, and \lstinline|moebiusViewer| package
hierarchies are part of the jtem project~\cite{jtem}. Of these classes, the
author has written \lstinline|ChebyshevApproximation|, \lstinline|Clausen|, and
\lstinline|HermitianCircle|.
 
The \lstinline|render|
package contains Ken Perlin's 3D renderer~\cite{perlin}. Originally, the
author used Oorange, an experimental programming tool for Java and a part of
the jtem project, to experiment with the circle pattern classes. Here, we use
the BeanShell~\cite{beanshell}, a Java scripting tool, to illustrate how the
classes work.

\section{Getting started}
\label{sec:getting_started}

\subsection*{Unix}

Download the zip-file {\tt circlepatternsoftware.zip} from the URL\\
\url{http://www.math.tu-berlin.de/~springb/software} and unzip it.
Change into the directory {\tt circlepatternsoftware} and execute the shell
script {\tt bin/bsh} with one of the example BeanShell scripts in the {\tt
  examples} directory as parameter:

\begin{verbatim}
> bin/bsh examples/cube.bsh
\end{verbatim}
You should see some text output like this
\begin{verbatim}
Surface has 6 faces, 12 edges, and 8 vertices.
Minimizing ... Minimum found after 237ms. 
Laying out circles ... Done.
\end{verbatim}
and one or more windows displaying circle patterns should open (see
figure~\ref{fig:bsh_examples}).

\subsection*{Windows}

Download the file {\tt circlepatternsoftware.zip} and unzip it as above. Open
the {\tt circlepatternsoftware} folder and, in it, the folders {\tt bin}
and {\tt examples}. Drag and drop an example BeanShell script (such as {\tt
  cube.bsh}) from the {\tt examples} folder onto the {\tt bsh.bat} batch file
in the {\tt bin} folder.

Alternatively, you may start a command shell and invoke the batch file {\tt
  bsh.bat} with one of the example BeanShell scripts as argument.

\subsection*{What the scripts do}

The scripts {\tt bsh} and {\tt bsh.bat} only call the BeanShell interpreter
with the right class path, passing on the argument. Should they not work
properly---for example, because you have a Mac---here is what to do. Make
sure that the directory {\tt cls} and the jar-file(s) in the directory {\tt
  lib} are in the class path. (Of course, no other versions of the
BeanShell, Ken Perlin's renderer, or the jtem classes should be in the
class path.) Invoke a Java virtual machine and run the class {\tt
  bsh.Interpreter} with a BeanShell script as argument.

\subsection*{Running the BeanShell in interactive mode}

If you execute the scripts {\tt bsh} or {\tt bsh.bat} without arguments, the
BeanShell will run in interactive mode. You can run BeanShell scripts using
the BeanShell command {\tt source()}.  Listing~\ref{lst:bsh_session} shows
an example session with the BeanShell in interactive mode. After the
BeanShell is invoked, it prompts the user to enter commands. In the example
session, the script {\tt cube.bsh} is run and then the gradient of the
functional is printed out.

\begin{listing}[p]
\begin{minipage}{\textwidth}
{\tt
> {\bf bin/bsh}\hspace*{\fill}\\
BeanShell 1.2.7 - by Pat Niemeyer (pat@pat.net)\\
bsh \% {\bf source("examples/cube.bsh");}\\
Surface has 6 faces, 12 edges, and 8 vertices.\\
Minimizing ... Minimum found after 232ms. \\
Laying out circles ... Done.\\
bsh \% {\bf for (int i = 0; i $<$ 6; i++) \{\\
System.out.println(data.getGradient(i));\\
\}}\\
4.5469956688748425E-8\\
1.9296370279420216E-8\\
1.9296370279420216E-8\\
1.9296370279420216E-8\\
1.9296370279420216E-8\\
4.5469954468302376E-8\\
bsh \% {\bf exit();}\\
> 
} 
\end{minipage}
\caption{Example session: Running the BeanShell interactively.}
\label{lst:bsh_session}
\end{listing}

\section{The example scripts}
\label{sec:example_scripts}

\begin{figure}[p]%
\begin{tabular}{ccc}
\includegraphics[width=0.3\textwidth,keepaspectratio]{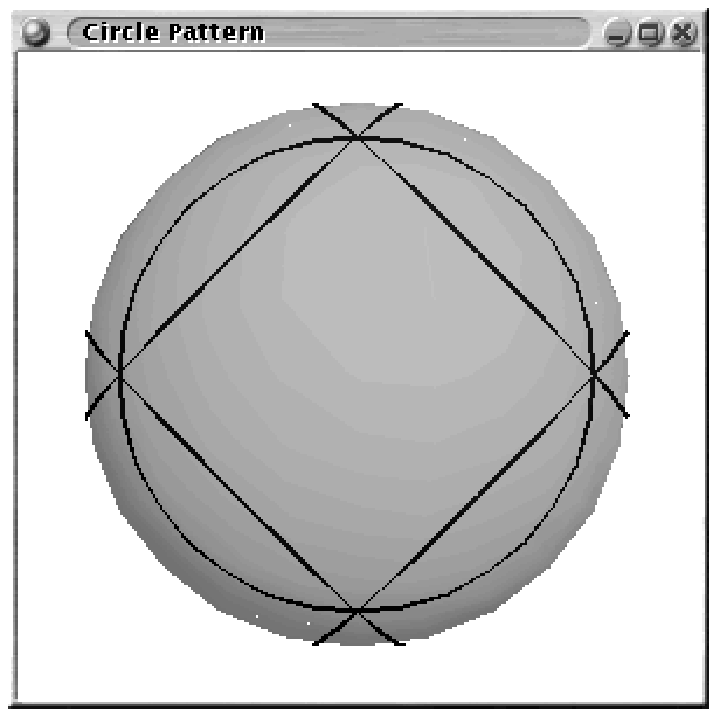}&
\includegraphics[width=0.3\textwidth,keepaspectratio]{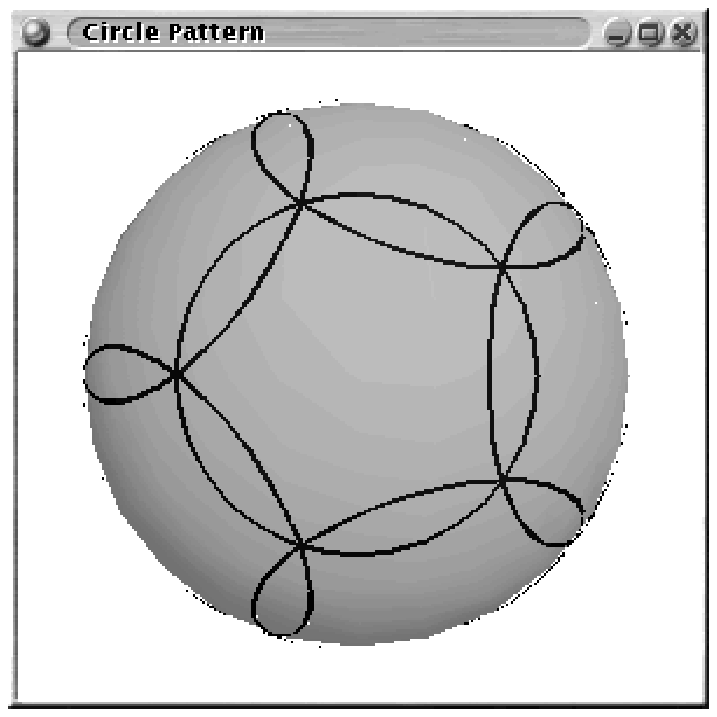}&
\includegraphics[width=0.3\textwidth,keepaspectratio]{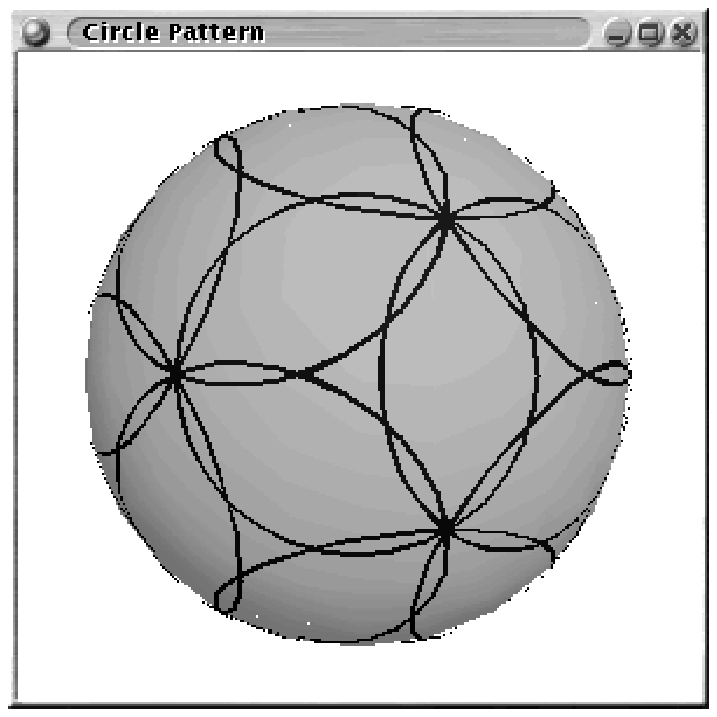}\\
{\small\tt cube.bsh} & 
{\small\tt dodecahedron.bsh} & 
{\small\tt icosahedron.bsh} \\[0.2cm]
\includegraphics[width=0.3\textwidth,keepaspectratio]{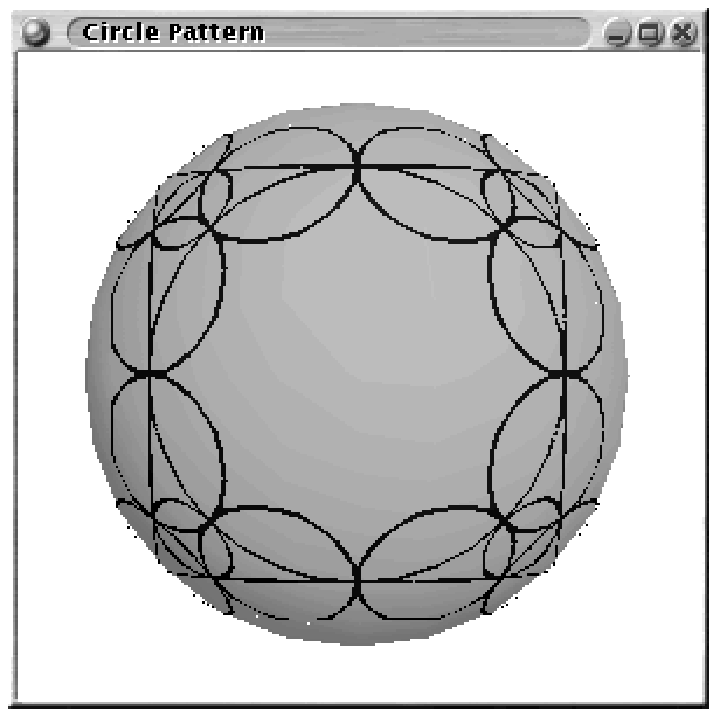}&
\multicolumn{2}{c}{%
\includegraphics[width=0.3\textwidth,keepaspectratio,clip]{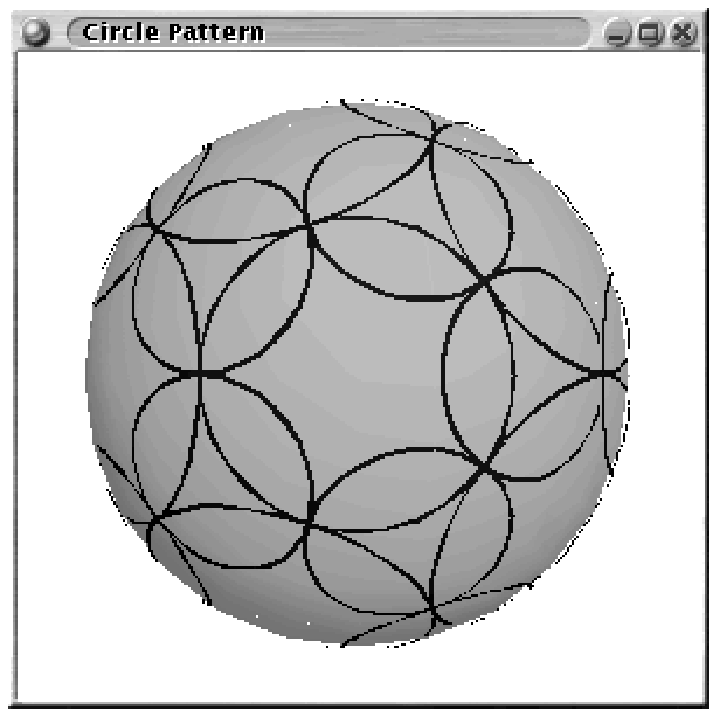}%
\includegraphics[width=0.3\textwidth,keepaspectratio,clip]{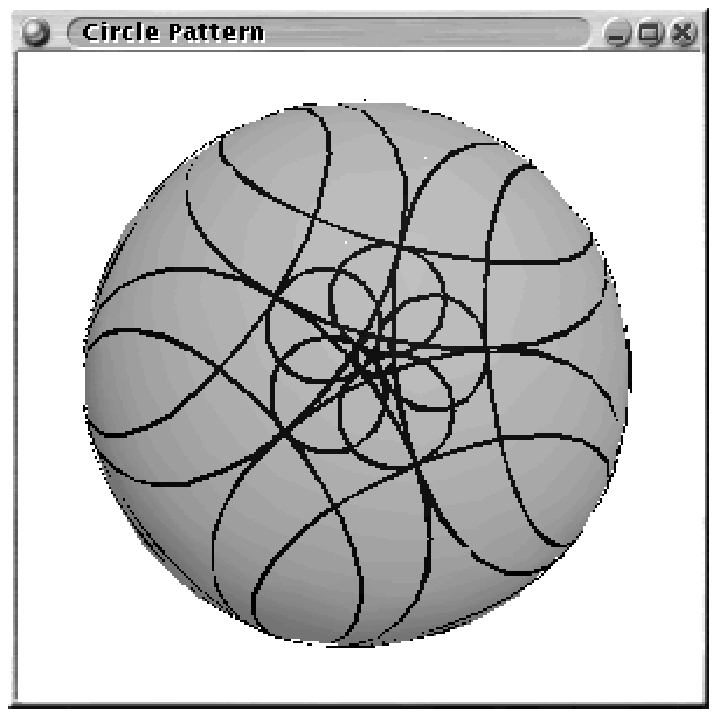}%
}\\
{\small\tt truncated\underline{~}cube.bsh} &
\multicolumn{2}{c}{\small\tt branched.bsh} \\[0.2cm]
\multicolumn{3}{c}{%
\includegraphics[width=0.3\textwidth,keepaspectratio]{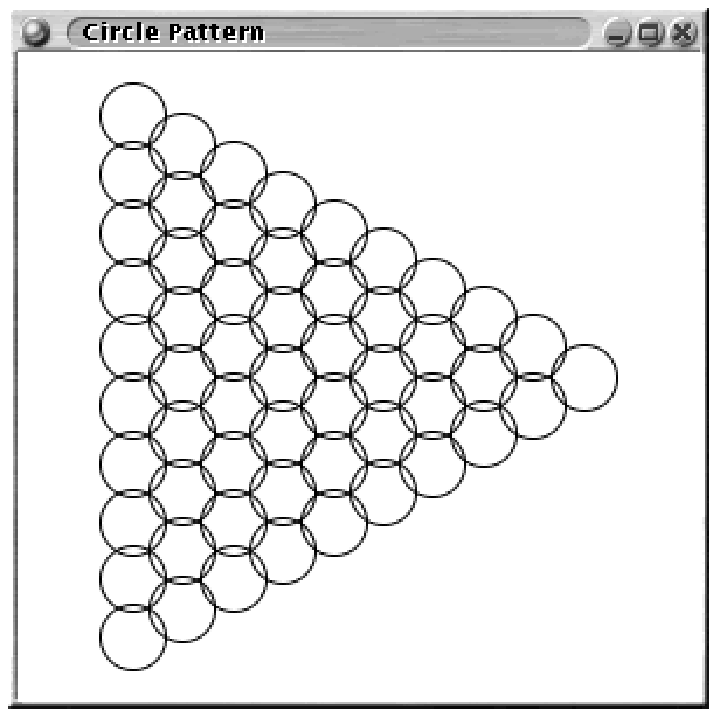}%
\includegraphics[width=0.3\textwidth,keepaspectratio]{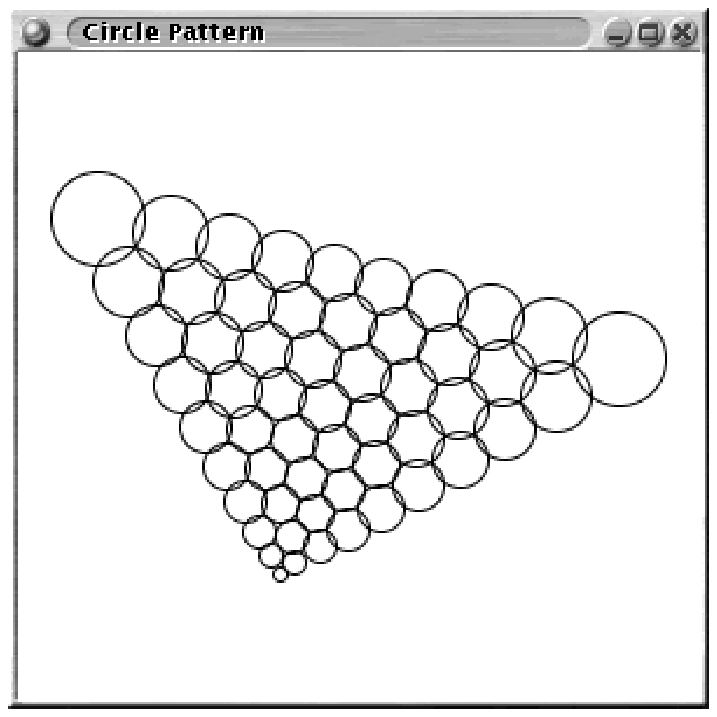}%
\includegraphics[width=0.3\textwidth,keepaspectratio]{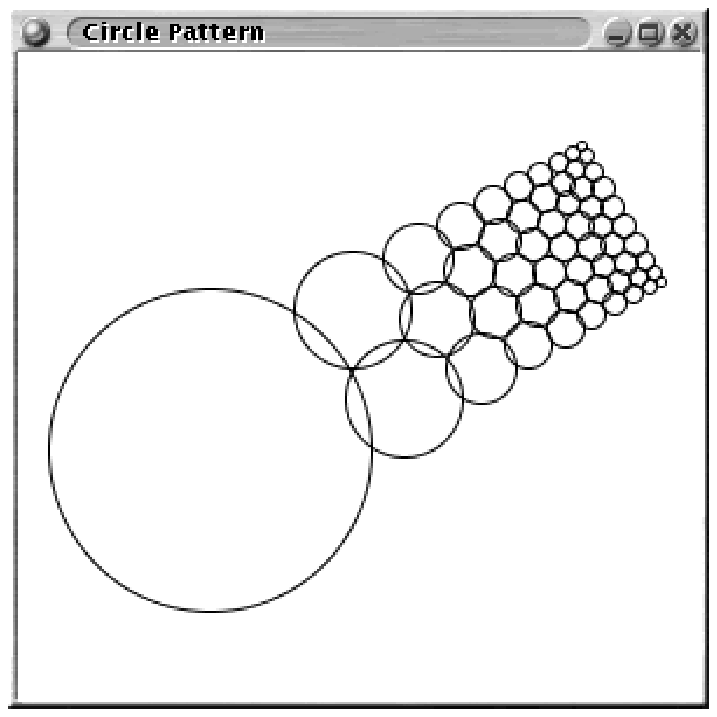}%
}\\
\multicolumn{3}{c}{{\small\tt hex\underline{~}grid.bsh}}\\[0.2cm]
\multicolumn{2}{c}{%
\parbox[b]{0.61\textwidth}{%
\includegraphics[width=0.3\textwidth]{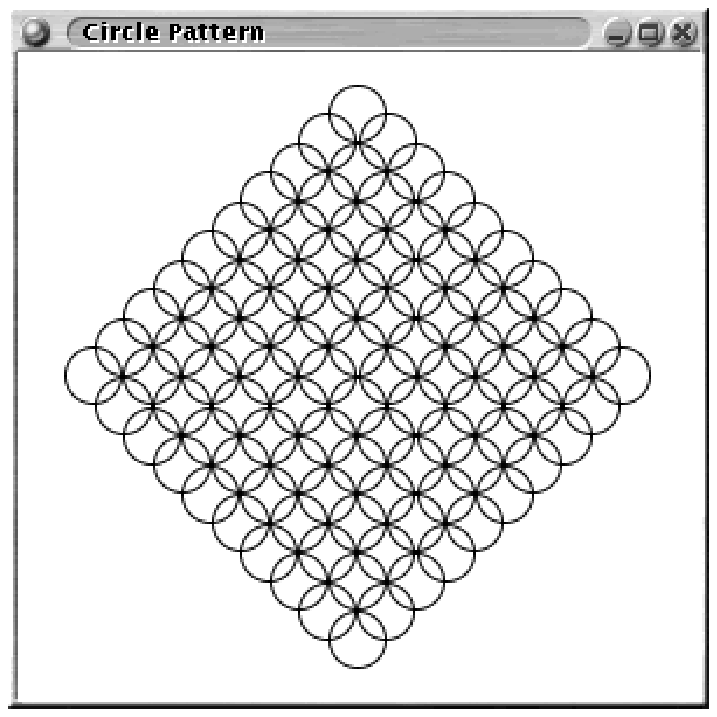}%
\includegraphics[width=0.3\textwidth]{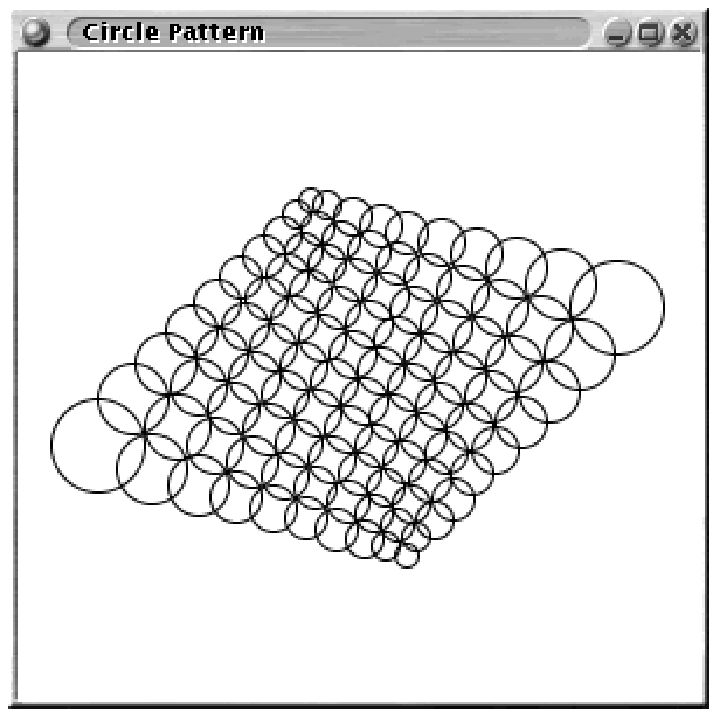}\\
\includegraphics[width=0.3\textwidth]{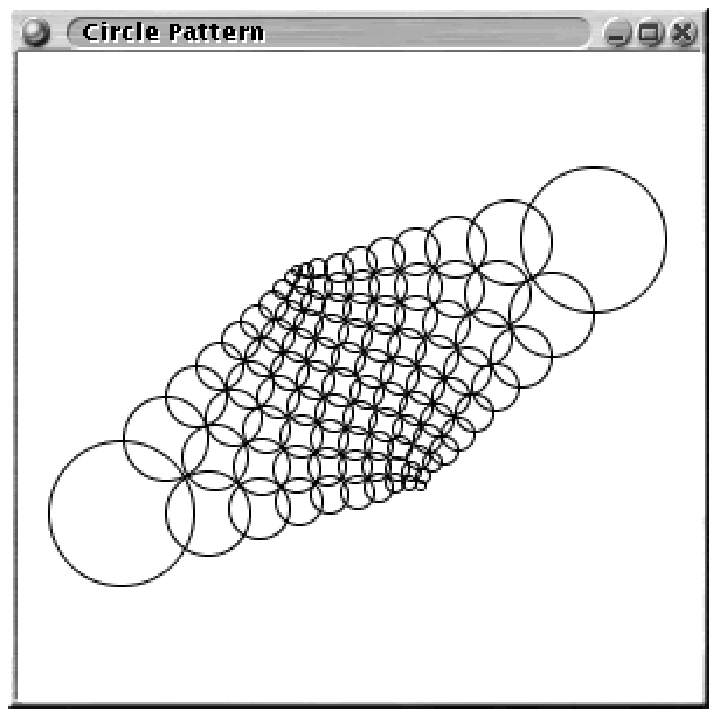}%
\includegraphics[width=0.3\textwidth]{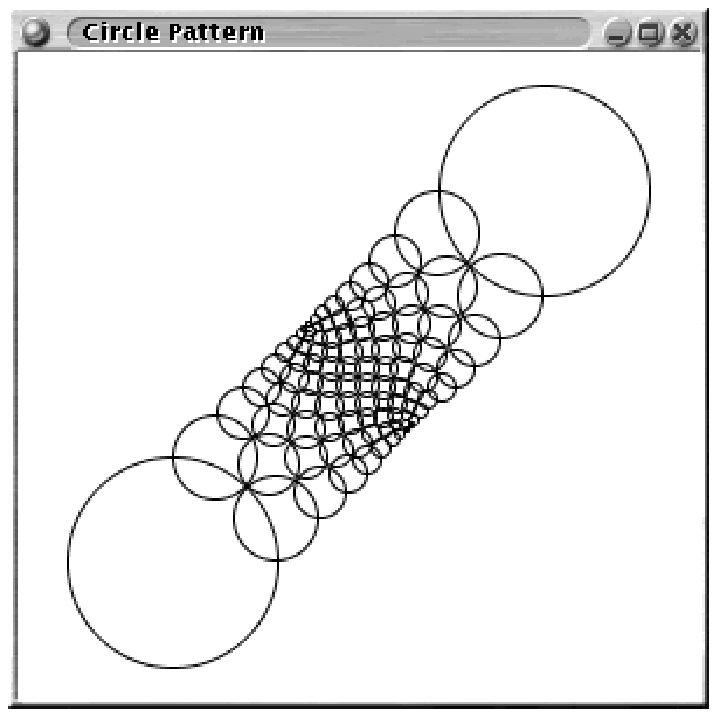}}} &
\includegraphics[width=0.3\textwidth]{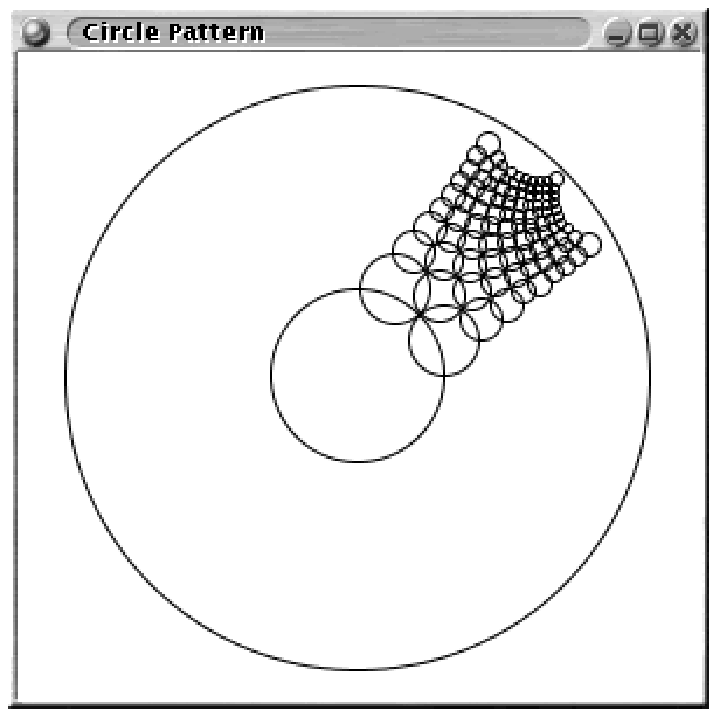}\\
\multicolumn{2}{c}{\small\tt quadmesh.bsh} & 
{\small\tt hyperbolic\underline{~}mesh.bsh}
\end{tabular}
\caption{Circle patterns produced by the example BeanShell scripts.}%
\label{fig:bsh_examples}%
\end{figure}

\begin{listing}[p]
\begin{lstlisting}[indent=12pt]
import cellularSurface.*;
import cellularSurface.examples.*;
import circlePattern.variational.*;
import circlePattern.viewer.perlin.*;
import circlePattern.viewer.moebius.*;

CellularSurface surface = new Cube();
GenericData data = new SphericalData();
GenericLayout layout = new SphericalLayout();

System.out.println("Surface has " + surface.getNumFaces() 
                   + " faces, " + surface.getNumEdges()
                   + " edges, and " + surface.getNumVertices() 
                   + " vertices.");

data.setSurface(surface);
data.assignAllCapitalPhi(2.0*Math.PI);
data.assignAllTheta(2.0/3.0*Math.PI);

//data.setRho(0, 1.0);

System.out.print("Minimizing ... ");
data.adjustRho();

System.out.print("Laying out circles ... ");
layout.layout(data);
System.out.println("Done.");

PerlinView.show(layout);
//MoebiusView.show(layout);

\end{lstlisting}
\caption{\tt cube.bsh}
\label{lst:cube.bsh}
\end{listing}

Figure~\ref{fig:bsh_examples} shows circle patterns produced by BeanShell
scripts in the {\tt examples} directory. The scripts {\tt cube.bsh}, {\tt
  dodecahedron.bsh} and {\tt icosahedron.bsh} produce the obvious trivial
circle patterns in the sphere. A slightly more interesting pattern is
produced by {\tt truncated\underline{~}cube.bsh.} All circles intersect at
right angles. Half the circles form a packing with the combinatorics of a
truncated cube. The script {\tt branched.bsh} produces two circle patterns.
The first is the orthogonal pattern consisting of two packings with
icosahedral and dodecahedral combinatorics. The second is a branched pattern
with the same combinatorics and intersection angles. There are two opposite
branched faces with cone angle $4\pi$.

We will illustrate how the scripts work by taking a closer
look at {\tt cube.bsh}, see listing~\ref{lst:cube.bsh}. First, the relevant
classes are imported.
\begin{lstlisting}[first=1,last=5]
import cellularSurface.*;
import cellularSurface.examples.*;
import circlePattern.variational.*;
import circlePattern.viewer.perlin.*;
import circlePattern.viewer.moebius.*;

CellularSurface surface = new Cube();
GenericData data = new SphericalData();
GenericLayout layout = new SphericalLayout();

System.out.println("Surface has " + surface.getNumFaces() 
                   + " faces, " + surface.getNumEdges()
                   + " edges, and " + surface.getNumVertices() 
                   + " vertices.");

data.setSurface(surface);
data.assignAllCapitalPhi(2.0*Math.PI);
data.assignAllTheta(2.0/3.0*Math.PI);

//data.setRho(0, 1.0);

System.out.print("Minimizing ... ");
data.adjustRho();

System.out.print("Laying out circles ... ");
layout.layout(data);
System.out.println("Done.");

PerlinView.show(layout);
//MoebiusView.show(layout);

\end{lstlisting}
The package \lstinline|cellularSurface| contains only one class,
\lstinline|CellularSurface|, which implements a combinatorial model for cell
decompositions of surfaces. The package \lstinline|cellularSurface.examples|
contains classes like \lstinline|Cube|, which inherit from
\lstinline|CellularSurface| but have different constructors. Thus,
\begin{lstlisting}[first=7,last=7]
import cellularSurface.*;
import cellularSurface.examples.*;
import circlePattern.variational.*;
import circlePattern.viewer.perlin.*;
import circlePattern.viewer.moebius.*;

CellularSurface surface = new Cube();
GenericData data = new SphericalData();
GenericLayout layout = new SphericalLayout();

System.out.println("Surface has " + surface.getNumFaces() 
                   + " faces, " + surface.getNumEdges()
                   + " edges, and " + surface.getNumVertices() 
                   + " vertices.");

data.setSurface(surface);
data.assignAllCapitalPhi(2.0*Math.PI);
data.assignAllTheta(2.0/3.0*Math.PI);

//data.setRho(0, 1.0);

System.out.print("Minimizing ... ");
data.adjustRho();

System.out.print("Laying out circles ... ");
layout.layout(data);
System.out.println("Done.");

PerlinView.show(layout);
//MoebiusView.show(layout);

\end{lstlisting}
assigns to the variable \lstinline|surface|\/ a \lstinline|CellularSurface|
representing the combinatorial type of a cube. 

The package \lstinline|circlePattern.variational| contains the classes
\lstinline|SphericalData| and \lstinline|SphericalLayout|.  The class
\lstinline|SphericalData| holds the intersection angles, cone angles, and
variables $\rho$\/ for a spherical circle pattern. It can compute the correct
$\rho$ by minimizing the spherical circle pattern functional $\Ssph(\rho)$.
Given an instance of \lstinline|SphericalData|, the class
\lstinline|SphericalLayout| calculates the positions of the centers and
intersection points of a spherical circle pattern. (To construct euclidean or
hyperbolic circle patterns, use the analogous classes
\lstinline|EuclideanData|, \lstinline|EuclideanLayout|,
\lstinline|HyperbolicData|, and \lstinline|HyperbolicLayout|.)
\begin{lstlisting}[first=8,last=9]
import cellularSurface.*;
import cellularSurface.examples.*;
import circlePattern.variational.*;
import circlePattern.viewer.perlin.*;
import circlePattern.viewer.moebius.*;

CellularSurface surface = new Cube();
GenericData data = new SphericalData();
GenericLayout layout = new SphericalLayout();

System.out.println("Surface has " + surface.getNumFaces() 
                   + " faces, " + surface.getNumEdges()
                   + " edges, and " + surface.getNumVertices() 
                   + " vertices.");

data.setSurface(surface);
data.assignAllCapitalPhi(2.0*Math.PI);
data.assignAllTheta(2.0/3.0*Math.PI);

//data.setRho(0, 1.0);

System.out.print("Minimizing ... ");
data.adjustRho();

System.out.print("Laying out circles ... ");
layout.layout(data);
System.out.println("Done.");

PerlinView.show(layout);
//MoebiusView.show(layout);

\end{lstlisting}
The next lines print out some information about the \lstinline|surface|.
\begin{lstlisting}[first=11,last=14]
import cellularSurface.*;
import cellularSurface.examples.*;
import circlePattern.variational.*;
import circlePattern.viewer.perlin.*;
import circlePattern.viewer.moebius.*;

CellularSurface surface = new Cube();
GenericData data = new SphericalData();
GenericLayout layout = new SphericalLayout();

System.out.println("Surface has " + surface.getNumFaces() 
                   + " faces, " + surface.getNumEdges()
                   + " edges, and " + surface.getNumVertices() 
                   + " vertices.");

data.setSurface(surface);
data.assignAllCapitalPhi(2.0*Math.PI);
data.assignAllTheta(2.0/3.0*Math.PI);

//data.setRho(0, 1.0);

System.out.print("Minimizing ... ");
data.adjustRho();

System.out.print("Laying out circles ... ");
layout.layout(data);
System.out.println("Done.");

PerlinView.show(layout);
//MoebiusView.show(layout);

\end{lstlisting}
Next, we set up the \lstinline|data|. Since three edges meet at each vertex,
the exterior intersection angles $\theta$ are set to $\frac{2}{3}\pi$.
All $\Phi$ are set to $2\pi$ because there are no boundary faces and the
pattern is unbranched.
\begin{lstlisting}[first=16,last=18]
import cellularSurface.*;
import cellularSurface.examples.*;
import circlePattern.variational.*;
import circlePattern.viewer.perlin.*;
import circlePattern.viewer.moebius.*;

CellularSurface surface = new Cube();
GenericData data = new SphericalData();
GenericLayout layout = new SphericalLayout();

System.out.println("Surface has " + surface.getNumFaces() 
                   + " faces, " + surface.getNumEdges()
                   + " edges, and " + surface.getNumVertices() 
                   + " vertices.");

data.setSurface(surface);
data.assignAllCapitalPhi(2.0*Math.PI);
data.assignAllTheta(2.0/3.0*Math.PI);

//data.setRho(0, 1.0);

System.out.print("Minimizing ... ");
data.adjustRho();

System.out.print("Laying out circles ... ");
layout.layout(data);
System.out.println("Done.");

PerlinView.show(layout);
//MoebiusView.show(layout);

\end{lstlisting}
Then, \lstinline|data| calculates the correct $\rho$, and from this,
\lstinline|layout| calculates the coordinates of the centers and intersection
points.
\begin{lstlisting}[first=22,last=27]
import cellularSurface.*;
import cellularSurface.examples.*;
import circlePattern.variational.*;
import circlePattern.viewer.perlin.*;
import circlePattern.viewer.moebius.*;

CellularSurface surface = new Cube();
GenericData data = new SphericalData();
GenericLayout layout = new SphericalLayout();

System.out.println("Surface has " + surface.getNumFaces() 
                   + " faces, " + surface.getNumEdges()
                   + " edges, and " + surface.getNumVertices() 
                   + " vertices.");

data.setSurface(surface);
data.assignAllCapitalPhi(2.0*Math.PI);
data.assignAllTheta(2.0/3.0*Math.PI);

//data.setRho(0, 1.0);

System.out.print("Minimizing ... ");
data.adjustRho();

System.out.print("Laying out circles ... ");
layout.layout(data);
System.out.println("Done.");

PerlinView.show(layout);
//MoebiusView.show(layout);

\end{lstlisting}
Finally, we display the resulting pattern in a 3-dimensional spherical view.
\begin{lstlisting}[first=29,last=29]
import cellularSurface.*;
import cellularSurface.examples.*;
import circlePattern.variational.*;
import circlePattern.viewer.perlin.*;
import circlePattern.viewer.moebius.*;

CellularSurface surface = new Cube();
GenericData data = new SphericalData();
GenericLayout layout = new SphericalLayout();

System.out.println("Surface has " + surface.getNumFaces() 
                   + " faces, " + surface.getNumEdges()
                   + " edges, and " + surface.getNumVertices() 
                   + " vertices.");

data.setSurface(surface);
data.assignAllCapitalPhi(2.0*Math.PI);
data.assignAllTheta(2.0/3.0*Math.PI);

//data.setRho(0, 1.0);

System.out.print("Minimizing ... ");
data.adjustRho();

System.out.print("Laying out circles ... ");
layout.layout(data);
System.out.println("Done.");

PerlinView.show(layout);
//MoebiusView.show(layout);

\end{lstlisting}
Dragging the mouse rotates the sphere. 

If you uncomment the last line,
\begin{lstlisting}[first=30,last=30]
import cellularSurface.*;
import cellularSurface.examples.*;
import circlePattern.variational.*;
import circlePattern.viewer.perlin.*;
import circlePattern.viewer.moebius.*;

CellularSurface surface = new Cube();
GenericData data = new SphericalData();
GenericLayout layout = new SphericalLayout();

System.out.println("Surface has " + surface.getNumFaces() 
                   + " faces, " + surface.getNumEdges()
                   + " edges, and " + surface.getNumVertices() 
                   + " vertices.");

data.setSurface(surface);
data.assignAllCapitalPhi(2.0*Math.PI);
data.assignAllTheta(2.0/3.0*Math.PI);

//data.setRho(0, 1.0);

System.out.print("Minimizing ... ");
data.adjustRho();

System.out.print("Laying out circles ... ");
layout.layout(data);
System.out.println("Done.");

PerlinView.show(layout);
//MoebiusView.show(layout);

\end{lstlisting}
another window will pop up with a view of the same circle pattern projected
stereographically to the plane. Dragging the mouse will translate the
pattern. Shift click inside the window to select between translate, rotate,
scale, info-coord, and three-point-transform mode. In three-point-transform
mode, you can M\"obius transform the image by dragging three points. Press
the i-key to perform an inversion on the unit circle.

The resulting circle pattern is fairly symmetric---all radii are about
equal---even though the solution is only unique up to M\"obius
transformations. This is so because initially, before the spherical
functional is minimized, all $\rho$ are equal to $0$. If you uncomment line
20,
\begin{lstlisting}[first=20,last=20]
import cellularSurface.*;
import cellularSurface.examples.*;
import circlePattern.variational.*;
import circlePattern.viewer.perlin.*;
import circlePattern.viewer.moebius.*;

CellularSurface surface = new Cube();
GenericData data = new SphericalData();
GenericLayout layout = new SphericalLayout();

System.out.println("Surface has " + surface.getNumFaces() 
                   + " faces, " + surface.getNumEdges()
                   + " edges, and " + surface.getNumVertices() 
                   + " vertices.");

data.setSurface(surface);
data.assignAllCapitalPhi(2.0*Math.PI);
data.assignAllTheta(2.0/3.0*Math.PI);

//data.setRho(0, 1.0);

System.out.print("Minimizing ... ");
data.adjustRho();

System.out.print("Laying out circles ... ");
layout.layout(data);
System.out.println("Done.");

PerlinView.show(layout);
//MoebiusView.show(layout);

\end{lstlisting}
you will get a more unsymmetric solution.

\section{Class overview}
\label{sec:computer_overview}

\newenvironment{topoverview}%
{\begin{list}{}{%
\setlength{\labelwidth}{0pt}%
\setlength{\leftmargin}{0pt}%
\setlength{\itemindent}{0pt}%
\setlength{\listparindent}{0pt}%
\setlength{\itemsep}{\bigskipamount}%
}}%
{\end{list}}

\newenvironment{packageoverview}%
{\begin{list}{}{%
\setlength{\labelwidth}{0pt}%
\setlength{\leftmargin}{30pt}%
\setlength{\itemindent}{-30pt}%
\setlength{\listparindent}{0pt}%
\setlength{\itemsep}{\smallskipamount}%
}}%
{\end{list}}

\newenvironment{classoverview}%
{\begin{list}{}{%
\setlength{\labelwidth}{0pt}%
\setlength{\leftmargin}{0pt}%
\setlength{\itemindent}{-15pt}%
\setlength{\listparindent}{0pt}%
\setlength{\itemsep}{0pt}%
}}%
{\end{list}}

Below, we list the packages and most important classes with some comments.
For more detailed documentation, see sections~\ref{sec:class_cellularsurface}
and \ref{sec:circlepattern_classes}, and the code documentation in the {\tt
  javadocs} directory.

\bigskip\noindent
\hrulefill\ classes pertaining to cellular surfaces\  \hrulefill

\smallskip
\begin{packageoverview}
\item \lstinline|package cellularSurface| : contains only one class.
  \begin{classoverview}
  \item \lstinline|class CellularSurface| : implements a winged
    edge model for cell decompositions of surfaces.
  \end{classoverview}
\item \lstinline|package cellularSurface.examples| : contains classes which
  inherit from the class \lstinline|CellularSurface|.  Their constructors
  create specific cell decompositions.
  \begin{classoverview} 
  \item \lstinline|class Cube extends CellularSurface| : the combinatorial
    type of the cube.
  \item \lstinline|class Dodecahedron extends CellularSurface| : the
    combinatorial type of the dodecahedron.
  \item \lstinline|class ProjectivizedCube extends CellularSurface| : the
    combinatorial type of a cube with diametrically opposite points
    identified. A non-orientable cell decomposition with $3$ faces, $6$
    edges, and $4$ vertices.
  \item \lstinline|class ProjectivizedDodecahedron extends CellularSurface| :
    the combinatorial type of a dodecahedron with diametrically opposite
    points identified. A non-orientable cell decomposition with $6$ faces,
    $15$ edges, and $10$ vertices.
  \item \lstinline|class QuadMesh extends CellularSurface| : a cell
    decomposition of the disc consisting of a rectangular grid of
    quadrilateral faces. (No boundary edges.)
  \item \lstinline|class HexGrid extends CellularSurface| : a cell
    decomposition of the disc consisting of a triangular grid of hexagonal
    faces. (No boundary edges.)
  \end{classoverview}
\end{packageoverview}

\smallskip{}
\noindent
\hrulefill{} classes pertaining to circle patterns \hrulefill{}

\smallskip{}
\begin{packageoverview}
\item \lstinline|package circlePattern.variational| : 
  \begin{classoverview}
    \item \lstinline|class CPMath| : provides static methods for
    mathematical functions that are used in different classes.
  \item \lstinline|class EuclideanFunctional| : provides a static method to
    compute the value of the euclidean functional $\Seuc(\rho)$ and its
    gradient, given a \lstinline&CellularSurface&, intersection angles,
    cone angles, and $\rho$.
    \item \lstinline|class HyperbolicFunctional| : same for the hyperbolic
    functional $\Shyp(\rho)$.
    \item \lstinline|class SphericalFunctional| : same for the spherical
    functional $\Ssph(\rho)$.
    \item \lstinline|abstract class GenericData| : abstract superclass of the
    following three classes. Properties: a \lstinline&CellularSurface&,
    intersection angles, and cone angles.
    \item \lstinline|class EuclideanData extends GenericData| : minimizes the
    euclidean functional.
    \item \lstinline|class HyperbolicData extends GenericData| : minimizes the
    hyperbolic functional.
    \item \lstinline|class SphericalData extends GenericData| : minimizes the
    spherical functional.
    \item \lstinline|abstract class GenericLayout| : abstract superclass of the
      following three classes.
    \item \lstinline|class EuclideanLayout extends GenericLayout| : responsible
      for laying out the circles, given an instance of
      \lstinline|EuclideanData| with correct $\rho$.
    \item \lstinline|class HyperbolicLayout extends GenericLayout| : same for
      \lstinline|HyperbolicData|. 
    \item \lstinline|class SphericalLayout extends GenericLayout| : same for
      \lstinline|SphericalData|.
  \end{classoverview}
\item \lstinline|package circlePattern.viewer.moebius| : helper classes to
  view circle patterns in the \lstinline|MoebiusViewer|. 
  \begin{classoverview}
  \item \lstinline|class CirclePatternShape| : This is an adapter class to
    view circle patterns in the \lstinline|MoebiusViewer|. Holds a reference
    to a \lstinline|GenericLayout| and implements the interface
    \lstinline|MoebiusShape|.
  \item \lstinline|class MoebiusView| : the static method
    \lstinline|show(GenericLayout layout)| opens a frame with a
    \lstinline|MoebiusViewer| displaying the circle pattern
    \lstinline|layout|.
  \end{classoverview}
\item \lstinline|package circlePattern.viewer.perlin| : contains only one
  class.
  \begin{classoverview}
  \item \lstinline|class PerlinView| : The static method
    \lstinline|show(GenericLayout layout)| opens a frame showing a
    3-dimensional spherical view of the circle pattern \lstinline|layout|\/
    in Ken Perlin's renderer.
  \end{classoverview}
\end{packageoverview}

\smallskip{}
\noindent
\hrulefill{} 
\underline{m}athematical \underline{f}oundation \underline{c}lasses
\hrulefill{} 

\smallskip{}
\begin{packageoverview}
\item \lstinline|package mfc.number| : provides the class \lstinline|Complex|
  for complex numbers.
\item \lstinline|package mfc.vector| : provides the classes \lstinline|Real3|
  and \lstinline|Complex2| for real $3$-vectors and complex $2$-vectors.
\item \lstinline|package mfc.matrix| : provides classes for different types
  of complex matrices.
  \begin{classoverview}
  \item \lstinline|abstract class AbstractComplex2By2| : abstract
    superclass of the other matrix classes. Provides functionality, but access
    is restricted. 
  \item \lstinline|class Complex2By2 extends AbstractComplex2By2| : class for
    general complex $2\times2$ matrices.
  \item \lstinline|class HermitianComplex2By2 extends AbstractComplex2By2| :
    class for Hermitian complex $2\times2$ matrices.
  \end{classoverview}
\item \lstinline|package mfc.group| : provides classes for matrix groups.
  \begin{classoverview}
  \item \lstinline|class Moebius extends AbstractComplex2By2| : class for
    elements of the group of orientation preserving conformal self-maps of
    the Riemann sphere.
  \end{classoverview}
\item \lstinline|package mfc.geometry| : provides classes for some geometric
  objects. 
  \begin{classoverview}
  \item \lstinline|class ComplexProjective1 extends Complex2| : homogeneous
    coordinates of 1-dimensional complex projective space.
  \item \lstinline|class HermitianCircle extends HermitianComplex2By2| :
    Oriented circles in the Riemann sphere are described by Hermitian
    $2\times 2$ matrices with negative determinant.
  \end{classoverview}
\end{packageoverview}

\smallskip\noindent
\hrulefill{} classes for numerical calculations \hrulefill

\smallskip
\begin{packageoverview}
\item \lstinline|package numericalMethods.function| : provides interfaces
  which compensate the lack of function pointers in Java.
  \begin{classoverview}
  \item \lstinline|interface DoubleParametrized|
  \item \lstinline|interface DoubleValued|
  \item \lstinline|interface DoubleArrayParametrized| 
  \item \lstinline|interface DoubleArrayValued|
  \end{classoverview}
\item \lstinline|package numericalMethods.calculus.functionApproximation|
  \begin{classoverview}
  \item \lstinline|class ChebyshevApproximation| : provides static methods
    for Chebyshev interpolation and integration.
  \end{classoverview}
\item \lstinline|package numericalMethods.calculus.specialFunctions|
  \begin{classoverview}
  \item \lstinline|class Clausen| : implements Clausen's integral.
  \end{classoverview}
\item \lstinline|package numericalMethods.calculus.minimizing| : provides
  methods for the unconstrained minimization of functions of one or many
  variables.
\end{packageoverview}

\smallskip\noindent
\hrulefill\ the 2D M\"obius viewer\ \hrulefill

\smallskip
\begin{packageoverview}
  \item \lstinline|package moebiusViewer| (and sub-packages).
\end{packageoverview}

\smallskip\noindent
\hrulefill\ Ken Perlin's 3D renderer\ \hrulefill

\smallskip
\begin{packageoverview}
  \item \lstinline|package render| : Some classes have been modified slightly
  by the author to allow resizing of the \lstinline|RenderApplet|. The
  changes are documented in the code (search for ``Springborn'').
\end{packageoverview}

\section{The class \lstinline|CellularSurface|}
\label{sec:class_cellularsurface}

In appendix~\ref{app:combi_top}, {\em non-oriented cellular surfaces}\/ were
defined in terms of a finite set $\vecE$\/ of oriented edges and three
permutations ($\iota$, $\sigma$, and $\tau$) of $\vecE$. The class
\lstinline$CellularSurface$ is an implementation of this combinatorial model.
Here, we describe its fundamental features. See also the documentation in the
{\tt javadocs} directory.

The oriented edges $\vece\in\vecE$\/ are numbered consecutively, starting at
$0$. The edge numbering complies with the following convention: If edge
$\vece$\/ has number $n$, then $\iota\vece$\/ has number $n\xor 2$ and
$\tau\vece$\/ has number $n\xor 1$, where $\xor$ denotes bitwise exclusive or.
Hence, edge-numbers $4k$, $4k+1$, $4k+2$, and $4k+3$ all belong to the same
edge of the cell decomposition. The edges with numbers $4k$ and $4k+2$
($4k+1$ and $4k+3$) are images of each other under $\iota$, that is, they
correspond to the different orientations of the same unoriented edge. The
edges with numbers $4k$ and $4k+1$ ($4k+2$ and $4k+3$) are images of each
other under $\tau$, that is, they correspond to the two preimages in the
oriented double cover of the same oriented edge in the non-oriented
cell-decomposition.

The faces are also numbered consecutively, starting at
$0$, and such that if a face $f$ has number $n$, then the oppositely oriented
face $\tau f$ has number $n\xor 1$. 

The same is true for the vertices. If the vertex $v$ has number $n$, then the
vertex $\tau v$ has number $n\xor 1$.

The default constructor of \lstinline|CellularSurface| produces an `empty'
cell decomposition with $0$ faces, $0$ vertices, and $0$ edges.

\begin{figure}[p]%
\input{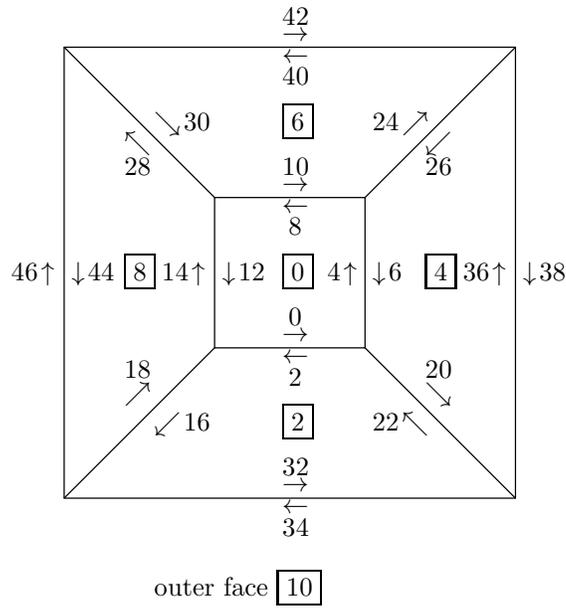}\\%
\vspace{3\baselineskip}%
\input{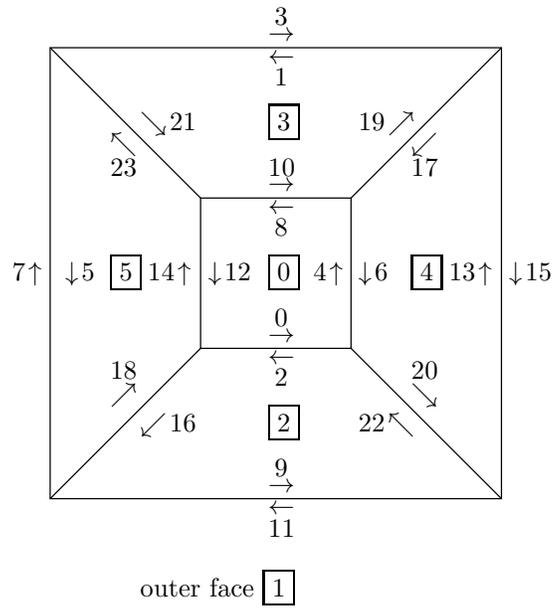}%
\caption{\lstinline|Cube| and \lstinline|ProjectivizedCube|.}%
\label{fig:cellsurf_cube}%
\end{figure}

Figure~\ref{fig:cellsurf_cube} 
shows the edge numbers and face numbers of the
cellular surfaces produced by the constructors of the classes
\lstinline|Cube| and \lstinline|ProjectivizedCube| in the package
\lstinline|cellularSurface.examples|. The first represents the combinatorial
type of the cube. Note that orientable surfaces may be defined using only
even indices. On the other side of the surface, all elements have then odd
indices. The second represents the combinatorial type of a `projectivized
cube', that is, a cube with diametrically opposite points identified. These
two \mbox{\lstinline|CellularSurface|\hspace{2pt}s} are constructed with
calls to the method 
\lstset{labelstep=0, indent=0pt}
\begin{lstlisting}{}
void buildFromFaceBoundaries(int[ ][ ] faceBoundaries) .
\end{lstlisting}
\lstinline|Cube| is produced by a call of \lstinline|buildFromFaceBoundaries|
with the double \lstinline|int| array 
\begin{lstlisting}{}
new int[ ][ ] {
   {0, 4, 8, 12},
   {2, 16, 32, 22},
   {6, 20, 36, 26},
   {10, 24, 40, 30},
   {14, 28, 44, 18},
   {34, 46, 42, 38}
}
\end{lstlisting}
as argument. The argument 
\begin{lstlisting}{}
new int[ ][ ] {
   {0, 4, 8, 12},
   {2, 16, 9, 22},
   {6, 20, 13, 17}
}
\end{lstlisting}
is passed to produce \lstinline|ProjectivizedCube|. Each row defines a face
by listing the oriented edge numbers of the boundary in cyclic order. The
first row defines face number $0$, the second row defines face number $2$,
the third row defines face number $4$ and so forth. (This determines also the
face boundaries of the oppositely oriented odd numbered faces.)  The vertex
indices are generated automatically. 

The methods 
\begin{lstlisting}{}
int getNumFaces()
int getNumVertices()
int getNumEdges()
\end{lstlisting}
return the numbers of non-oriented faces, vertices and edges. The face
indices and vertex indices range from $0$ to 
\begin{lstlisting}{}
2 * getNumFaces() - 1
\end{lstlisting}
and from $0$ to 
\begin{lstlisting}{}
2 * getNumVertices() - 1
\end{lstlisting}
while the edge indices
range from $0$ to 
\begin{lstlisting}{}
4 * getNumEdges() - 1
\end{lstlisting}
because four consecutive indices correspond to a single edge of the cell
decomposition. 

Often, one wants to loop over all non-oriented faces or all non-oriented
edges of a cell decomposition represented by a 
\lstinline|CellularSurface surf|.
This is achieved by loops like
\begin{lstlisting}{}
for (int face = 0, max = 2 * surf.getNumFaces(); face < max; face += 2) {
. . .
}
\end{lstlisting}
and
\begin{lstlisting}{}
for (int edge = 0, max = 4 * surf.getNumEdges(); edge < max; edge += 4) {
. . .
}
\end{lstlisting}
while 
\begin{lstlisting}{}
for (int edge = 0, max = 4 * surf.getNumEdges(); edge < max; edge += 2) {
. . .
}
\end{lstlisting}
loops over all oriented edges of a non-oriented cell decomposition. 

Suppose an oriented edge $\vece$\/ has index \lstinline|i|. When called
with \lstinline|i| as argument,\\[\baselineskip]
\begin{tabular}{ll}
\parbox{0.4\textwidth}{the method} 
& returns the index of (see figure~\ref{fig:cellsurf_edges})\\
\hline
\parbox{0.5\textwidth}{\lstinline|int nextEdgeOfLeftFace(int)|} 
& the edge $\sigma\vece$. \\
\parbox{0.5\textwidth}{\lstinline|int previousEdgeOfLeftFace(int i)|} 
& the edge $\sigma^{-1}\vece$. \\
\parbox{0.5\textwidth}{\lstinline|int nextEdgeOfRightFace(int)|} 
& the edge $\iota\sigma^{-1}\iota\vece$. \\
\parbox{0.5\textwidth}{\lstinline|int previousEdgeOfRightFace(int)|} 
& the edge $\iota\sigma\iota\vece$. \\
\parbox{0.5\textwidth}{\lstinline|int leftEdgeOfInitialVertex(int)|} 
& the edge $\iota\sigma^{-1}\vece$. \\
\parbox{0.5\textwidth}{\lstinline|int rightEdgeOfInitialVertex(int)|} 
& the edge $\sigma\iota\vece$. \\
\parbox{0.5\textwidth}{\lstinline|int leftEdgeOfTerminalVertex(int)|} 
& the edge $\iota\sigma\vece$. \\
\parbox{0.5\textwidth}{\lstinline|int rightEdgeOfTerminalVertex(int)|} 
& the edge $\sigma^{-1}\iota\vece$. \\
\parbox{0.5\textwidth}{\lstinline|int leftFace(int)|} 
& the face on the left side of $\vece${}. \\
\parbox{0.5\textwidth}{\lstinline|int rightFace(int)|} 
& the face on the right side of $\vece${}. \\
\parbox{0.5\textwidth}{\lstinline|int initialVertex(int)|} 
& the initial vertex of $\vece${}. \\
\parbox{0.5\textwidth}{\lstinline|int terminalVertex(int)|} 
& the terminal vertex of $\vece${}.
\end{tabular}\\[\baselineskip]
All these methods are fast, because they consist essentially in an array
lookup. There are also methods like
\begin{lstlisting}{}
int edgeWithLeftFace(int face)
int edgeWithTerminalVertex(int vertex)
\end{lstlisting}
which return an edge in the boundary of \lstinline|face| and an edge ending
in \lstinline|vertex|, respectively. These are computationally more
expensive, because they involve a looping over the edges and checking, say,
whether \lstinline|leftFace(edge)| returns \lstinline|face|. 
The same holds for 
\begin{lstlisting}{}
int boundaryLength(int face).
\end{lstlisting}

To describe cell decompositions of surfaces with boundary, there is a special
index
\begin{lstlisting}{}
static final int NO_ELEMENT = -1
\end{lstlisting}
to indicate that there is no such edge/face/vertex. The class
\lstinline|CellularSurface| may be used to describe both {\em cellular
  surfaces with holes and punctures}\/ and {\em cellular surfaces with
  boundary faces and boundary vertices}, as defined in
section~\ref{sec:surfaces_w_boundary}. For example, after calling
\begin{lstlisting}{}
buildFromFaceBoundaries(new int[ ][ ]{{0,4,8}}) 
\end{lstlisting}
the \lstinline|CellularSurface| describes a cell decomposition of the disc
with one face, three boundary edges, and three vertices. The methods
\lstinline|leftFace(int)| and \lstinline|rightFace(int)| return the following
values:\\
\begin{center}
 \begin{tabular}{c|c|c|c|c|c|c|c|c|c|c|c|c}
\parbox{1em}{\lstinline|i|} 
& 0 & 1 & 2 & 3 & 4 & 5 & 6 & 7 & 8 & 9 & 10 & 11\\ 
\hline
\parbox{5.5em}{\lstinline|leftFace(i)|} 
& 0 & -1 & -1 & 1 & 0 & -1 & -1 & 1 & 0 & -1 & -1 & 1 \\
\parbox{5.5em}{\lstinline|rightFace(i)|} 
& -1 & 1 & 0 & -1 & -1 & 1 & 0 & -1 & -1 & 1 & 0 & -1
\end{tabular} 
\end{center}

\vspace{\baselineskip} The \lstinline|circlePattern|-classes expect
\lstinline|CellularSurface|\hspace{1pt}s which encode cellular surfaces with
boundary faces and boundary vertices. In this case every edge has a left face
and a right face and an initial and a terminal vertex. But some faces may
have non-closed boundary. For example, consider the example \lstinline|Cube|
(see figure~\ref{fig:cellsurf_cube}, {\em top}) with face number 10 and its
boundary edges removed. Such a \lstinline|CellularSurface| is produced by the
method call
\begin{lstlisting}{}
buildFromFaceBoundaries(new int[ ][ ] {
  {0, 4, 8, 12},
  {2, 16, -1, 22},
  {6, 20, -1, 26}, 
  {10, 24, -1, 30},
  {14, 28, -1, 18}
});
\end{lstlisting}
It has $5$ faces, $8$ edges, and $8$ vertices. The method
\lstinline|nextEdgeOfLeftFace(int)| will return $-1$ when called with
argument $16$, $20$, $24$, or $28$ (and also when called with $17$, $21$,
$25$, or $29$).

\begin{figure}[tb]%
\begin{tabular}{@{\extracolsep{0.5cm}}ccc}
\includegraphics{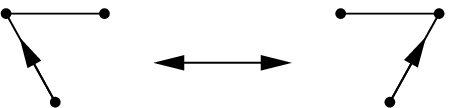} 
& \includegraphics{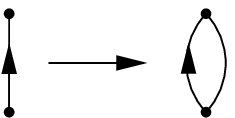}
& \includegraphics{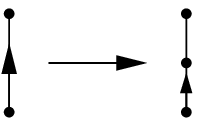} \\
\lstinline|slideEdgeRight| and \lstinline|slideEdgeLeft| 
& \lstinline|splitEdgeAlong|
& \lstinline|splitEdgeAcross|
\end{tabular}%
\caption{Simple moves.}%
\label{fig:cellsurf_moves}%
\end{figure}

The method 
\begin{lstlisting}{}
void copy(CellularSurface surface)
\end{lstlisting}
builds a copy of \lstinline|surface| and the method
\begin{lstlisting}{}
void buildPoincareDual(CellularSurface surface)
\end{lstlisting}
builds the Poincar\'e dual of \lstinline|surface|. For example, a
\lstinline|CellularSurface| representing the combinatorial type of an
icosahedron may be obtained in this way:
\begin{lstlisting}{}
import cellularSurface.CellularSurface;
import cellularSurface.examples.Dodecahedron;
CellularSurface surface = new Dodecahedron();
surface.buildPoincareDual(surface);
\end{lstlisting}

There are several methods to modify a \lstinline|CellularSurface|.
Figure~\ref{fig:cellsurf_moves} illustrates the simple moves.
The methods
\begin{lstlisting}{}
void slideEdgeLeft(int edge)
void slideEdgeRight(int edge)
\end{lstlisting}
do not change the numbers of faces, edges and vertices, while 
\begin{lstlisting}{}
int splitEdgeAlong(int edge)
\end{lstlisting}
introduces a new edge and a new face, and
\begin{lstlisting}{}
int splitEdgeAcross(int edge)
\end{lstlisting}
introduces a new edge and a new vertex. These last two methods return the
index of the new edge. An edge may be contracted with 
\begin{lstlisting}{}
void contractEdge(int edge).
\end{lstlisting}
More complicated moves such as truncating a vertex with
\begin{lstlisting}{}
void truncateVertex(int vertex)
\end{lstlisting}
and building the medial decomposition (see figure~\ref{fig:medial}) with 
\begin{lstlisting}{}
void medialGraph()
\end{lstlisting}
are composed of many of the simple moves.

\section{The \lstinline|circlePattern|-classes}
\label{sec:circlepattern_classes}

Here, we give a brief description of the classes which perform the
construction of circle patterns. They are contained in the package
\lstinline|circlePattern.variational|. The construction proceeds in two
steps. First, one of the classes \lstinline|EuclideanData|,
\lstinline|HyperbolicData|, or \lstinline|SphericalData| computes the radii
of the circles, and then then corresponding class
\lstinline|EuclideanLayout|, \lstinline|HyperbolicLayout|, or
\lstinline|SphericalLayout| computes the centers and intersection points of
the circles in the respective constant curvature space.

\subsection*{The {---\hspace{-1pt}---}\mbox{\lstinline|Data|} classes}

The classes \lstinline|EuclideanData|, \lstinline|HyperbolicData|, and
\lstinline|SphericalData| perform the minimization of the respective
functionals to compute the radii. They all inherit from the abstract
superclass \lstinline|GenericData| which provides the member variables with
\lstinline|get|/\lstinline|set| methods. 

\smallskip
The \lstinline|protected| member variables, which can be accessed
with get/set methods, are:

\smallskip
\paragraph*{$\bullet$ \lstinline|CellularSurface surface|} The cell
decomposition which determines the combinatorics of the circle pattern. Faces
of the decomposition correspond to circles of the pattern. The
\lstinline|surface|\/ may be non-orientable. It may be closed or a cellular
surface with boundary faces and boundary vertices (see
sections~\ref{sec:class_cellularsurface} and~\ref{sec:surfaces_w_boundary}).
Boundary edges are not allowed.

\smallskip
\paragraph*{$\bullet$ \lstinline|double[ ] capitalPhi|} The cone/boundary
angles $\Phi$. Before \lstinline|surface|\/ is set, \lstinline|capitalPhi|\/
may be \lstinline|null|. When \lstinline|surface|\/ is set,
\lstinline|capitalPhi|\/ is assigned and array with length
\lstinline|surface.getNumFaces()|. There is one angle $\Phi_f$ for each
non-oriented face $f$. The oriented faces with indices $2n$ and $2n+1$ have
the cone/boundary angle \lstinline|capitalPhi[n]|.

\smallskip
\paragraph*{$\bullet$ \lstinline|double[ ] theta|} The exterior intersection
angles $\theta$ (see figure~\ref{fig:intersection_angle}). Before
\lstinline|surface|\/ is set, \lstinline|theta|\/ may be \lstinline|null|.
When \lstinline|surface|\/ is set, \lstinline|theta|\/ is assigned an array
with length \lstinline|surface.getNumEdges()|. There is one intersection
angle $\theta_e$ for each non-oriented edge $e$. The oriented edges with
indices $4n$, $4n+1$, $4n+2$ and $4n+3$ have the intersection angle
\lstinline|theta[n]|.

\smallskip
\paragraph*{$\bullet$ \lstinline|double tolerance|} The error tolerance for the
minimal value of the functional. The default is
\lstinline|DEFAULT_TOLERANCE|.

\smallskip
\paragraph*{$\bullet$ \lstinline|double[ ] rho|} The variables $\rho$. Before
\lstinline|surface|\/ is set, \lstinline|rho|\/ may be \lstinline|null|. When
\lstinline|surface|\/ is set, $\rho$ is assigned an array with length
\lstinline|surface.getNumFaces()|. There is one variable $\rho_f$ for each
non-oriented face $f$. To the oriented faces with indices $2n$ and $2n+1$
belongs the variable \lstinline|rho[n]|. Apart from the set methods,
the values
\lstinline|rho[n]|
are changed by the method \lstinline|adjustRho()|.

\smallskip
\paragraph*{$\bullet$ \lstinline|double value|} The value of the respective
circle pattern functional. This is set by \lstinline|evaluateFunctional()|
and \lstinline|adjustRho()|. There is no set method.

\smallskip
\paragraph*{$\bullet$ \lstinline|double[ ] gradient|} The gradient of the
respective circle pattern functional. Before \lstinline|surface|\/ is set,
\lstinline|rho|\/ may be \lstinline|null|. When \lstinline|surface|\/ is set,
\lstinline|rho|\/ is assigned an array with length
\lstinline|surface.getNumFaces()|. To the oriented faces with indices $2n$
and $2n+1$ belongs the gradient entry \lstinline|gradient[n]|. The values
\lstinline|gradient[n]| are set by the methods
\lstinline|evaluateFunctional()| and \lstinline|adjustRho()|. There are no
set methods.

\smallskip{} The functionality is provided by the following methods, which
are abstract in \lstinline|GenericData|:

\smallskip
\paragraph*{$\bullet$ \lstinline|void evaluateFunctional()|} From the values of
\lstinline|surface|, \lstinline|capitalPhi|, \lstinline|theta|, and
\lstinline|rho|, this method calculates the value and gradient of the
functional $\Seuc$, $\Shyp$, or $\Ssph$ (depending on the subclass of
\lstinline|GenericData|) and sets \lstinline|value| and \lstinline|gradient|
accordingly.

\smallskip
\paragraph*{$\bullet$ \lstinline|void adjustRho()|} Minimizes the  functional
$\Seuc$, $\Shyp$, or $\Ssph$ (depending on the subclass of
\lstinline|GenericData|) using a conjugate gradient method. (In the case of
the spherical functional $\Ssph$, the minimax procedure described in
section~\ref{sec:sph_func} is performed.) It attempts to determine the
minimal value up to a maximal error of \lstinline|tolerance|. If it returns
successfully, the member variables \lstinline|rho|, \lstinline|value|, and
\lstinline|gradient| hold the approximated position of the minimum, the
approximated minimal value and the gradient at the approximated minimum
(which should be close to zero).

\smallskip
\paragraph*{$\bullet$ %
\lstinline[deleteemph={[4]radius}]|double radius(double rho)|} Performs the
variable transformation from the variables $\rho$ to the radii $r$. The
relation between $r$ and $\rho$ is given by
equations~\eqref{eq:rho_of_r_euc}, \eqref{eq:rho_hyp}, or \eqref{eq:rho_sph},
depending on the subclass of \lstinline|GenericData|.

\subsection*{The {---\hspace{-1pt}---}\mbox{\lstinline|Layout|} classes} Once
the radii of the circles have been computed, the classes
\lstinline|EuclideanLayout|, \lstinline|HyperbolicLayout|, and
\lstinline|SphericalLayout| compute the centers and intersection points of
the circles. They all inherit from the abstract superclass
\lstinline|GenericLayout|, which provides the member variables and almost all
of the functionality.

\smallskip 
The \lstinline|protected| member variables, which can only be read 
with get methods, are:

\smallskip
\paragraph*{$\bullet$ \lstinline|ComplexProjective1[ ] centerPoint|}
Complex homogeneous coordinates for the centers of the circles. After a call
of the method \lstinline[deleteemph={[4]layout}]|layout(data)| (see below), 
\lstinline|centerPoint| is an array of length 
\lstinline|2*data.getNumFaces()|. There is one entry for each oriented
face. If the oriented double cover of the surface is not connected (for
example, if the surface is orientable),  then
some entries of \lstinline|centerPoint|\/ will be null. 

\smallskip
\paragraph*{$\bullet$ \lstinline|ComplexProjective1[ ] vertexPoint|} Complex
homogeneous coordinates for the intersection points. After a call of the
method \lstinline[deleteemph={[4]layout}]|layout(data)|,
\lstinline|vertexPoint| is an array of length
\lstinline|2*data.getNumVertices()|. There is one entry for each vertex of the
oriented double cover. If the oriented double cover is not connected (for
example, if the surface is orientable), then some entries of 
\lstinline|vertexPoint|\/ will be null.

\smallskip
\paragraph*{$\bullet$ \lstinline|HermitianCircle[ ] circle|} The circles of
the circle pattern, described in a M\"obius-invariant fashion as Hermitian
$(2\times2)$-matrices with negative determinant. After a call
of the method \lstinline[deleteemph={[4]layout}]|layout(data)|, 
\lstinline|circle| is an array of length 
\lstinline|2*data.getNumFaces()|. There is one entry for each oriented
face. If the oriented double cover of the surface is not connected (for
example, if the surface is orientable),  then
some entries of \lstinline|circle|\/ will be null. 

\smallskip
The member variables are set by the method 
\begin{lstlisting}[deleteemph={[4]layout}]{}
void layout(GenericData data)
\end{lstlisting}
of class \lstinline|GenericLayout|. The member variables \lstinline|surface|,
\lstinline|theta|, and \lstinline|rho| of the argument \lstinline|data| are
expected to hold correct values for a circle pattern of the respective
geometry. This means, the \lstinline[deleteemph={[4]layout}]|layout| method
of \lstinline|EuclideanLayout| expects to be passed an instance of
\lstinline|EuclideanData| with \lstinline|rho| set by a successful call of
its method \lstinline|adjustRho()|, etc.

The layout algorithm is fairly unsophisticated: Initialize member variables
\lstinline|centerPoint|, \lstinline|vertexPoint|, and \lstinline|circle| to
arrays of the correct size with all entries equal to \lstinline|null|. Then,
place the first two points: set the \lstinline|centerPoint|\/ of the left
face of the edge with index \lstinline|INITIAL_EDGE| to the origin;
and set the \lstinline|vertexPoint|\/ of the initial vertex of the same edge
an appropriate distance away on the positive real axis. After the initial
points are placed, keep searching for more points whose position can be
determined until no more can be found.

Almost all the work is done in the abstract superclass
\lstinline|GenericLayout|. Its \lstinline[deleteemph={[4]layout}]|layout|\/
method is not abstract. But internally, it calls the two abstract methods of
\lstinline|GenericLayout|:
\begin{lstlisting}{}
abstract void assignOtherFixPoint(ComplexProjective1 center, 
                                    ComplexProjective1 otherFixPoint);
abstract void assignCircle(HermitianCircle circle, 
                            ComplexProjective1 center,
                            double radius);
\end{lstlisting}
All the subclasses \lstinline|EuclideanLayout|, \lstinline|HyperbolicLayout|,
and \lstinline|SphericalLayout| do is to implement these two abstract
methods:

The method \lstinline|assignOtherFixPoint| assigns the argument
\lstinline|otherFixPoint|\/ depending on the argument \lstinline|center|. In
the class \lstinline|EuclideanFunctional|, the argument
\lstinline|otherFixPoint| is always set to the infinite point of the
euclidean plane. In the class \lstinline|HyperbolicFunctional|,
\lstinline|otherFixPoint| is set to image of \lstinline|center| under an
inversion on the circle which is the boundary of hyperbolic space. In the
class \lstinline|SphericalFunctional|, \lstinline|otherFixPoint| is set to
the point which is diametrically opposite to  \lstinline|center| on the
Riemann sphere.

The method \lstinline|assignCircle| assigns \lstinline|circle| by the circle
which, in the respective geometry, has center \lstinline|center| and radius
\lstinline|radius|.

\section{Computing Clausen's integral}
\label{sec:computing_clausen}

In 1832, Clausen first tabulated the integral \eqref{eq:clausen_int}, which
now bears his name \cite{clausen32}. He expands the functions
\begin{equation*}
  \Cl(\pi x) + \pi\Big(2\log\frac{2+x}{2-x}+x\log\big((2+x)(2-x)|x|\big)\Big)
\end{equation*}
and
\begin{equation*}
  \Cl(\pi(1-x)) - \pi\log\frac{1+x}{1-x}-\pi x\log\big((1+x)(1-x)\big)
\end{equation*}
into power series of $x$, which he uses to calculate $\Cl(x)$ in the
intervals $[0,\frac{\pi}{2}]$ and $[\frac{\pi}{2},\pi]$ to 16 decimals.

We proceed in a similar fashion, except that we use Chebyshev series. For
$x\in[-\pi, \pi]$, the function
\begin{equation*}
  h(x) = -\log\bigg(\frac{2\sin\frac{x}{2}}{x}\bigg)
\end{equation*}
is analytic and
\begin{equation*}
  \Cl(x) = \int_0^x h(\xi)\,d\xi-x (\log|x|-1).
\end{equation*}
The function $h(x)$ is approximated in the interval $[-\pi$, $\pi]$ by a
Chebyshev series, which can easily be integrated. 

The algorithms for fitting, integrating and evaluating Chebyshev series are
described in the {\it Numerical Recipes} \cite{numerical_recipes}. They are
implemented in static methods of the class
\lstinline|ChebyshevApproximation|, which is contained in the package
\lstinline|numericalMethods.calculus.functionApproximation|.

Using these methods for Chebyshev approximation, the class \lstinline|Clausen| in package 
\lstinline|numericalMethods.calculus.specialFunctions| 
provides the static method 
\begin{lstlisting}{}
double cl2(double x)
\end{lstlisting}
to evaluate Clausen's integral. It also provides Catalan's constant $G=\Cl(\frac{\pi}{2})$:
\begin{lstlisting}{}
static final double CATALAN
\end{lstlisting}

\appendix

\chapter{Proof of the trigonometric relations of
  lemma~$\text{\ref{lem:phi_hyp}}$ and lemma~$\text{\ref{lem:phi_sph}}$}
\label{app:trig}

In this appendix, we derive the below formulae for the remaining angles in a
hyperbolic or spherical triangle when two sides and the enclosed angle
between them are given.  Suppose the given sides of the triangle are $r_1$
and $r_2$, the included angle between them is $\theta$, its complement is
$\theta^* = \pi - \theta$, and the remaining angles are $\varphi_1$ and
$\varphi_2$, as shown in figure~\ref{fig:phiRhoTriangle}.

In the spherical case,
\begin{equation}\label{eq:phiOfRhoSph}
  \varphi_1 = \frac{1}{2i}\log
  \frac{1+e^{\rho_2-\rho_1+i\theta^*}}{1+e^{\rho_2-\rho_1-i\theta^*}}
  - \frac{1}{2i}\log
  \frac{1-e^{\rho_1+\rho_2+i\theta^*}}{1-e^{\rho_1+\rho_2-i\theta^*}}\;,
\end{equation}
where we assume $0<r_j<\pi$, and
\begin{equation*}
  \rho_j=\log\tan\frac{r_j}{2}.
\end{equation*}
In the hyperbolic case,
\begin{equation}\label{eq:phiOfRhoHyp}
  \varphi_1 = \frac{1}{2i}\log
  \frac{1+e^{\rho_2-\rho_1+i\theta^*}}{1+e^{\rho_2-\rho_1-i\theta^*}}
  - \frac{1}{2i}\log
  \frac{1+e^{\rho_1+\rho_2+i\theta^*}}{1+e^{\rho_1+\rho_2-i\theta^*}}\;,
\end{equation}
where
\begin{equation*}
  \rho_j=\log\tanh\frac{r_j}{2}.
\end{equation*}
These equations~\eqref{eq:phiOfRhoSph} and~\eqref{eq:phiOfRhoHyp} imply
equations~\eqref{eq:phi_of_rho_sph} and~\eqref{eq:phi_of_rho_hyp}.

To derive such formulae from known relations of non-euclidean trigonometry
can be cumbersome. Instead, we present a method using M\"obius transformations
which has proved to be helpful. We identify the group of M\"obius
transformations of $\Chat=\C\cup\{\infty\}$ with the projectivized linear
group $\PGL(2,\C)$. A M\"obius transformation
\begin{equation*}
  z \mapsto \frac{a z + b}{c z + d}\;,\quad ad -bc \neq 0
\end{equation*}
corresponds to
\begin{equation*}
  \begin{bmatrix}a&b\\c&d\end{bmatrix}
  :=\C\cdot\begin{pmatrix}a&b\\c&d\end{pmatrix}\in\PGL(2,\C).
\end{equation*}

\section{The spherical case}

The unit sphere $S^2=\big\{(x_1, x_2, x_3)\in\R^3 \;\big|\; {x_1}^2 + {x_2}^2
+ {x_3}^2 = 1\big\}$ is mapped to the extended complex plane
$\Chat=\C\cup\{\infty\}$ by stereographic projection
\[ (x_1, x_2, x_3) \mapsto z = \frac{x_1 + ix_2}{1 - x_3}\,. \]
The orientation preserving isometries of the sphere (the rotations)
correspond to M\"obius transformations
$\big[\begin{smallmatrix}a&b\\c&d\end{smallmatrix}\big]$ with $a=\bar{d}$ and
$b=-\bar{c}$.  Rotations around the $x_1$-, $x_2$-, and $x_3$-axis by an
angle $\alpha$ correspond to
\begin{equation*}
  \begin{split}
    R_1(\alpha) &=
    \begin{bmatrix} 
      \cos\frac{\alpha}{2} & i\sin\frac{\alpha}{2} \\
      i\sin\frac{\alpha}{2} & \cos\frac{\alpha}{2}
    \end{bmatrix}, \\
    R_2(\alpha) &=
    \begin{bmatrix} 
      \cos\frac{\alpha}{2} & -\sin\frac{\alpha}{2} \\
      \sin\frac{\alpha}{2} & \cos\frac{\alpha}{2}
    \end{bmatrix}, \\
    R_3(\alpha) &=
    \begin{bmatrix} 
      e^{i\alpha/2} & 0 \\
      0 & e^{-i\alpha/2}
    \end{bmatrix}
    =
    \begin{bmatrix} 
      e^{i\alpha} & 0 \\
      0 & 1
    \end{bmatrix},
  \end{split}
\end{equation*}
respectively. If $\alpha$ is not an odd multiple of $\pi$, then
\begin{equation*}
  \begin{split}
    &R_1(\alpha) =
    \begin{bmatrix} 
      1 & i\tan\frac{\alpha}{2} \\
      i\tan\frac{\alpha}{2} & 1
    \end{bmatrix}, \\
    &R_2(\alpha) =
    \begin{bmatrix} 
      1 & -\tan\frac{\alpha}{2} \\
      \tan\frac{\alpha}{2} & 1
    \end{bmatrix}.
  \end{split}
\end{equation*}
The quantities $r_1$, $r_2$, $l$, and $\varphi_1$, $\varphi_2$,
$\theta=\pi-\theta^*$ are the sides and angles of a spherical triangle as in
figure~\ref{fig:phiRhoTriangle} if and only if
\begin{equation*}
  R_2(r_1)R_3(\theta^*)R_2(r_2) = R_3(\varphi_1)R_2(l)R_3(\varphi_2).
\end{equation*}
With $\rho_j=\log\tan(r_j/2)$, the left hand side of this equation equals
\begin{equation*}
  \begin{split}
    \LHS &=
    \begin{bmatrix}1&-e^{\rho_1}\\e^{\rho_1}&1\end{bmatrix}
    \begin{bmatrix}e^{i\theta^*}&0\\0&1\end{bmatrix}
    \begin{bmatrix}1&-e^{\rho_2}\\e^{\rho_2}&1\end{bmatrix} \\
    &=
    \begin{bmatrix}
      e^{i\theta^*}-e^{\rho_1+\rho_2} & -e^{\rho_2+i\theta^*}-e^{\rho_1} \\
      e^{\rho_1+i\theta^*}+e^{\rho_2} & 1-e^{\rho_1+\rho_2+i\theta^*}
    \end{bmatrix} \\
    &=
    \begin{bmatrix}
      \frac{\displaystyle e^{i\theta^*}-e^{\rho_1+\rho_2}} {\displaystyle
        1-e^{\rho_1+\rho_2+i\theta^*}} & -\frac{\displaystyle
        e^{\rho_2+i\theta^*}+e^{\rho_1}}
      {\displaystyle 1-e^{\rho_1+\rho_2+i\theta^*}}    \\
      \frac{\displaystyle e^{\rho_1+i\theta^*}+e^{\rho_2}} {\displaystyle
        1-e^{\rho_1+\rho_2+i\theta^*}} & 1
    \end{bmatrix}.
  \end{split}
\end{equation*}
With $\lambda=\log\tan(l/2)$, the right hand side equals
\begin{equation*}
  \begin{split}
    \RHS &=
    \begin{bmatrix}e^{i\varphi_1}&0\\0&1\end{bmatrix}
    \begin{bmatrix}1&-e^{\lambda}\\e^{\lambda}&1\end{bmatrix}
    \begin{bmatrix}e^{i\varphi_2}&0\\0&1\end{bmatrix} \\
    &=
    \begin{bmatrix}
      e^{i(\varphi_1+\varphi_2)} & -e^{\lambda+i\varphi_1} \\
      e^{\lambda+i\varphi_2} & 1
    \end{bmatrix}.
  \end{split}
\end{equation*}
Hence,
\begin{equation*}
  \varphi_1 = \arg\frac{\displaystyle e^{\rho_2+i\theta^*}+e^{\rho_1}}
                       {\displaystyle 1-e^{\rho_1+\rho_2+i\theta^*}},
\end{equation*}
which implies equation \eqref{eq:phiOfRhoSph}.

\section{The hyperbolic case}

We use the Poincar\'e disc model of hyperbolic 2-space: $H^2=\big\{z\in\C
\;\big|\; |z| < 1\big\}$ with metric $|ds| = 2|dz|/(1-|z|^2)$. The
orientation preserving isometries of $H^2$ are the M\"obius transformations
$\big[\begin{smallmatrix}a&b\\c&d\end{smallmatrix}\big]$ with $a=\bar{d}$ and
$b=\bar{c}$. Let
\begin{equation*}
  \begin{split}
    T(s) &= \begin{bmatrix}\cosh\frac{s}{2}&\sinh\frac{s}{2}\\
      \sinh\frac{s}{2}&\cosh\frac{s}{2}
            \end{bmatrix}
            =
    \begin{bmatrix}1&\tanh\frac{s}{2}\\\tanh\frac{s}{2}&1\end{bmatrix}, \\
    R(\alpha) &=
    \begin{bmatrix}e^{i\alpha/2}&0\\0&e^{-i\alpha/2}\end{bmatrix}
    =\begin{bmatrix}e^{i\alpha}&0\\0&1\end{bmatrix}.
  \end{split}
\end{equation*}
The transformation $T(s)$ translates all points on the line $\big\{z\in H^2
\big| \im z = 0 \big\}$ towards $z=1$ by the distance $s$. The transformation
$R(\alpha)$ is a rotation around $z=0$ by the angle $\alpha$.

The quantities $r_1$, $r_2$, $l$, and $\varphi_1$, $\varphi_2$,
$\theta=\pi-\theta^*$ are the sides and angles of a hyperbolic triangle as in
figure~\ref{fig:phiRhoTriangle} if and only if
\begin{equation*}
  T(r_1)R(\theta^*)T(r_2) = R(\varphi_1)T(l)R(\varphi_2).
\end{equation*}
With $\rho_j=\log\tanh(r_j/2)$, the left hand side of this equation equals
\begin{equation*}
  \begin{split}
    \LHS &=
    \begin{bmatrix}1&e^{\rho_1}\\e^{\rho_1}&1\end{bmatrix}
    \begin{bmatrix}e^{i\theta^*}&0\\0&1\end{bmatrix}
    \begin{bmatrix}1&e^{\rho_2}\\e^{\rho_2}&1\end{bmatrix} \\
    &=
    \begin{bmatrix}
      e^{i\theta^*}+e^{\rho_1+\rho_2} & e^{\rho_2+i\theta^*}+e^{\rho_1} \\
      e^{\rho_1+i\theta^*}+e^{\rho_2} & 1+e^{\rho_1+\rho_2+i\theta^*}
    \end{bmatrix} \\
    &=
    \begin{bmatrix}
      \frac{\displaystyle e^{i\theta^*}+e^{\rho_1+\rho_2}} {\displaystyle
        1+e^{\rho_1+\rho_2+i\theta^*}} & \frac{\displaystyle
        e^{\rho_2+i\theta^*}+e^{\rho_1}}
      {\displaystyle 1+e^{\rho_1+\rho_2+i\theta^*}}    \\
      \frac{\displaystyle e^{\rho_1+i\theta^*}+e^{\rho_2}} {\displaystyle
        1+e^{\rho_1+\rho_2+i\theta^*}} & 1
    \end{bmatrix}.
  \end{split}
\end{equation*}
With $\lambda=\log\tanh(l/2)$, the right hand side equals
\begin{equation*}
  \begin{split}
    \RHS &=
    \begin{bmatrix}e^{i\varphi_1}&0\\0&1\end{bmatrix}
    \begin{bmatrix}1&e^{\lambda}\\e^{\lambda}&1\end{bmatrix}
    \begin{bmatrix}e^{i\varphi_2}&0\\0&1\end{bmatrix} \\
    &=
    \begin{bmatrix}
      e^{i(\varphi_1+\varphi_2)} & e^{\lambda+i\varphi_1} \\
      e^{\lambda+i\varphi_2} & 1
    \end{bmatrix}.
  \end{split}
\end{equation*}
Hence,
\begin{equation*}
  \varphi_1 = \arg\frac{\displaystyle e^{\rho_2+i\theta^*}+e^{\rho_1}}
                       {\displaystyle 1+e^{\rho_1+\rho_2+i\theta^*}},
\end{equation*}
which implies equation \eqref{eq:phiOfRhoHyp}.

\chapter{The dilogarithm function and Clausen's integral}
\label{app:dilog}

In this appendix, we collect all relevant facts about the dilogarithm
function and Clausen's integral that is relevant for this paper. A more
thorough treatment and an extensive bibliography are contained in Lewin's
monograph \cite{lewin81}.

For $|z|\leq 1$, the dilogarithm function is defined by the power series
\begin{equation*}
\Li(z)=\frac{z}{1^2}+\frac{z^2}{2^2}+\frac{z^3}{3^2}+\ldots\,.
\end{equation*}
For $|z|<1$,
\begin{equation*}
-\log(1-z)=\frac{z}{1}+\frac{z^2}{2}+\frac{z^3}{3}+\ldots,
\end{equation*}
and hence
\begin{equation}\label{dilogIntegral}
\Li(z)=-\int_{0}^{z}\frac{\log(1-\zeta)}{\zeta}\,d\zeta\,.
\end{equation}
In the light of this integral representation, one sees that the dilogarithm
can be continued analytically to the complex plane cut from $1$ to $\infty$
along the positive real axis.

Take logarithms on both sides of the identity
\begin{equation*}
1-z^n = (1-z) (1-\omega z) (1-\omega^2 z) \ldots(1-\omega^{n-1} z),
\end{equation*}
where $\omega = e^{2\pi i/n}$ is the fundamental $n^{\text{th}}$ root of
unity, then divide by $z$ and integrate to obtain
\begin{equation}\label{dilogzton}
\frac{1}{n}\,\Li(z^n) = \Li(z) + \Li(\omega z) +
\Li(\omega^2 z) + \ldots \Li(\omega^{n-1} z)
\end{equation}
for $|z|\leq 1$. Both sides of the equation can be continued analytically to
the complex plane with radial cuts outward from the
$n^{\text{th}}$ roots of unity to infinity.

Clausen's integral $\Cl(x)$ can be defined by the imaginary part of the
dilogarithm on the unit circle:
\begin{equation*}
\begin{split}
  \Cl(x)&=\im \Li(e^{ix}) \\
  &=\frac{1}{2i}\big(\Li(e^{ix})-\Li(e^{-ix})\big).
\end{split}
\end{equation*}
(The name ``Clausen's {\em integral}\thinspace" comes from the integral
representation~\eqref{eq:clausen_int}.)  We consider Clausen's integral as a
real valued function of a real variable. It is $2\pi$-periodic and odd.  The
power series representation of the dilogarithm yields
the Fourier series representation for Clausen's integral,
\begin{equation*}
 \Cl(x)=\sum_{n=1}^{\infty} \frac{\sin(nx)}{n^2}\,.
\end{equation*}
Substitute $\zeta=e^{i\xi}$ in the integral representation of the dilogarithm
(\ref{dilogIntegral}) to obtain, for $0\leq x \leq 2\pi$,
\begin{equation}\label{eq:clausen_int}
\Cl(x)=-\int_{0}^{x}\log\left(2 \sin\frac{\xi}{2}\right)\,d\xi\,.
\end{equation}
By periodicity, for all $x\in\mathds R$, 
\begin{equation}\label{eq:clausen_diff}
  \frac{d}{dx}\Cl(x)=-\log\big|2 \sin{\textstyle\frac{x}{2}}\big|.
\end{equation}
Clausen's integral is almost the same as Milnor's Lobachevski
function~\cite{milnor82}
\begin{equation*}
  \lobachevski(x)=-\int_{0}^{x}\log|2 \sin \xi|\,d\xi=\frac{1}{2}\Cl(2x)\,.
\end{equation*}
Lobachevski~\cite{lobatschefskij04} himself defined a function $L(x)$ (Lobachevski's Lobachevski
function) by
\begin{equation*}
L(x)=-\int_{0}^{x}\log(\cos \xi)\,d\xi\,,
\end{equation*}
such that 
\begin{equation*}
L(x)=\frac{1}{2}\Cl(2x-\pi)+x \log 2\,.
\end{equation*}

From equation (\ref{dilogzton}) one obtains
\begin{equation*}
\frac{1}{n}\Cl(nx)=
\sum_{k=0}^{n-1}\Cl(x+2\pi k/n)
\,.
\end{equation*}
In particular, using that $\Cl(x)$ is $2\pi$-periodic and odd, one obtains
the double-angle formula for Clausen's integral:
\begin{equation}\label{eq:clausen2x}
\frac{1}{2}\Cl(2x)=\Cl(x)-\Cl(\pi-x)\,.
\end{equation}

We will now derive a formula expressing the imaginary part of the dilogarithm
in terms of Clausen's integral not only on the unit circle, but anywhere in
the complex plane. Suppose $x$ is real and $0<\theta<2\pi$. Substitute
$\zeta=e^{\xi+i\theta}$ in (\ref{dilogIntegral}), to obtain
\begin{equation}\label{eq:ImLiIntegral}
\begin{split}
  \im \Li(e^{x+i\theta})&= \frac{1}{2i}\big(\Li(e^{x+i\theta})
  -\Li(e^{x-i\theta})\big) \\
  &=\frac{1}{2i}\int_{-\infty}^{x}
  \log\left(\frac{1-e^{\xi-i\theta}}{1-e^{\xi+i\theta}}\right)\,d\xi\,.
\end{split}
\end{equation}
Now substitute
\begin{equation*}
\eta=\frac{1}{2i}
\log\left(\frac{1-e^{\xi-i\theta}}{1-e^{\xi+i\theta}}\right)\,,
\end{equation*}
and note that inversely
\begin{equation*}
\xi=\log\left(\frac{2\sin\eta}{2\sin(\eta+\theta)}\right)\,.
\end{equation*}
Finally, one obtains
\begin{equation}\label{eq:LiCl}
\im \Li(e^{x+i\theta})
=yx + \frac{1}{2}\Cl(2y) - \frac{1}{2}\Cl(2y+2\theta)
+ \frac{1}{2}\Cl(2\theta)\,,
\end{equation}
where
\begin{equation*}
y=\frac{1}{2i}
\log\left(\frac{1-e^{x-i\theta}}{1-e^{x+i\theta}}\right)\,.
\end{equation*}
From this, we derive the formula
\begin{multline}\label{eq:LiClp}
  \im \Li(e^{x+i\theta}) + \im \Li(e^{-x+i\theta})= \\
  px + \Cl(p+\theta^*) + \Cl(-p+\theta^*) - \Cl(2\theta^*)\,,
\end{multline}
where $\theta^*=\pi-\theta$, and
\begin{equation*}
p=\frac{1}{2i}\log
\frac{(1+e^{x+i\theta^*})(1+e^{-x-i\theta^*})}
{(1+e^{x-i\theta^*})(1+e^{-x+i\theta^*})}\,.
\end{equation*}
Finally, $p$ and $x$ are related by
\begin{equation*}
\tan\Big(\frac{p}{2}\Big)=
\tanh\Big(\frac{x}{2}\Big)\tan\Big(\frac{\theta^*}{2}\Big)\,.
\end{equation*}

\chapter[Volume of an infinite triply orthogonal hyperbolic tetrahedron]{The volume of a triply orthogonal hyperbolic tetrahedron with a
  vertex at infinity}
\label{app:triply_ortho_tetra}

Milnor~\cite[pp.~19f]{milnor82} calculated the volume of a triply orthogonal
hyperbolic tetrahedron (three-dimensional hyperbolic orthoscheme) with {\em
  two}\/ vertices at infinity by a straightforward integration. By the same
computation, one derives a formula for the volume of a triply orthogonal
hyperbolic tetrahedron with one vertex at infinity. For reference, we present
this computation here.

\begin{figure}[tb]
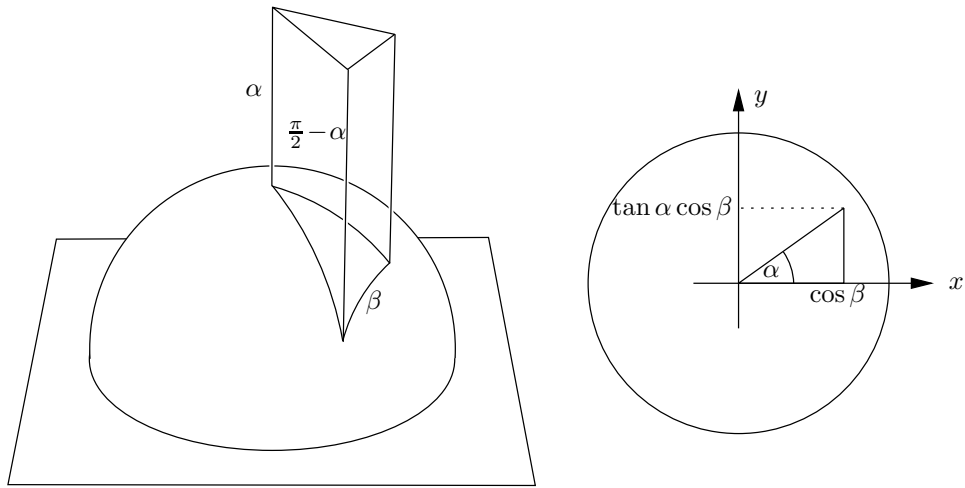
%
\centering%
\input{orthoscheme.tex}%
\hfill%
\raisebox{0.7cm}{\input{orthoscheme_domain.tex}}%
\caption{{\em Left:}~A triply orthogonal hyperbolic tetrahedron with one
  vertex at infinity (truncated at a horosphere). {\em Right:}~The domain of integration in the $xy$-plane.}%
\label{fig:orthoscheme}%
\end{figure}

Figure~\ref{fig:orthoscheme} shows a triply orthogonal tetrahedron with one
vertex at infinity in the Poincar\'e half-space model of hyperbolic space. In
the figure, the tetrahedron is truncated by a horosphere centered at the
infinite vertex. Three of the six dihedral angles are $\frac{\pi}{2}$.
Because the sum of dihedral angles is $\pi$ at an infinite vertex and larger
than $\pi$ at a finite vertex, the three remaining angles are $\alpha$,
$\frac{\pi}{2}-\alpha$, and $\beta$, as shown in the figure, where
\begin{equation*}
  0<\alpha\leq\beta<\frac{\pi}{2}.
\end{equation*}
If $\alpha=\beta$, the tetrahedron has two infinite vertices.

\begin{proposition}
  The hyperbolic volume of the triply orthogonal hyperbolic tetrahedron with
  a vertex at infinity is
  \begin{equation}
    \label{eq:vol_orthoscheme}
    V=\frac{1}{8}\Big(
    2\Cl(\pi-2\alpha)+\Cl(2\alpha-2\beta)+\Cl(2\alpha+2\beta)
    \Big).
  \end{equation}
\end{proposition}

Setting $\alpha=\beta$ and using the double-angle formula for Clausen's
integral~\eqref{eq:clausen2x}, one obtains the formula for the volume of a
triply orthogonal hyperbolic tetrahedron with {\em two}\/ vertices at
infinity and {\em characteristic angle} $\alpha$:
\begin{equation}
  \label{eq:Valpha}
  V=\frac{1}{4}\Cl(2\alpha).
\end{equation}

\begin{proof}[Proof of the proposition]
  In the Poincar\'e half-space model, hyperbolic space is
  $H^3=\{(x,y,z)\in\mathds{R}^3\,|\,z>0\}$ with metric
  $ds^2=\frac{1}{z^2}(dx^2+dy^2+dz^2)$. The volume form is therefore
  $\frac{1}{z^3}dx\,dy\,dz$. Hence,
  \begin{equation*}
    V = 
    \int_{x=0}^{\cos\beta}
    \int_{y=0}^{x\tan\alpha}
    \int_{z=\sqrt{1-x^2-y^2}}^{\infty}
    \,\frac{1}{z^3}\,dz\,dy\,dx,
  \end{equation*}
  see figure~\ref{fig:orthoscheme} {\em (right)}. The first two integrations
  are readily performed:
  \begin{equation*}
    \begin{split}
      V &= \frac{1}{2}
      \int_{x=0}^{\cos\beta}
      \int_{y=0}^{x\tan\alpha}
      \frac{1}{1-x^2-y^2}\,dy\,dx \\
      &= \frac{1}{2}
      \int_{x=0}^{\cos\beta}
      \frac{1}{\sqrt{1-x^2}}
      \artanh\bigg(\frac{x\tan\alpha}{\sqrt{1-x^2}}\bigg)\,dx.
    \end{split}
  \end{equation*}
  Substitute $x=\cos t$, and note that
  \begin{equation*}
    \begin{split}
      \artanh\Big(\frac{\tan\alpha}{\tan t}\Big) 
      &= \frac{1}{2}\log
      \frac{1+\frac{\sin\alpha}{\cos\alpha}\frac{\cos t}{\sin t}}
           {1-\frac{\sin\alpha}{\cos\alpha}\frac{\cos t}{\sin t}}\\
      &= \frac{1}{2}\log
      \frac{\cos\alpha\sin t + \sin\alpha\sin t}
           {\cos\alpha\sin t - \sin\alpha\sin t} \\
      &= \frac{1}{2}\log\frac{\sin(t+\alpha)}{\sin(t-\alpha)}.
    \end{split}
  \end{equation*}
  Using equation~\eqref{eq:clausen_int}, one obtains
  \begin{equation*}
    \begin{split}
      V &= -\frac{1}{4}\int_{t=\frac{\pi}{2}}^{\beta}
      \log\frac{\sin(t+\alpha)}{\sin(t-\alpha)} \\
      &=
      \frac{1}{8}
      \Big[\Cl(2t+2\alpha)-\Cl(2t-2\alpha)\Big]_{t=\frac{\pi}{2}}^{\beta},
    \end{split}
  \end{equation*}
  and hence equation~\eqref{eq:vol_orthoscheme}.
\end{proof}

\chapter[Topology and homology of cellular surfaces]{The combinatorial topology and homology of cellular surfaces}
\label{app:combi_top}

This appendix is concerned with a combinatorial model for finite cell
decompositions of surfaces, and with their $\Z_2$-homology. Equivalent
combinatorial models were described by Jacques~\cite{jacques70},
Tutte~\cite{tutte73}, and Jackson and Visentin~\cite{jackson_visentin01}.
They are sometimes called `winged-edge models'. We present as much of the
$\Z_2$-homology theory as needed in the proof of
lemma~\ref{lem:general_euler} in section~\ref{sec:proof_lem:general_euler}.
Of course, all of this is well known.

\section{Cellular surfaces}
\label{sec:cellular_surface}

\begin{figure}[tb]%
\centering%
\input{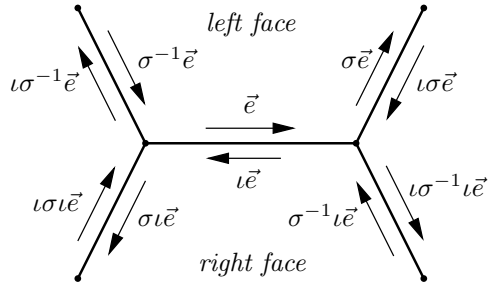}%
\caption{Elements of the group generated by the permutations $\iota$ and
  $\tau$ acting on an oriented edge $\vece$.}%
\label{fig:cellsurf_edges}
\end{figure}

Consider cell decompositions of a compact oriented 2-manifold without
boundary. There is one such cell decomposition with no edges: The sphere with
one face and one vertex. The others can be
described by two permutations of the oriented edges: The first, $\iota$,
reverses the orientation of an edge. The second, $\sigma$, maps each oriented
edge to the next oriented edge in the boundary of the face on its left side.
(See figure~\ref{fig:cellsurf_edges}.)

\begin{definition}
  An {\em oriented cellular surface}\/ $\Sigma=(\vecE,\iota,\sigma)$ is a
  finite set $\vecE$\/ together with an involution
  $\iota:\vecE\rightarrow\vecE$\/ without fixed points (that is,
  $\iota^2=\id$ and $\iota\vece\neq\vece$\/ for all $\vece\in\vecE$) and a
  permutation $\sigma:\vecE\rightarrow\vecE$.
  
  The elements $\vece\in\vecE$ are the {\em oriented edges} of $\Sigma$. The
  {\em unoriented edges}\/ are unordered pairs
  $|\vece\,|=\{\vece,\iota\vece\}$.
  
  The {\em faces}\/ of $\Sigma$ are the orbits in $\vecE$\/ of the group
  generated by $\sigma$. The {\em face on the left side} of an oriented edge
  $\vece\in\vecE$ is the orbit of $\vece$. The {\em face on the right side}
  of $\vece$ is the orbit of $\iota\vece$.
 
  The {\em vertices}\/ of $\Sigma$ are the orbits in $\vecE$\/ of the group
  generated by $\iota\sigma^{-1}$. The {\em initial vertex}\/ of an
  oriented edge $\vece\in\vecE$\/ is the orbit of $\vece$. The {\em terminal
    vertex}\/ of $\vece$\/ is the orbit of $\iota\vece$.

  An oriented cellular surface is {\em connected}\/ if the group generated by
  $\iota$ and $\sigma$ acts transitively on $\vece$.
\end{definition}

\begin{figure}[tb]%
\centering%
\input{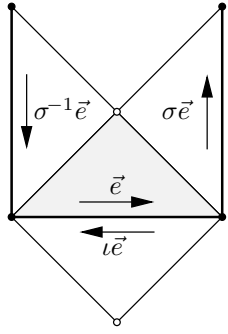}%
\caption{How the `winged edges' are glued together.}%
\label{fig:winged_edges}%
\end{figure}

Every oriented cellular surface corresponds to a cell decomposition of a
compact oriented 2-manifold without boundary. Indeed, the cell decomposition
can be constructed from the cellular surface as follows. For each oriented
edge $\vece\in\vecE$, take one triangle (a `winged edge') and glue them
together as shown in figure~\ref{fig:winged_edges}.

\begin{definition}
  Let $\Sigma_1=(\vecE_1,\iota_1,\sigma_1)$ and
  $\Sigma_2=(\vecE_2,\iota_2,\sigma_2)$ be two oriented cellular surfaces. An
  {\em orientation preserving homeomorphism} $f$\/ from $\Sigma_1$ to
  $\Sigma_2$ is a bijection $f:\vecE_1\rightarrow\vecE_2$\/ with
  \begin{equation*}
    f\circ\iota_1=\iota_2\circ f
    \quad\text{and}\quad 
    f\circ\sigma_1=\sigma_2\circ f.
  \end{equation*}
  An {\em orientation reversing homeomorphism} $f$\/ from $\Sigma_1$ to
  $\Sigma_2$ is a bijection $f:\vecE_1\rightarrow\vecE_2$ with
  \begin{equation*}
    f\circ\iota_1=\iota_2\circ f
    \quad\text{and}\quad 
    f\circ\iota_1\sigma_1^{-1}\iota_1=\sigma_2\circ f.
  \end{equation*}
\end{definition}

An orientation preserving homeomorphism $f$ maps the orbits of $\iota_1$,
$\sigma_1$, and $\iota_1\circ\sigma_1^{-1}$ to the orbits of $\iota_2$,
$\sigma_2$, and $\iota_2\circ\sigma_2^{-1}$. Hence, it induces bijections
between the non-oriented edges, the faces, and the vertices of $\Sigma_1$ and
$\Sigma_2$. The face on the left side of $\vece\in\vecE_1$ is mapped to the
face on the left side of $f(\vece)$, and the initial vertex of $\vece$
is mapped to the initial vertex of $f(\vece)$.

An orientation reversing homeomorphism $f$ also induces bijections between
the non-oriented edges, the faces, and the vertices of the two cellular
surfaces. The orbits of $\iota_1$ and $\iota_1\circ\sigma_1^{-1}$ are mapped
to the orbits of $\iota_2$ and $\iota_2\circ\sigma_2^{-1}$ (although the
latter are traversed in the opposite direction). The orbits of $\sigma_1$ are
not mapped to orbits of $\sigma_2$, but to the images of such orbits under
$\iota_2$. Nonetheless, this puts the faces in a one-to-one correspondence:
The orientation reversing homeomorphism $f$ maps the face on the left side of
$\vece\in\vecE_1$ to the face on the right side of $f(\vece)$.

Cell decompositions of non-oriented (and possibly non-orientable) 2-manifolds
are described in terms of their orientable double cover.
 
\begin{definition}
  A {\em non-oriented cellular surface}\/ $\Sigma=(\vecE,\iota,\sigma,\tau)$
  is an oriented cellular surface $\Sigma_0=(\vecE,\iota,\sigma)$ together
  with an orientation reversing homeomorphism $\tau$ from $\Sigma_0$ onto
  itself, which is an involution without fixed oriented edges. The {\em
  faces}, {\em oriented edges}, and {\em vertices} of $\Sigma$ are pairs of
  faces, oriented edges, and vertices of $\Sigma_0$, which are mapped to each
  other by $\tau$.
\end{definition}

In other words, $\tau:\vecE\rightarrow\vecE$ is another permutation with
$\tau^2=\id$,\, $\tau(\vece)\neq\vece$ for all $\vece\in\vecE$,\,
$\tau\iota=\iota\tau$, and $\tau\iota\sigma^{-1}\iota=\sigma\tau$.

\begin{definition}
  The {\em Poincar\'e dual cellular surface}\/ of an oriented cellular
  surface $\Sigma=(\vecE,\iota,\sigma)$ is the oriented cellular surface
  $\Sigma^*=(\vecE,\iota,\sigma^*)$ with $\sigma^*=\iota\sigma^{-1}$. The
  Poincar\'e dual cellular surface of a non-oriented cellular surface
  $\Sigma=(\vecE,\iota,\sigma,\tau)$ is the non-oriented cellular surface
  $\Sigma^*=(\vecE,\iota,\sigma^*,\tau^*)$ with $\sigma^*$ as in the oriented
  case and $\tau^*=\iota\tau$.
\end{definition}

There is a one-to-one correspondence between the unoriented edges of $\Sigma$
and $\Sigma^*$, between the vertices of $\Sigma$ and the faces of $\Sigma^*$,
and between the the faces of $\Sigma$ and the vertices of $\Sigma^*$. The
Poincar\'e dual of the Poincar\'e dual is essentially the same cellular
surface. The map $\vece\mapsto\iota\vece$\/ is an orientation preserving
homeomorphism from $\Sigma$ to $\Sigma^{**}$.

\section{Surfaces with boundary}
\label{sec:surfaces_w_boundary}

Combinatorial models for cell decompositions of surfaces with boundary are
less canonical. In this section we describe a few ways in which the above
definition for cellular surfaces may be adapted to deal with bounded
surfaces. We treat only the oriented case; unoriented surfaces may be
described in terms of their oriented double cover, as in the previous
section.

The most straightforward way to define oriented cellular surfaces with
boundary is to mark some faces as holes:

\begin{definition}
  An {\em oriented cellular surface with holes}\/
  $\Sigma=(\vecE,\iota,\sigma,H)$ is an oriented cellular surface
  $(\vecE,\iota,\sigma)$ together with a subset $H$\/ of the set of faces,
  such that, if an oriented edge $\vece\in\vecE$ is contained in a face of
  $H$, then $\iota\vece$ is not contained in a face of $H$.
\end{definition}

There is no Poincar\'e duality for cellular surfaces with boundary by this
definition. As a remedy, one may consider surfaces with holes and punctures:

\begin{definition}
  An {\em oriented cellular surface with holes and punctures}
  $$\Sigma=(\vecE, \iota, \sigma, H, P)$$ is an oriented cellular surface
  $(\vecE,\iota,\sigma)$ together with a subset $H$\/ of the set of faces and
  a subset $P$\/ of the set of vertices, such that, if an oriented edge
  $\vece\in\vecE$ is contained in a face of $H$ or a vertex of $P$, then
  $\iota\vece$ is not contained in a face of $H$ or a vertex in $P$.
\end{definition}

(If Poincar\'e duality is to be upheld, and adjacent holes are not
allowed, then adjacent punctures must be forbidden, too. Alternatively, one
may allow adjacent holes and punctures.)

In the computer implementation described in chapter~\ref{cha:computer}, a
slightly more exotic combinatorial model is used. It is obtained by deleting
all boundary edges. As a consequence, some faces have non-closed boundary;
and therefore, for some oriented edges, there {\em is no}\/ next edge in the
boundary of the left face.  This model is well suited for applications where
one is not interested in boundary edges at all, as in the case of circle
patterns.

\begin{definition}
  An {\em oriented cellular surface with boundary faces and boundary
  vertices} $\Sigma=(\vecE, \iota, \sigma)$ is a finite set $\vecE$ together
  with an involution $\iota:\vecE\rightarrow\vecE$ without fixed points and a
  map $\sigma:\vecE\rightarrow\vecE\,\cup\{\emptyset\}$ which satisfies the
  following condition: 
  \begin{equation*}
    \forall\;\vece_1,\vece_2\in\vecE:\quad
    \sigma(\vece_1)=\sigma(\vece_2)\neq\emptyset
    \quad\Longrightarrow\quad
    \vece_1=\vece_2.
  \end{equation*}
\end{definition}

Although $\sigma$ is not invertible, there is a `pseudo-inverse'
$\sigma^{-1}:\vecE\rightarrow\vecE\,\cup\{\emptyset\}$ which satisfies
$\sigma^{-1}(\sigma(\vece\,))=\vece$ for all $\vece\in\vecE$ with
$\sigma(\vece\,)\neq\emptyset$, and $\sigma(\sigma^{-1}(\vece\,))=\vece$ for
all $\vece\in\vecE$ with $\sigma^{-1}(\vece\,)\neq\emptyset$. Faces and
vertices are defined as `orbits' of $\sigma$ and $\iota\sigma^{-1}$ as in the
case of ordinary oriented cellular surfaces. Every edge has a left and a
right face and an initial and a terminal vertex, and the Poincar\'e dual
cellular surface is well defined.

\section{$\Z_2$-Homology}
\label{sec:homology}

Let $\Sigma$ be an oriented cellular surface, or, equivalently, a cell
decomposition of a compact $2$-manifold without boundary. Let $F$, $E$, and
$V$ be the sets of faces, non-oriented edges, and vertices. 

\subsection{Homology groups}
\label{sec:homology_groups}

The {\em chain-spaces}\/ $C_0(\Sigma)$, $C_1(\Sigma)$, and $C_2(\Sigma)$ are
the $\Z_2$ vector spaces generated by $V$, $E$, and $F$; that is, they are
the vector spaces of formal sums of vertices, edges, and faces with
coefficients in the field $\Z_2=\Z/2\Z$. The {\em boundary operators}
$\partial:C_2(\Sigma)\rightarrow C_1(\Sigma)$ and
$\partial:C_1(\Sigma)\rightarrow C_0(\Sigma)$ are the $\Z_2$-linear maps which
are defined as follows. A face $f\in F$\/ of a cellular surface was defined
as an orbit $f=\{\vece_1,\ldots,\vece_n\}$ of oriented edges. The boundary of
$f$ is $\partial f=|\vece_1|+\ldots+|\vece_n|$, where $|\vece\,|$ is the
non-oriented edge corresponding to $\vece$. For an edge $e\in E$ with end
vertices $v_1,v_2\in V$, $\partial e=v_1+v_2$.

The {\em chain complex}\/ $C_*(\Sigma)$ is the following sequence of chain
spaces and linear maps:
\begin{equation*}
  0\overset{\partial}{\longrightarrow} 
  C_2(\Sigma)\overset{\partial}{\longrightarrow} C_1(\Sigma)
  \overset{\partial}{\longrightarrow} C_0(\Sigma)
  \overset{\partial}{\longrightarrow} 0
\end{equation*}
The {\em cycle spaces}\/ are the kernels $Z_n(\Sigma)=\{c\in
C_n(\Sigma)\,|\,\partial c=0\}$ for $n\in\{0,1,2\}$. The {\em boundary spaces}
are the images $B_2(\Sigma)=0$ and $B_n(\Sigma)=\partial C_{n+1}(\Sigma)$ for
$n\in\{0,1\}$.  The {\em homology groups}\/ are the quotient spaces
$H_n(\Sigma)=Z_n(\Sigma)/B_n(\Sigma)$ and their dimensions are the {\em Betti
numbers} $h_n(\Sigma)=\dim_{\Z_2}H_n(\Sigma)$.

The second Betti number $h_2(\Sigma)$ is the number of equivalence classes of
faces under the equivalence relation
\begin{equation*}
  f_1\sim f_2 \quad\Longleftrightarrow\quad 
  \parbox{0.35\textwidth}{Every closed $2$-chain which contains $f_1$ also contains $f_2$.}
\end{equation*}
The zeroth Betti number $h_0(\Sigma)$ is the number of equivalence classes of
vertices under the equivalence relation
\begin{equation*}
  v_1\sim v_2 \quad\Longleftrightarrow\quad 
  \parbox{0.35\textwidth}{The vertices $v_1$ and $v_2$ can be joined by a path.}
\end{equation*}
Using subdivision arguments one can show that $h_2(\Sigma)=h_0(\Sigma)$. Of
course, this is just the number of connected components of $\Sigma$.

\subsection{Poincar\'e duality}
\label{sec:poincare_duality}

The Poincar\'e dual cellular surface $\Sigma^*$ has face set $F^*=V$, edge
set $E^*=E$, and vertex set $V^*=F$. The chain spaces are therefore
$C_n(\Sigma^*)=C_{2-n}(\Sigma)$ for $n\in\{0,1,2\}$. Denote the boundary
operators of the chain complex $C_*(\Sigma^*)$ by $\partial^*$. Let the
chain spaces be equipped with scalar products 
$\langle\,\cdot\,,\,\cdot\,\rangle$
such that $F$, $E$, and $V$ are orthonormal bases.

\begin{lemma}[Poincar\'e duality]
  The boundary operators $\partial^*$ of the chain complex $C_*(\Sigma^*)$
  are the adjoint operators of the boundary operators $\partial$ of the
  chain complex $C_*(\Sigma)$ with respect to the scalar product
  $\langle\,\cdot\,,\,\cdot\,\rangle$. 

  This means, for $f\in F=V^*$, $e\in
  E=E^*$, and $v\in V=F^*$,
  \begin{equation*}
    \langle\partial f,e\rangle=\langle f,\partial^* e\rangle
    \quad\text{and}\quad
    \langle\partial e,v\rangle=\langle e,\partial^* v\rangle.
  \end{equation*}
\end{lemma}

\begin{proposition}[Poincar\'e duality for cellular surfaces]
  The Betti numbers of $\Sigma$ and $\Sigma^*$ satisfy
  \begin{equation*}
    h_n(\Sigma^*)=h_{2-n}(\Sigma).
  \end{equation*}
\end{proposition}

\begin{proof}
  The dimensions of $H_n(\Sigma^*)$ and $H_{2-n}(\Sigma)$ are equal because
  they are dual vector spaces. This follows from the following facts from
  linear algebra. 

  Let $V$ be a vector space, and let $U$ and $W$ be subspaces such that
  $W\subset U\subset V$. Then the dual space of the
  quotient $U/W$ is 
  \begin{equation*}
    (U/W)^*=W^\perp/U^\perp,
  \end{equation*}
  where $W^\perp\subset V^*$ and $U^\perp\subset V^*$ are the spaces of linear
  functionals on $V$ which vanish on $W$ and $U$, respectively.
  
  Let $f:V\longrightarrow W$ be a linear map between two vector spaces, and
  let $f^*:W^*\longrightarrow V^*$ be the adjoint map. Then the kernels and
  images of $f$ and $f^*$ are related by
  \begin{equation*}
    \image f^* = (\kernel f)^\perp,\qquad \kernel f^*=(\image f)^\perp.
  \end{equation*}

  Hence,
  \begin{multline*}
    H_n(\Sigma)^*
    =\big(Z_n(\Sigma)/B_n(\Sigma)\big)^*
    =B_n(\Sigma)^\perp/Z_n(\Sigma)^\perp =\\
    =Z_{2-n}(\Sigma^*)/B_{2-n}(\Sigma^*)=H_{2-n}(\Sigma^*).
  \end{multline*}
\end{proof}

\subsection{Relative homology}
\label{sec:relative_homology}

Define a {\em subcomplex}\/ $\Upsilon$ of a cellular surface $\Sigma$ as a
triple $\Upsilon=(F_\Upsilon, E_\Upsilon, V_\Upsilon)$ of subsets
$F_\Upsilon\subset F$, $E_\Upsilon\subset E$, and $V_\Upsilon\subset V$, such
that the following holds: If a face in $F_\Upsilon$ contains an oriented edge
$\vece$, then the corresponding unoriented edge $|\vece\,|\in E_\Upsilon$;
and the end-vertices of any edge in $E_\Upsilon$ are contained in
$V_\Upsilon$.

The {\em relative chain complex}\/ $C_*(\Sigma,\Upsilon)$ consists of the
quotient spaces $C_n(\Sigma,\Upsilon)=C_n(\Sigma)/C_n(\Upsilon)$ and the
boundary operators induced by the boundary operators of $C_*(\Sigma)$.  The
corresponding homology groups $H_n(\Sigma,\Upsilon)$ are the {\em relative
  homology groups}\/ of the pair $(\Sigma,\Upsilon)$.

The composition of the inclusion map $C_n(\Upsilon)\rightarrow C_n(\Sigma)$,
followed by the natural projection $C_n(\Sigma)\rightarrow
C_n(\Sigma,\Upsilon)$ is the zero map. Hence, the short sequence
\begin{equation*}
  0\longrightarrow C_*(\Upsilon)\longrightarrow C_*(\Sigma)
  \longrightarrow C_*(\Sigma,\Upsilon)\longrightarrow 0
\end{equation*}
is exact. This means that the columns of the commutative
diagram shown in figure~\ref{fig:relative_homology} are exact.
\begin{figure}[tb]%
\centering%
\begin{equation*}
  \begin{CD}
                @.      0      @.           0 @. 0\\
    @. @VVV @VVV @VVV\\
    0 @>>>  C_2(\Upsilon)     @>\partial>> C_1(\Upsilon)        @>\partial>> 
    C_0(\Upsilon) @>>> 0 \\
    @. @VVV @VVV @VVV\\
    0 @>>> C_2(\Sigma)        @>\partial>> C_1(\Sigma)        @>\partial>>
    C_0(\Sigma) @>>> 0 \\
    @. @VVV @VVV @VVV \\
    0 @>>> C_2(\Sigma,\Upsilon) @>\partial>> C_1(\Sigma,\Upsilon) @>\partial>> 
    C_0(\Sigma,\Upsilon) @>>> 0\\
    @. @VVV @VVV @VVV \\
    @. 0 @. 0 @. 0    
  \end{CD}
\end{equation*}%
\caption{The short exact sequence $0\rightarrow C_*(\Upsilon)\rightarrow C_*(\Sigma)\rightarrow C_*(\Sigma,\Upsilon)\rightarrow 0$ of chain complexes.}%
\label{fig:relative_homology}%
\end{figure}
By the zigzag-lemma, there is a long exact sequence of homology groups
\begin{multline*}
  0\longrightarrow 
  H_2(\Upsilon)\longrightarrow
  H_2(\Sigma)\longrightarrow
  H_2(\Sigma,\Upsilon)\longrightarrow\\
  \longrightarrow H_1(\Upsilon)\longrightarrow 
  H_1(\Sigma)\longrightarrow 
  H_1(\Sigma,\Upsilon)\longrightarrow \\
  \longrightarrow H_0(\Upsilon)\longrightarrow H_0(\Sigma)\longrightarrow
  H_0(\Sigma,\Upsilon)\longrightarrow 0.
\end{multline*}

\subsection{Lefschetz duality}
\label{sec:lefschetz_duality}

Let $\Upsilon=(F_\Upsilon, E_\Upsilon, V_\Upsilon)$ be a subcomplex of the
cellular surface  $\Sigma$ as in the previous section. Then the complement of
$\Upsilon$ in $\Sigma^*$,
\begin{equation*}
  \Sigma^*\setminus\Upsilon:=(F^*\setminus V_\Upsilon, E^*\setminus
  E_\Upsilon,V^*\setminus F_\Upsilon),
\end{equation*}
is a subcomplex of $\Sigma^*$.

\begin{proposition}[Lefschetz duality for cellular surfaces]
  The Betti numbers of the subcomplex $\Sigma^*\setminus\Upsilon$ and the
  pair $(\Sigma,\Upsilon)$ satisfy
  \begin{equation*}
    h_n(\Sigma^*\setminus\Upsilon)=h_{2-n}(\Sigma,\Upsilon).
  \end{equation*}
\end{proposition}

\begin{proof}
  As in the case of Poincar\'e duality, the dimensions of the homology groups
  $H_n(\Sigma^*\setminus\Upsilon)$ and $H_{2-n}(\Sigma,\Upsilon)$ are equal
  because they are dual vector spaces. This is proved in the same way.
\end{proof}

\subsection{Proof of lemma~\ref{lem:general_euler}}
\label{sec:proof_lem:general_euler}

Consider the long exact homology sequence of relative homology for a
$1$-dimensional subcomplex $\Gamma$, that is, a cellularly embedded graph.
Let $E_{\Gamma}\subset E$ and $V_{\Gamma}\subset V$ be the edge and vertex
sets of $\Gamma$.  Since $C_2(\Gamma)$ is zero, the short exact sequence
\begin{equation*}
  0\longrightarrow C_*(\Gamma)\longrightarrow C_*(\Sigma)
  \longrightarrow C_*(\Sigma,\Gamma)\longrightarrow 0
\end{equation*}
of chain-complexes is equivalent to the commutative diagram shown in
figure~\ref{fig:relative_homology_gamma}, in which the columns are exact.
\begin{figure}[tb]%
\centering%
\begin{equation*}
  \begin{CD}
                @.            @.           0 @. 0\\
    @. @. @VVV @VVV\\
                @.   0        @>\partial>> C_1(\Gamma)        @>\partial>> 
    C_0(\Gamma) @>>> 0 \\
    @. @VVV @VVV @VVV\\
    0 @>>> C_2(\Sigma)        @>\partial>> C_1(\Sigma)        @>\partial>>
    C_0(\Sigma) @>>> 0 \\
    @. @VVV @VVV @VVV \\
    0 @>>> C_2(\Sigma,\Gamma) @>\partial>> C_1(\Sigma,\Gamma) @>\partial>> 
    C_0(\Sigma,\Gamma) @>>> 0\\
    @. @VVV @VVV @VVV \\
    @. 0 @. 0 @. 0
  \end{CD}
\end{equation*}%
\caption{The short exact sequence $0\rightarrow C_*(\Gamma)\rightarrow 
  C_*(\Sigma)\rightarrow C_*(\Sigma,\Gamma)\rightarrow 0$.}
\label{fig:relative_homology_gamma}
\end{figure}
This gives rise to the following long exact homology sequence.
\begin{multline*}
  0\longrightarrow H_2(\Sigma)\longrightarrow
  H_2(\Sigma,\Gamma)\longrightarrow
  H_1(\Gamma)\longrightarrow 
  H_1(\Sigma)\longrightarrow 
  H_1(\Sigma,\Gamma)\longrightarrow \\
  \longrightarrow H_0(\Gamma)\longrightarrow H_0(\Sigma)\longrightarrow
  H_0(\Sigma,\Gamma)\longrightarrow 0
\end{multline*}
Since the sequence is exact, the corresponding dimensions satisfy
\begin{equation*}
  h_2(\Sigma)-h_2(\Sigma,\Gamma)+h_1(\Gamma)-h_1(\Sigma)+h_1(\Sigma,\Gamma) -
  h_0(\Gamma)+h_0(\Sigma)-h_0(\Sigma,\Gamma) = 0.    
\end{equation*}
Now assume that $\Sigma$ is connected and $\Gamma$ is non-empty, that is,
$V_{\Gamma}\neq\emptyset$. Then $h_0(\Sigma,\Gamma)=0$, because from every vertex of
$\Sigma$ there is a path to a vertex of the graph $\Gamma$. Also,
\begin{equation*}
  h_1(\Gamma)-h_0(\Gamma)=|E_{\Gamma}|-|V_{\Gamma}|
\end{equation*}
and
\begin{equation*}
  h_2(\Sigma)-h_1(\Sigma)+h_0(\Sigma)=|F|-|E|+|V|,
\end{equation*}
and hence,
\begin{equation*}
  h_2(\Sigma,\Gamma)-|E_{\Gamma}|+|V_{\Gamma}| = |F|-|E|+|V| + h_1(\Sigma,\Gamma).
\end{equation*}
Now, $h_2(\Sigma,\Gamma)$ is the number of regions into which the graph
$\Gamma$ separates the cellular surface $\Sigma$. By Lefschetz-duality
$h_1(\Sigma,\Gamma)=h_1(\Sigma^*\setminus\Gamma)$. Suppose $\Gamma$ separates
$\Sigma^*$ into the components $\Pi_1,\ldots,\Pi_n$; then
$h_1(\Sigma^*\setminus\Gamma)=\sum h_1(\Pi_j)$.

This completes the proof of lemma~\ref{lem:general_euler}.

\backmatter{}

\refstepcounter{chapter}
\let\chaptername

\bibliographystyle{plain} 
\bibliography{diss}

\end{document}